\documentclass[a4paper,10pt]{article}
\usepackage{amssymb}
\usepackage{cite}
\usepackage{color}
\usepackage{amsmath}
\usepackage{color,enumitem,graphicx}
\usepackage{subfigure}
 \usepackage{placeins}
\usepackage[export]{adjustbox}
\allowdisplaybreaks

\newcommand{\cA}{{\cal A}}

\newcommand{\cK}{{\cal K}}

\newcommand{\cJ}{{\cal J}}

\newcommand{\R}{\mathbb{R}}
\newcommand{\Z}{\mathbb{Z}}
\newcommand{\N}{\mathbb{N}}

\newcommand{\cI}{{\cal I}}
\newcommand{\cP}{{\cal P}}

\newcommand{\cuad}{{\sqcap\kern-.68em\sqcup}}

\newcommand{\norm}[1]{\|#1\|}

\newcommand{\bQ}{\mathbb{Q}}
\newcommand{\bK}{\mathbb{K}}

\newcommand{\bA}{\mathbb{A}}

 \newcommand{\cM}{{\cal M}}
\newcommand{\bF}{\mathbb{F}}
\newcommand{\bT}{\mathbb{T}}

\newtheorem{theorem}{Theorem}[section]
\newtheorem{proposition}{Proposition}[section]

\newtheorem{lemma}{Lemma}[section]
\newtheorem{corollary}{Corollary}[section]
\newtheorem{remark}{Remark}[section]
\newcommand{\bremark}{\begin{remark} \em}
\newcommand{\eremark}{\end{remark} }

\begin{document}

\begin{center}{\bf  \large On positive  solutions of Lane-Emden equations
\\[2mm]
 on the integer lattice graphs    \\[2mm]
 }\bigskip

{\small
  Huyuan Chen\qquad Bobo Hua\qquad Feng Zhou  }  \\[10pt]

 \bigskip

\begin{abstract}
In this paper, we study the existence and nonexistence of positive solutions to the Lane–Emden equation  
$$ -\Delta u = Q(x)\,|u|^{p-2}u $$  
on three distinct domains: the $d$-dimensional integer lattice graph $\mathbb{Z}^d$, the discrete half-space, and a discrete quadrant-type domain—each equipped with the zero Dirichlet boundary condition in the latter two cases. Here, $d \geq 3$, $p > 1$, and $Q(x)$ is a Hardy-type positive potential satisfying $Q(x) \sim (1+|x|)^{-\alpha}$ for some $\alpha \in [0,\infty)$.  

We precisely characterize the Sobolev supercritical regime—the set of parameter pairs $(\alpha,p)$ lying strictly above the Sobolev critical curve and super linear —where existence of positive solutions is established via variational methods. Conversely, in the Serrin subcritical regime—i.e., the region strictly below the Serrin critical curve—we prove nonexistence by a refined iterative analysis of the asymptotic decay of solutions at infinity when $\alpha\in[0,2)$, culminating in a contradiction. Notably, between the Sobolev and Serrin critical curves lies a transitional region whose classification remains open across all three domains: $\mathbb{Z}^d$, the half-space, and the quadrant.  For $\alpha>2$, the existence holds for the whole range $p>1$. 
\end{abstract}

  \end{center}
  {\small
    \tableofcontents }\vspace{5mm}
   \smallskip

   \noindent {\small {\bf AMS Subject Classifications}:     35R02;  35J20; 35A01.
  \smallskip

\vspace{2mm}

\setcounter{equation}{0}
\section{Introduction}

Let $\Z^d$ be the \textit{d-dimensional integer lattice graph}   consisting of the set of vertices $\mathbb{Z}^d$, the edge weight be defined by
\[\omega:\Z^d\times \Z^d\to [0,+\infty),\]
\begin{equation*}
    \omega_{xy}= \begin{cases}
1\ \,\, \text{ if } \displaystyle |x-y|_{_Q}:=  \sum^d_{k=1}|x_k-y_k|=1,\\
0\ \,\, \text{ otherwise}
\end{cases}
\end{equation*}
$$x\sim y\quad {\rm if}\quad \omega_{xy}=1,$$
so that the Laplace be defined as
 $$\Delta_{\Z^d} u(x)=\sum_{y\sim x}\big(u(y)-u(x)\big)\quad\ \text{ for all }x\in\mathbb{Z}^d.$$
Our first purpose in this article is to prove the nonexistence of  solution of semilinear   elliptic  equation
in the whole integer lattice space
 \begin{equation}\label{eq 1.1}
 -\Delta_{\Z^d}  u=Q |u|^{p-2}u  \quad
    {\rm in}\ \  \Omega\subset \Z^d,
\end{equation}
subject to the zero Dirichlet boundary condition $u=0$ on $\partial \Omega$ if $\Omega\not=\Z^d$, 
where  $d\geq 3$, $p>1$ and $Q\in C(\Z^d)$ is a nonnegative Hardy-type potential.  \smallskip

%  Assume that
% \begin{equation}\label{exp 1}
%Q\in L^{q_0,\infty}(\Z^d),\quad q_0\in[1,+\infty],
% \end{equation}
% where $L^{q_0,\infty}(\Z^d)$ is the weak $L^{q_0}$ space.

 The Lane-Emden equation is a classical model of semilinear elliptic differential equations that arises in astrophysics for describing the structure of a self-gravitating, spherically symmetric polytropic fluid in hydrostatic equilibrium. The standard form of the equation is given by
\begin{equation}\label{eq 1.1-ou}
 -\Delta_{\mathbb{R}^d} u=|u|^{p-2}u \quad
    {\rm in}\ \  \mathbb{R}^d,
\end{equation}
where $p > 1$ and
$$\Delta_{\mathbb{R}^d} u(x) = \sum^d_{i=1} \partial_{ii} u(x).$$
When $p \in \left(1, \frac{2d}{d-2}\right)$, Eq. (\ref{eq 1.1-ou}) admits no positive solutions due to the Pohozaev identity. When $p = \frac{2d}{d-2}$, Eq. (\ref{eq 1.1-ou}) has exactly the following family of solutions:
$$
    u_{\lambda,\bar x}(x) = c_d \frac{\lambda^{\frac{d-2}{2}}}{(\lambda^2 + |x - \bar x|^2)^{\frac{d-2}{2}}}, \;\; x \in \mathbb{R}^d,
$$
for any $\bar x \in \mathbb{R}^d$ and $\lambda > 0$, which are the well-known Aubin-Talenti bubble solutions. When $p > \frac{2d}{d-2}$, Eq. (\ref{eq 1.1-ou}) admits infinitely many positive solutions by variational method or
the shooting method and phase-plane analysis, see \cite{St,GN,JK},  book \cite{struwe} and the references therein.

When a potential term is introduced, Eq. (\ref{eq 1.1-ou}) can be generalized as
\begin{equation}\label{eq 1.1-Eu}
 -\Delta_{\mathbb{R}^d} u = Q |u|^{p-2}u \quad
    {\rm in}\ \  \mathbb{R}^d,
\end{equation}
which may be regarded as a modified version of the classical Lane-Emden equation. This generalization is frequently used in mathematical physics and astrophysics to describe density distributions under various gravitational or thermodynamic conditions. When $p = \frac{2d}{d-2}$, Ni \cite{N} established a connection between this equation and conformal geometry, where $Q$ represents the scalar curvature of a given Riemannian manifold. He proved that Eq. (\ref{eq 1.1-Eu}) has no positive solutions if
$$
\bar Q(r) \geq C r^{\alpha_0}
$$
for some $\alpha_0 > 2$, where
$$
\bar Q(r) = \left(\frac{1}{\omega_d r^{d-1}} \int_{|x|=r} \frac{d\omega(x)}{Q(x)^{\frac{d-2}{4}}} \right)^{-\frac{4}{d-2}}.
$$
When $Q$ is nonnegative, radially symmetric, and non-increasing, Eq. (\ref{eq 1.1-Eu}) possesses infinitely many positive solutions. Bianchi et al. \cite{B} demonstrated the existence of a positive radial solution that asymptotically behaves like the standard Aubin-Talenti bubble at infinity, assuming that $Q$ is radially symmetric, decreasing, and satisfies $Q(r) \to Q_\infty > 0$ as $r \to \infty$. Cao and Peng \cite{CP} further established the existence of a positive radial solution that decays polynomially at infinity under the same assumptions but with $Q_\infty = 0$.

\smallskip

In the lattice graph, the Lane-Emden type equation can be expressed as
\begin{equation}\label{eq 1.1-zd}
-\Delta_{\mathbb{Z}^d} u = |u|^{p-2}u \quad \text{in } \mathbb{Z}^d.
\end{equation}
  Gu-Huang-Sun \cite{GHS} established that there are no positive solutions to (\ref{eq 1.1-zd}) when $d \geq 3$ and $p \leq \frac{d}{d-2} + 1$. On the other hand, Hua-Li \cite{HL} proved that (\ref{eq 1.1-zd}) admits a positive solution when $p > \frac{2d}{d-2}$. Moreover, \cite{GHS} pointed out that the question of nonexistence of positive solutions in the range $\frac{d}{d-2} + 1 < p \leq \frac{2d}{d-2}$ remains open.

    In a general graph $(G,E)$,   elliptic equations on graphs attracts  more and more attention  recently.
 Particularly, semilinear elliptic problem on graphs
   $$\Delta u+ f(x,u)=0\quad {\rm in}\ \, G $$
  has been studied in \cite{GS,GL,GHS,HLW,HLW1} for the existence of solutions,
  in \cite{GHJ,BMP} for the Liouville properties and books \cite{HLY,Grigbook}.

To investigate the existence of positive solutions to (\ref{eq 1.1}), we make the following notations:
\\[1mm]
 Set
  \begin{align}\label{exp alpha}
  \alpha:=\sup\big\{\tilde \alpha \in\R:\,  \limsup_{|x|\to+\infty} Q(x)|x|^{\tilde \alpha} <+\infty   \big\}\in(-\infty,+\infty],
  \end{align}
  \begin{align}\label{exp 2star}
  2^*_{\alpha}:=\frac{2(d-\alpha)}{d-2} .
  \end{align}
Obviously we have
\begin{align}\label{exp 2compar}
2<1+\frac{d-\alpha}{d-2}<2^*_{\alpha} \quad {\rm for} \;\; \alpha<2,  \;\; {\rm and} \;\;   2^*_{\alpha} < 1+\frac{d-\alpha}{d-2} <2\;\; {\rm  for}
\;\; \alpha>2.
  \end{align}
Then we impose the following assumptions:
\begin{itemize}
\item[ $(\bA_1)$] Let
$$\alpha\in[0,+\infty], \qquad p\in(2,+\infty)\cap  \big( 2^*_{\alpha},+\infty\big)   $$
and if $\alpha=0$, we assume more that
  $$\limsup_{|x|\to+\infty} Q(x) <+\infty. $$

 \item[$(\bA_2)$] Let
  \begin{align}\label{ass-3.1}
\lim_{|x|\to+\infty} Q(x)|x|^{\alpha}=0   \quad \text{for  $\alpha\in[0,+\infty)$}
\end{align}
 and
 $$p\in (2,+\infty)\cap [ 2^*_{\alpha}, +\infty\big).    $$

 %\item[$(A1)$] Let    $Q\in L^{q_0, \infty}(\Z^d)$ with $q_0\in[1,+\infty]$ and
% $$p>\max\Big\{2,\,  \frac{2d}{d-2} \frac{q_0-1}{q_0}\Big\};   $$

\end{itemize}

Now we can state our first results on the whole space $\Z^d$.
\begin{theorem}\label{teo 1}
$(i)$  Assume  that $d\geq 3$,  either  $(\bA_1)$ or $(\bA_2)$ holds.
  Then equation (\ref{eq 1.1})
 has at least one nontrivial positive solution $u\in L^{p}(\Z^d,Qdx)$.
Furthermore, if $Q\geq C>0$ for some $C>0$,  then
  $$\lim_{|x|\to+\infty} u(x) = 0 .  $$

$(ii)$ When $\alpha>2$ and
$$ \limsup_{|x|\to+\infty} Q(x)|x|^{\alpha}<+\infty,   $$
 then for $p\in (1,2)$,     Eq. (\ref{eq 1.1})
has  a  positive solution. \smallskip

$(iii)$ When $\alpha < 2 $ and
$$ \liminf_{|x|\to+\infty} Q(x)|x|^{\alpha}>0, $$

if  $p\in(1,1+ \frac{d-\alpha}{d-2} ]$,
%or $p=1+ \frac{d-\alpha}{d-2} >2$,
 then  Eq. (\ref{eq 1.1})
has  no positive solution.

    \end{theorem}

 \begin{remark}
 Set  $Q(x)=(1+|x|)^{-\alpha}$ in $\Z^d$, we have the following conclusions: \\
  $(a)$ In $(ii)$ when $\alpha > 2$, and for all $p \in (1,2)$, a positive solution is derived by the method of super and sub-solutions.  When $\alpha <2$ and   $(A2)$ hold, then for $p \in (2^*_{\alpha}, +\infty)$, there exists at least one nontrivial solution, which has been obtained by variational method. \\
  %$1+\frac{d-\alpha}{d-2}$ is the Serrin Exponent.  Eq.(\ref{eq 1.1}) has a unique positive solution if $1<1+\frac{d-\alpha}{d-2}<2$, does no positive solution %if $1+\frac{d-\alpha}{d-2}>2$.
  %Moreover,
  % $2<1+\frac{d-\alpha}{d-2}<2^*_{\alpha}$ for $\alpha<2$ and   $1+\frac{d-\alpha}{d-2}> 2^*_{\alpha}$ for $\alpha>2$.\\
% $(b)$ When $p>2^*_{\alpha}$, the existence of solution is derived by variational method. \\
% $(c)$  When $1+\frac{d-\alpha}{d-2}<p<2$ for $\alpha>2$, the  unique positive solution is derived by the method of super and sub solutions.
 $(b)$  In general case, it is open for the existence of solutions of (\ref{eq 1.1}) in the  zone of $(\alpha,p)$:
   $$\Big\{(\alpha,p)\in[0,+\infty)\times (1,+\infty):\ \alpha\in (0,2), \  1+\frac{d-\alpha}{d-2}< p \leq 2^*_{\alpha}
   %\leq 1+\frac{d+2-2\alpha}{d-2}
   \Big\}. $$
 $(c)$  We emphasize that when $\alpha=0$,   Eq.(\ref{eq 1.1}) admits positive solution for the critical case $p=\frac{2d}{d-2}$ along with the assumption that
   $\lim_{|x|\to+\infty}Q(x)=0, $
and for all supercritical cases $p >  \frac{2d}{d-2}$
provided the potential $Q$ satisfies $\limsup_{|x|\to+\infty}Q(x) < \infty $. \\
$(d)$ When $\alpha=+\infty$, we let $2^*_{\alpha}=-\infty$.  For instance, $Q(x)\leq c_1 e^{-c_2(1+|x|)}$ for some $c_1,c_2>0$ or $Q$ is compact supported.\smallskip

\end{remark}

     Our second purpose in this article is to prove the existence of  solution of semilinear   elliptic  equation
in the half integer lattice graph
 \begin{equation}\label{eq 1.1-half}
 \left\{%\arraycolsep=1pt
\begin{array}{lll}
 -\Delta_{\Z^d}  u=Q |u|^{p-2}u  \quad
 &{\rm in}\ \  \Z^d_+ , \\[2mm]
 \phantom{ \quad\quad \   }
 \displaystyle u=0    &{\rm on}\ \,  \partial\Z^d_+,
 \end{array}
 \right.
\end{equation}
where  $d\geq 2$, $p>2$,  $Q\in C(\Z^d_+)$ is nonnegative, nontrivial and
$$\Z^d_+=\{(x_1,x')\in\Z^d:\, x_1>0\}. $$

 % Assume that
 %\begin{equation}\label{exp 1}
%Q\in L^{q_0,\infty}(\Z^d),\quad q_0\in[1,+\infty],
% \end{equation}
 % where $L^{q_0,\infty}(\Z^d)$ is the weak $L^{q_0}$ space.

\begin{theorem}\label{teo 1-hf}
 $(i)$ Assume  that $d\geq 3$,  either  $(\bA_1)$ or $(\bA_2)$ holds.
  Then   equation (\ref{eq 1.1-half})
 has at least one nontrivial positive solution $u\in L^{p}(\Z^d_+,Qdx)$.
Furthermore, if $Q\geq C>0$ for some $C>0$,  then
  $$\lim_{x\in\Z^d_+,  |x|\to+\infty} u(x)  = 0 .  $$

  $(ii)$ When $\alpha > 2$   and
$$ \limsup_{|x|\to+\infty} Q(x)|x|^{\alpha}<+\infty,   $$
 then for $p\in(1,2)$,    Eq. (\ref{eq 1.1-half})
has  a positive solution. \smallskip

$(iii)$ When $\alpha < 2$ and
$$ \liminf_{|x|\to+\infty} Q(x)|x|^{\alpha}>0,    $$

if  $p\in(1,1+ \frac{d+1-\alpha}{d-1} )$,
%or $p=1+ \frac{d-\alpha}{d-1} >2$,
 then   Eq. (\ref{eq 1.1-half})
has  no positive solution.
       \end{theorem}

     Our  aim of this article is to prove the existence of a solution of  the Lame-Emden equation in quadrant type domain.
 \begin{equation}\label{eq 1.1-quad}
 \left\{%\arraycolsep=1pt
\begin{array}{lll}
 -\Delta_{\Z^d}  u=Q |u|^{p-2}u  \quad
 &{\rm in}\ \  \Z^d_{*} , \\[2mm]
 \phantom{ --- \! }
 \displaystyle u=0    &{\rm on}\ \,  \partial \Z^d_{*} ,
 \end{array}
 \right.
\end{equation}
where  $d\geq 3$, $p>2$,  $Q\in C(\Z^d_*)$ is nonnegative, nontrivial and
$$ \Z^d_{*} =\{(x_1,x_2,x')\in\Z^d:\, x_1,x_2>0\}. $$

    \begin{theorem}\label{teo 1-qd}
$(i)$  Assume  that $d\geq 3$,  either  $(\bA_1)$ or $(\bA_2)$ holds.
  Then   equation (\ref{eq 1.1-quad})
 has at least one nontrivial positive solution $u\in L^{p}(\Z^d_{*},Qdx)$.
  $$\lim_{x\in \Z^d_{*},  |x|\to+\infty} u(x)  = 0 .  $$

   $(ii)$ When $\alpha>2$ and
$$ \limsup_{|x|\to+\infty} Q(x)|x|^{\alpha}<+\infty,   $$
 then for $p\in (1,2)$,     Eq.(\ref{eq 1.1-quad})
has  a positive solution. \smallskip

$(iii)$ When $\alpha < 2$ and
$$ \liminf_{|x|\to+\infty} Q(x)|x|^{\alpha}>0, $$
then for  $p\in(1,1+ \frac{d+2-\alpha}{d} )$,
    Eq.(\ref{eq 1.1-quad}) has no positive solution.

    \end{theorem}

Finally, we consider the particularly linear case, i.e. $p=2$ in the whole space, half space and quadrant type domain. 

\begin{theorem}\label{teo 1-eig}
  Assume  that $d\geq 3$, $p=2$,  $\Omega=\Z^d,\ \Z^d_+$ or $\Z^d_*$,
   and nonnegative function $Q\in C(\Omega)$ satisfies  
   $$\lim_{x\in\Omega,|x|\to+\infty} Q(x)|x|^2=0. $$
 
Let
$$%\begin{equation}\label{ei-0-4}
   m_\Omega:=\big|\big\{x\in\Z^d:\, \exists y\in\Omega,\,  Q(x)>0\big\} \big|\in[1, +\infty], 
$$%  \end{equation}
then   
   \begin{itemize}
\item[(i)] 
problem 
 \begin{equation}\label{eq 1.1-eig}
 \left\{%\arraycolsep=1pt
\begin{array}{lll}
 -\Delta_{\Z^d}  u=\lambda Qu  \quad
 &{\rm in}\ \  \Omega , \\[2mm]
 \phantom{ --- \! }
 \displaystyle u=0    &{\rm on}\ \,  \Z^d\setminus  \Omega
 \end{array}
 \right.
\end{equation}
admits   a  solution $(\lambda_1,\phi_1)\in (0,+\infty)\times L^2(\Z^d)$ with $\phi_1\geq0$, 
 where
 $$\lambda_1=\Big(\sup_{\|v\|_{L^2(\Z^d)}=1}\int_{\Z^d}(Q^{\frac12}v) \Phi_{d,\Omega}\ast(Q^{\frac12}v) dx\Big)^{-1}>0. $$

Moreover,  problem \eqref{eqq-3.2-l} admits a sequence of the  eigenvalue $\lambda_{k}$ for any integer $k\in [2,m_\Omega]\cap (0,+\infty)$ if $m_0\in[2,+\infty]$, 
 with corresponding eigenfunctions $\phi_{k}$
 which can be characterized as
\begin{equation}\label{Lambda-k-s}
\lambda_{k} =\Big( \sup_{v\in \cP_{k} }(Q^{\frac12}v) \Phi_{d,\Omega}\ast(Q^{\frac12}v)  dx\Big)^{-1},
\end{equation}
where $\Phi_{d,\Omega}$ is the Green kernel of $-\Delta$ in $\Omega$ under the zero Dirichlet boundary condition and
\[
\cP_{k} := \Big\{ u\in  L^2(\Z^d): \int_{\Z^d}u\phi_{j}dx  =0 \ \text{ for \ } j=1,2,\cdots k-1 \ \text{ and } \ \|u\|_{L^2(\Z^d)}=1\Big\}.
\]

\item[(ii)] For $k\in[2,m_\Omega]\cap (0,+\infty)$ if $m_0\geq 2$,    
\[
 \lambda_{k-1} \le \lambda_{k}  
\]
 and
$$
  \lim_{k\to \infty}\lambda_{k} = +\infty\quad  {\rm if}\ \, m_0=+\infty.
$$

Moreover, for any $k\in\Z$, 
$$\lim_{|x|\to+\infty}\phi_k(x) = 0.  $$

\item[(iii)] There hold %The sequence of eigenfunctions  $\{\phi_{k}\}_{k\in\N}$ corresponding to eigenvalues $\{\lambda_{k}(\Z^d)\}$ form a complete orthonormal basis of $L^2(\Z^d)$ and 
$$ \int_{\Z^d}\phi_k \phi_j dx=0\quad{\rm and}\quad   \int_{\Z^d} (Q^{\frac12}\phi_k) \Phi_{d,\Omega}\ast(Q^{\frac12}\phi_j)   dx=0\quad {\rm for}\ \, k\not=j.  $$

\end{itemize}

    \end{theorem}

Theorem \ref{teo 1} and  Theorem \ref{teo 1-hf}  and Theorem \ref{teo 1-qd}  
establish the existence of solutions when $p$ is Sobolev supercritical and superlinear, via variational methods. We also establish the existence when $p$ is Serrin critical, supercritical, or sublinear, using the method of super- and subsolutions. In the Serrin subcritical case, we obtain the nonexistence of positive solutions through an iterative method based on decay estimates of solutions. These non-existnce   is also considered by \cite{GHS}, but our method is totally different, so we present our proof here.

In particular, nonexistence also holds in the Serrin critical and superlinear case. Recall that {\it it is open for the existence of positive solutions of the model equation with  Hardy potential
$$-\Delta u=(1+|x|)^{-\alpha}|u|^{p-2}u $$
for $1+\frac{d-\alpha}{d-2} < p\leq  \frac{2(d-\alpha)}{d-2} $ with $\alpha \in[0, 2)$
  in the whole domain $\Z^d$,   for $1+\frac{d+1-\alpha}{d-1} < p\leq  \frac{2(d-\alpha)}{d-2} $ with $\alpha \in[0, 2)$in the half domain $\Z^d_+$ and  for $1+\frac{d+2-\alpha}{d} < p\leq  \frac{2(d-\alpha)}{d-2} $ with $\alpha \in[0, 2)$in the half domain $\Z^d_{*}$.
  }
 the regions of $(\alpha,p) \subset [0,+\infty)\times (1,+\infty)$ corresponding to existence and nonexistence  results from Theorem \ref{teo 1}, Theorem \ref{teo 1-hf} and Theorem \ref{teo 1-qd} are illustrated in the following figures.
 
 \FloatBarrier
 \begin{figure}[h!]

 \centering
   \subfigure{
  \begin{minipage}{43mm}
   \includegraphics[scale=0.12]{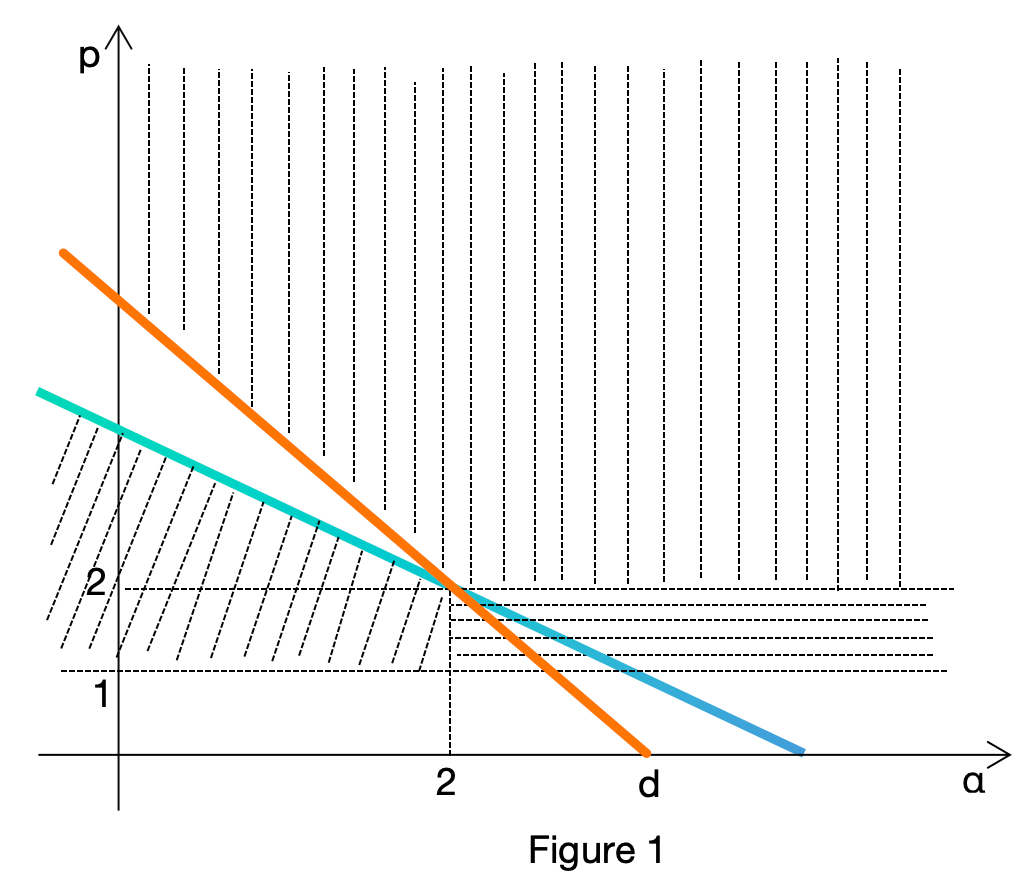}
 \end{minipage}
  }
    \subfigure{
    \begin{minipage}{43mm}
  \includegraphics[scale=0.12]{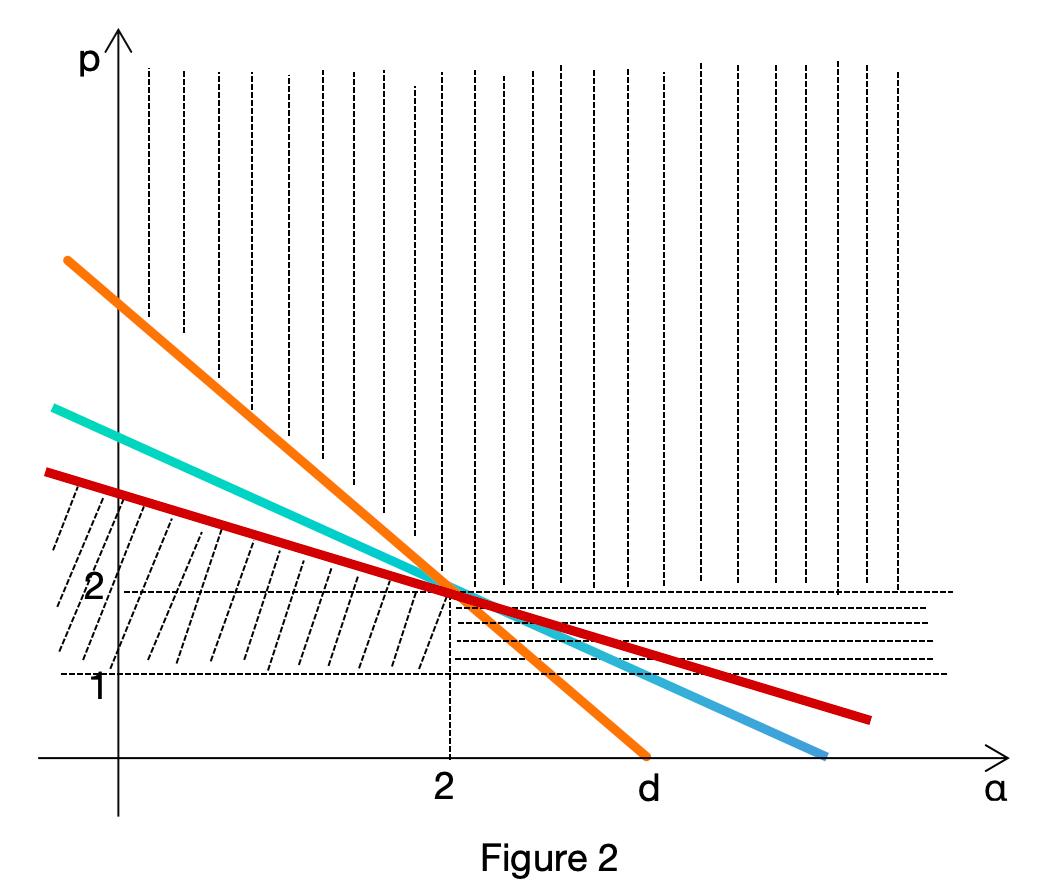}
  \end{minipage}
  }
    \subfigure{
    \begin{minipage}{43mm}
  \includegraphics[scale=0.12]{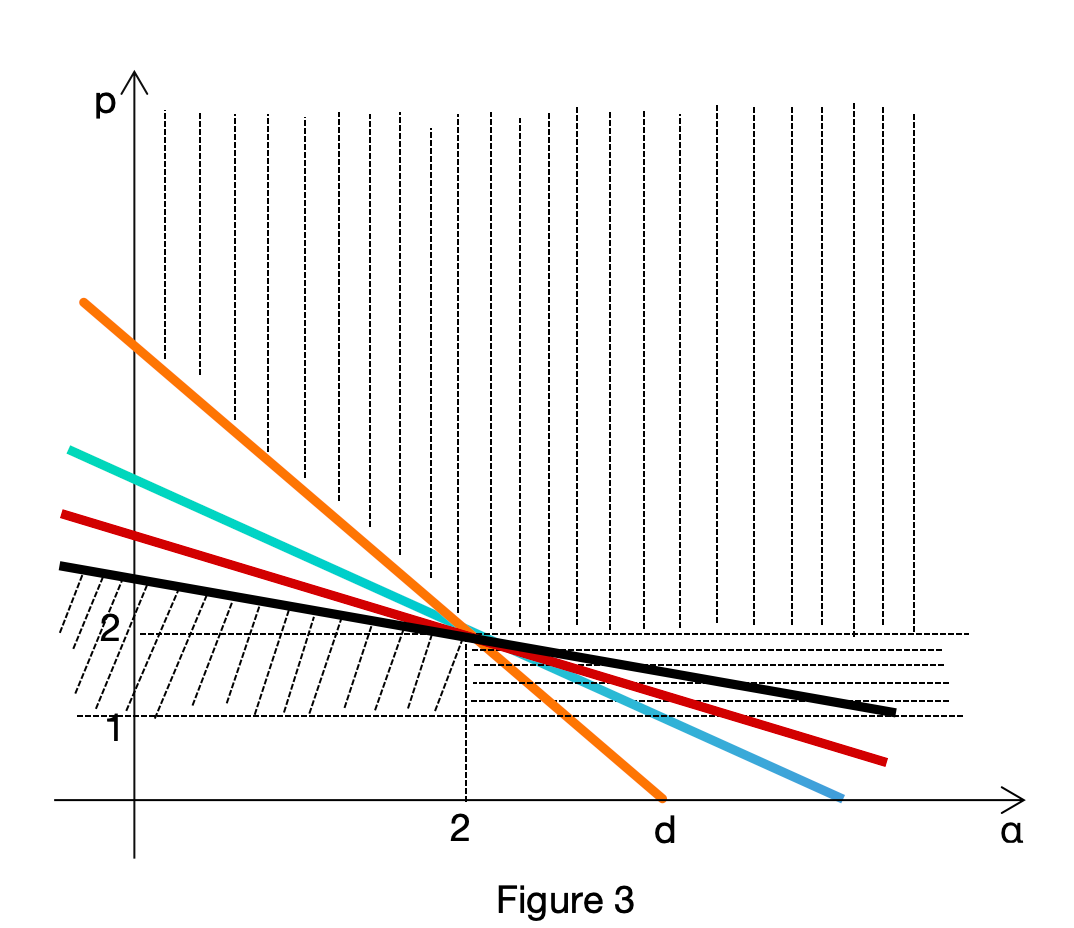}
 \end{minipage}
  }

  \vspace{10pt}
  {\footnotesize   %The blue line is the one of Serrin's exponent and the yellow is the line of Sobolev exponent.
  Figure 1,2,3 show the regions of $(\alpha,p)\subset [0,+\infty)\times [1,+\infty)$ in Theorem \ref{teo 1}, Theorem \ref{teo 1-hf} and Theorem \ref{teo 1-qd} in the domains  $\Z^d$, $\Z^d_+$  and $\Z^d_*$  respectively.
  
  Particularly, the blank regions $(\alpha,p)\subset [0,2)\times [1+\frac{d+k-\alpha}{d+k-2},2^*_\alpha)$ between the blue and  yellow lines in Figure 1 with $k=0$,    between the red and  yellow lines   in Figure 2 with $k=1$ and between the black and  yellow lines   in Figure 3  with $k=2$ are still open for the existence.  
  }
  \vspace{20pt}
 \end{figure}

  \FloatBarrier

 We emphasize that in the Sobolev supercritical case $p>\max\{2,2^*_{\alpha}\}$ or $p\in[ 2^*_{\alpha},+\infty)\cap (2,+\infty)$, our approach to derive the solution involves transforming the equations defined on three distinct domains into an integral equation by employing the corresponding fundamental solutions for these domains. Specifically, we consider the equation
$$u=\Phi_{d,\Omega} \ast (Q|u|^{p-2} u)\quad {\rm in}\ \, \Z^d,$$
where $\Phi_{d,\Omega}$ denotes the fundamental solution associated with the domains $\Omega$, here $\Phi_{d,\Omega}=0$ in $\Z^d\setminus \Omega$ if $\Omega\not=\Z^d$,   and $\ast$ stands the convolution for functions defined on the integers lattice graphs.
By introducing the substitution
$$v=Q^{\frac{1}{p'}} |u|^{p-2} u\quad {\rm in}\ \, \Z^d,$$
with $p'=\frac{p}{p-1}$, the conjugate exponent of $p$, then
the equation reduces to
$$|v|^{p'-2}v = Q^{\frac{1}{p}} \Phi_{d,\Omega} \ast (Q^{\frac{1}{p}} v) \quad {\rm in}\ \Z^d,$$
which possesses a variational structure. The corresponding energy functional is defined as
$$\cJ_0(v)=\frac{1}{p'}\int_{\Z^d} |v|^{p'} dx - \frac{1}{2} \int_{\Z^d} v\bK_{p,\Omega}(v) dx \quad {\rm for}\ v\in L^{p'}(\Z^d),$$
where $\bK_{p,\Omega}(v)=Q^{\frac{1}{p}} \Phi_{d,\Omega} \ast (Q^{\frac{1}{p}} v)$.
This framework allows us to apply the Mountain Pass Theorem to identify critical points of the energy functional $\cJ_0$. This variational formulation requires that $Q$ be bounded and nonnegative.

  The rest of this paper is organized as follows. In section 2, we analyze the basic properties of the related spaces and the estimates of the corresponding  Birman-Schwinger Operator.
 In Section 3, we show the existence of positive solution for the integral model, which is formed by the fundamental soluiotn.    Section 4 
 and Section 5 are devoted to  show the existence of positive solution in three types domains and the key point is to show the bounds of the fundamental solutions.  \smallskip

    \setcounter{equation}{0}
    \section{ Preliminary  }
\subsection{ Notations}
  In the sequel, we use the following notations:  $\Delta_{\Z^d}=\Delta$, and for $x \in \Z^d$,
$$|x|=\Big(\sum^d_{i=1}x_i^2\Big)^{\frac12}, \quad |x|_{_Q}=\sum^d_{i=1}|x_i|.$$
The gradient of $u$ at $x \in \Z^d$ is defined as
$$\nabla u(x)=\big( u(x+e_1)-u(x), \cdots, u(x+e_d)-u(x)   \big), $$
where $\{e_1, \cdots, e_d\}$ is the standard orthonormal basis of $\R^d$.
For $\emptyset\not= \Omega \subset \Z^d$,
$$\partial \Omega=\{y\in \Z^d\setminus \Omega:\exists \; x\in\Omega, x\sim y  \},\quad \bar \Omega=\partial \Omega\cup \Omega, \quad \Omega^c=\Z^d\setminus \Omega,$$
%$$\Omega^c=\Z^d\setminus \Omega\qquad {\rm for}\ \,\emptyset\not= \Omega \subset \Z^d,   $$
the ball
$$B_r(x^0)=\big\{x\in\Z^d: \text {$\exists\, n \, (\leq r)$ points $x^1,\cdots, x^n=x$ such that  $x^{i-1}\sim x^{i}$ for $i=1,\cdots, n$ } \big\},$$
$B_r=B_r(0)$
and the cube
 $$\bQ_\ell(x_0)=\Big\{x=(x_1,\cdots, x_d)\in\Z^d:\,  \sum_{i=1}^d|x_i-(x_0)_i|\leq \ell \Big\},\quad \ell>0.   $$
Let $C(\Z^d)$ with $d\geq 1$ be the set of all  functions $u:\Z^d\to\R$.  For $q\in[1,+\infty]$, let
$$L^q(\Z^d)=\{u\in C(\Z^d):\,  \|u\|_{L^q(\Z^d)}<+\infty \}$$
and
$$L^{q,+\infty}(\Z^d)=\{u\in C(\Z^d):\,  \|u\|_{L^{q,+\infty}(\Z^d)}<+\infty \}, $$
where
$$  \|u\|_{L^q(\Z^d)}=\big(\int_{\Z^d} |u(x)|^qdx\big)^{\frac1q}\ \ {\rm for}\ \,  q\in[1,+\infty),\qquad   \|u\|_{L^\infty(\Z^d)}= \sup_{x\in\Z^d} |u(x)|$$
and
$$  \|u\|_{L^{q,+\infty} (\Z^d)}= \sup_{\lambda>0} \big\{\lambda\cdot  \big|\{x\in\Z^d:\,  |u(x)| >\lambda  \}\big| ^{\frac1q}\big\}.  $$

A nonzero nonnegative function $\Phi_{d,\beta}:\Z^d\times\Z^d\to\R$ with $ 0<\beta<\frac{d}2$ and $d\geq 1$  satisfies that
 \begin{equation}\label{fund-up1}
  \Phi_{d,\beta}(x,y)= \Phi_{d,\beta}(y,x) \quad{\rm for}\ x,y\in\Z^d 
  \end{equation}
  and there is at least one point $\bar x\in\Z^d$ such that $\Phi_{d,\beta}(\bar x,\bar x)>0$.

For $f \in C(\Z^d)$, we denote
$$\Phi_{d,\beta}\ast f(x):=\int_{\Z^d}\Phi_{d,\beta}(x,y) f(y)dy $$
and let $\bK_{p,\beta}$ be the Birman-Schwinger operator \cite{FLW} associated to $\Phi_{d,\beta}$ and $Q$,
\begin{align}\label{def-K}
\bK_{p,\beta}[v]:=Q^{\frac1{p}}\Phi_{d,\beta}\ast  (Q^{\frac1{p}}v).
 \end{align}
Then we have
\begin{align}\label{eqq-00}
\int_{\Z^d} u\bK_{p,\beta}(v) dx=\int_{\Z^d} v\bK_{p,\beta}(u) dx,
 \end{align}
by the fact that
\begin{align*}
\int_{\Z^d} u\bK_{p,\beta}(v) dx&= \int_{\Z^d} (Q^{\frac1p}u)\, \Phi_{d,\beta}\ast (Q^{\frac1p}v ) dx
\\[1mm]&= \int_{\Z^d}\int_{\Z^d} (Q^{\frac1{p}}u)(x) (Q^{\frac1{p}}v)(y)\Phi_{d,\beta}(x,y)dxdy
%\\[1mm]&= \int_{\Z^d} (Q^{\frac1p}v)\, \Phi_{d,\beta}\ast (Q^{\frac1p}u ) dx
\\[1mm]&=\int_{\Z^d} v\bK_{p,\beta}(u) dx.
 \end{align*}
A Birman-Schwinger operator serves as a crucial tool in addressing elliptic problems involving polynomial nonlinearities and potentials. It is also widely employed in the study of spectral properties of operators, particularly within the contexts of quantum mechanics and the analysis of Schr\"odinger operators.

 For $\tau>0$, denote $h_\tau \in C^2(\R_+)$ by
$$h_\tau (t)= t^ {-\frac{\tau}2 }, \; \ \forall\, t>1. $$
Direct computation shows that
%\begin{align*}
$h_\tau'(t)=-\frac{\tau}2  t^{-\frac{\tau}2-1},\ \ h_\tau''(t)= \frac{\tau}2(\frac{\tau}2+1) t^{-\frac{\tau}2-2}, \ \forall\, t>1.$
%\end{align*}
Now we set
\begin{align}\label{fund-1}
\bar w_\tau(x)=h_\tau(|x|^2)\quad \text{ for $x\in\Z^d \setminus \{0\}$.}
\end{align}
Then for $x\in\Z^d$, $|x|$ large, we see that
\begin{align*}
 \Delta \bar w_\tau(x)   &= \sum_{y\sim x}\big(h_\tau(|y|^2)-h_\tau(|x|^2)\big)
 \\[1mm]&= \sum_{y\sim x} \Big[-\frac{\tau}2   |x| ^{-\tau-2}(|y|^2-|x|^2) +\frac14\tau(\frac{\tau}2+1)   |x| ^{-\tau-4}(|y|^2-|x|^2)^2  \Big]   (1+o(1))
  \\[1mm]&=-d\tau   |x| ^{- \tau -2}   + \frac14\tau(\frac{\tau}2+1)  |x| ^{- \tau -4}\big(8|x|^2+2d\big)(1+o(1))
   \\[1mm]&= \tau \big( \tau+2-d\big) |x| ^{- \tau   -2}+\frac14 d\tau( \tau +2)  |x| ^{- \tau -4}(1+o(1)),
\end{align*}
thus, for $|x|$ large
\begin{align}\label{est-1}
 -\Delta \bar w_\tau(x)   &=-\tau \big( \tau+2-d\big)|x|^{-\tau -2}-\frac14 d\tau( \tau+2) |x|^{-\tau-4} (1+o(1)).
\end{align}

\subsection{Basic properties}

We first state the maximum principle  in the discrete setting.
 \begin{theorem}\label{thm:max}
     Let  $\Omega\subset \Z^d$  be a connected domain    verifying
   either $ \partial  \Omega\not=\emptyset$ or $\Omega$ is unbounded,
     if  $u: \overline{\Omega}\to \R$ satisfies
\begin{equation}\label{eq 2.1 cm}
\left\{%\arraycolsep=1pt
\begin{array}{lll}
-\Delta u+\kappa  u     \geq  0   \quad
   &{\rm in}\ \  \Omega , \\[2mm]
 \phantom{ ----  }
u\geq 0 \quad &{\rm   in}\ \ \,    \partial \Omega, \\[1mm]
\liminf\limits_{x\in \Omega,\,  |x| \to \infty}u(x)\geq 0,
 \end{array}
 \right.
\end{equation} where $\kappa:\Omega\to[0,\infty)$ is a given function,
then $u\geq 0$ in $\Omega$. 

 Moreover,
either $u\equiv 0$ in $\Omega$ or  $u>0$ in $\Omega$.
 \end{theorem}
 \noindent{\bf Proof. }
 Without loss of generality, we prove it for an unbounded subset $\Omega.$ Since $\Omega$ is connected, so is $\overline{\Omega}.$ Suppose that the first assertion is not true, i.e. there exists $x_0\in\Omega$ such that $u(x_0)<0.$ Since $\liminf\limits_{x\in \Omega, |x| \to \infty }u(x)\geq 0
$ and $u|_{\partial\Omega}\geq 0,$ then $-\infty < \inf_{x\in \overline{\Omega}}u<0$ and
$$A:=\{x\in \overline\Omega:u(x)=\inf_{x\in\overline\Omega}u\}\neq\emptyset,\quad A\subsetneqq \Omega.$$
So there exists $x\in A$ such that there exists $y\in \Omega\setminus A$ and $y\sim x$,
then $u(y)>u(x)$ and $\Delta u(x)>0$,
while  by the equation,
    $$\Delta u(x)\leq \kappa(x)u(x)\leq 0.$$
 This is impossible.
 So $A=\emptyset$. So we obtain $u\geq 0$ in $\Omega$.

 Moreover, if there exists $\bar x\in\Omega$ such that $u(\bar x)=0,$ then by the same argument above, one can show that $u\equiv 0$ on $\overline\Omega.$ This proves the result. \hfill$\Box$\medskip

On the lattice space, the Hardy-type inequality  from \cite[Theorem 7.1]{KPP} (also see \cite{KPP1}) has 
the form that  
\begin{equation}\label{Htineq}
\int_{\Z^d} |\nabla u(x)|^2 dx= \sum_{x\in\Z^d}  \sum_{y\sim x} \big( \varphi(x)-\varphi(y)\big)^2 \geq \sum_{z\in\Z^d}  \bar{H}_0(z)\varphi(z)^2
 \end{equation}
 for all finitely supported function $\varphi$ on $\Z^d$, where 
  \begin{equation}\label{Green 1}
\bar H_0(x) =\big[\Delta\Phi_{d}^{\frac12}(x)\big]\Phi_{d}^{-\frac12}(x)  >0,
 \end{equation}
where   $\Phi_{d}$ is the fundamental solution of $-\Delta$ in $\Z^d$, 
i.e.
$$
 \left\{\arraycolsep=1pt
\begin{array}{lll}
-\Delta \Phi_{d} =\delta_0\quad
    {\rm in}\ \  \Z^d, \\[2mm]
 \phantom{   }
 \displaystyle \lim_{|x|\to+\infty}\Phi_{d}(x)=0. 
 \end{array}
 \right.
$$
Furthermore, we have that 
\begin{equation}\label{ome-1}
 \bar H_0(x)=\frac{(d-2)^2}{4}|x|^{-2}+O(\frac1{|x|^3})\quad {\rm as}\ \, |x|\to+\infty. 
 \end{equation}
  A survey on the Hardy-type inequality on graph could refer to \cite{Ff}
 and the reference therein. \smallskip

 \begin{lemma}\label{com lm-h1}
Assume that  $d\geq 3$ and $\mu> -1.$

Let $\Omega$ be a connected {\bf infinite set} of $\Z^d$ and $u: \Omega\to \R$ be a function fulfilling 
\begin{equation}\label{eq 2.1 cm}
\left\{\arraycolsep=1pt
\begin{array}{lll}
-\Delta u+ \mu \bar H_0  u     \geq  0   \quad
   &{\rm in}\ \  \Omega , \\[2mm]
 \phantom{ -----\  }
u\geq 0 \quad &{\rm   in}\ \ \,    \Z^d\setminus \Omega, 
\\[1mm]
\liminf\limits_{x\in \Omega,  |x| \to \infty}u(x)\geq 0,
 \end{array}
 \right.
\end{equation}
then $u\geq 0$ in $\Omega$.  Furthermore, 
either $u\equiv 0$ in $\Omega$ or  $u>0$ in $\Omega$. 
\end{lemma}
\noindent{\bf Proof. } 
By contradiction, if $\Omega^-:=\{x\in\Z^d:\,  u(x)<-\epsilon\}\not=\emptyset$,   then $u_{\epsilon,-}=-\max\{-(u(x)+\epsilon),0 \}$ has the support in $\Omega^-$ with finite vertices by the assumption that $\liminf\limits_{x\in \Omega,  |x| \to \infty}u(x)\geq 0$. 
  Direct computation shows that
\begin{align*} 
 &\quad\ \sum_{x,y\in \Z^d,\, x\sim y} \big( u(x)-u(y)\big) \big( u_{\epsilon,-}(x)-u_{\epsilon,-}(y)\big)  
  \\[1mm] &=\sum_{x,y\in \Z^d,\, x\sim y} \big( (u(x)+\epsilon) -(u(y)+\epsilon)\big) \big( u_{\epsilon,-}(x)-u_{\epsilon,-}(y)\big)  
 \\[1mm] &=  \sum_{x,y\in \Omega_-,\, x\sim y}   \big( u_{\epsilon,-}(x)-u_{\epsilon,-}(y)\big)^2+ \sum_{x\in \Omega_-,y\not\in \Omega_-,\, x\sim y}   \big( (u(x)+\epsilon)-(u(y)+\epsilon)\big) u_{\epsilon,-}(x) 
 \\[1mm] & \qquad -  \sum_{x\not\in \Omega_-,y \in \Omega_-,\, x\sim y}   \big( (u(x)+\epsilon)-(u(y)+\epsilon)\big) u_{\epsilon,-}(y)
  \\[1mm] &=  \sum_{x,y\in \Omega_-,\, x\sim y}   \big( u_{\epsilon,-}(x)-u_{\epsilon,-}(y)\big)^2+ 2 \sum_{x\in \Omega_-} \big( u_{\epsilon,-}(x)\big)^2 -2
  \sum_{x\in \Omega_-,y\not\in \Omega_-,\, x\sim y}(u(y)+\epsilon)u_{\epsilon,-}(x) 
  \\[1mm] &\geq  \sum_{x,y\in \Omega_-,\, x\sim y}   \big( u_{\epsilon,-}(x)-u_{\epsilon,-}(y)\big)^2+ 2 \sum_{x\in \Omega_-} \big( u_{\epsilon,-}(x)\big)^2
   \\[1mm] &= \sum_{x,y\in \Z^d,\, x\sim y}   \big( u_{\epsilon,-}(x)-u_{\epsilon,-}(y)\big)^2
 \end{align*}
 by the fact that $u+\epsilon \geq 0$ for $  \Z^d\setminus \Omega^-$.

If $u_-$ is zero, we are done. If not, multiply $ u_-$ in (\ref{eq 2.1 cm}) and sum over $\Omega$, then we obtain  a contradiction  that  
for $\mu>-1$, 
\begin{align*} 
 0&> \sum_{x\in\Omega}\big(-\Delta u(x)+ \mu\bar H_0(x)  u(x) \Big)u_-(x)
 \\[1mm]&=  \frac1{2}\sum_{x,y\in \Z^d,\, x\sim y} \big( u(x)-u(y)\big) \big( u_-(x)-u_-(y)\big)  + \mu \sum_{z\in\Z^d}  \bar H_0(z)u_-(z)^2 
\\[1mm]&\geq \frac1{2}\sum_{x,y\in \Z^d,\, x\sim y}   \big( u_-(x)-u_-(y)\big)^2 +\mu \sum_{z\in\Z^d}  \bar H_0(z)u_-(z)^2
 \\[1mm]&\geq0, 
 \end{align*}
where the last inequality holds by lthe Hardy inequality (\ref{Htineq}).  Thus, we have that $u+\epsilon\geq 0$ in $\Omega$
by the arbitrary of $\epsilon>0$. we get that $u \geq 0$ in $\Omega$. 

Now we assume that $x_0\in\Omega$ such that $u(x_0)=0$, then 
\begin{align*} 
 0\leq \big(-\Delta  + \mu H_0 \big) u(x_0) &= -\Delta u(x_0)
 =- \sum^d_{y\sim x_0} u(y) \leq 0, 
 \end{align*}
together with $u\geq 0$,  which leads to $u(y)=0$ for any $y\in \Omega$. 
As a consequence, we derive that either $u\equiv 0$ in $\Omega$ or  $u>0$ in $\Omega$.   \hfill$\Box$\medskip

We have also the following relationship between the different integrable functions spaces.

\begin{lemma}\label{lm 2.0}
$(i)$ If $u\in L^q(\Z^d)$ with $q\in[1,+\infty)$, then
$\displaystyle \lim_{|x|\to+\infty}u(x)=0. $

$(ii)$ For $1\leq q_1<q_2<+\infty$,
$$L^{q_1}(\Z^d)\subsetneqq L^{q_2}(\Z^d)\subsetneqq L^\infty(\Z^d)\subsetneqq C(\Z^d). $$

$(iii)$ For $1\leq q_1<+\infty$, we have that
$$L^{q_1}(\Z^d)\subsetneqq L^{q_1,\infty}(\Z^d).%{\color{blue} \subsetneqq L^{q_2}(\Z^d)}.
$$
Similar results hold for  $\Z^d_+$ and $\Z^d_*$.
\end{lemma}
\noindent{\bf Proof. }
Part $(i)$ and $(ii)$:   By contradiction,  let $u\in L^{q}(\Z^d)$ for   $q\in[ 1,\infty)$, and  assume that there is a sequence $(x_n)_n\subset \Z^d$
such that
$$|u(x_n)|\geq \sigma_0>0\quad \text{  for $n\geq n_0$} $$
for some  $\sigma_0>0$ and   $n_0>0$.
 Then   there holds
$$\int_{\Z^d} |u(x_n)|^q dx\geq \sum_{n\geq n_0} |u(x_n)|^q \geq \sigma_0 \sum_{n=n_0}^{+\infty} 1 =+\infty, $$
which implies that
\begin{equation}\label{limitzero}
\lim_{|x|\to+\infty}u(x)=0
\end{equation}
and
$L^{q}(\Z^d)\subset  L^\infty(\Z^d). $
Note that  $w_0(x)\equiv1$ for $x\in\Z^d$, then $w_0\in L^\infty(\Z^d)$
but $w_0\not\in L^{q}(\Z^d)$. Thus $L^{q}(\Z^d)\subsetneqq  L^\infty(\Z^d). $

 Now for $u\in L^{q_1}(\Z^d)\subset L^\infty(\Z^d)$, then
 $$\int_{\Z^d} |u(x)|^{q_2}dx\leq \|u\|_{L^\infty(\Z^d)}^{q_2-q_1} \int_{\Z^d} |u(x)|^{q_1}dx<+\infty, $$
 which leads to $u\in L^{q_2}(\Z^d)$. Thus, $L^{q_1}(\Z^d)\subset  L^{q_2}(\Z^d)$ and, obviously,  $L^{q_1}(\Z^d)\not=  L^{q_2}(\Z^d)$. \smallskip

 Part $(iii)$: For  given $u\in L^{q_1}(\Z^d)$ and any $\lambda>0$, let $E_\lambda=\{x\in\Z^d:\,  |u(x)| >\lambda \}$,
 then $|E_\lambda| < \infty$ by (\ref{limitzero}) and
 \begin{align*}
 \lambda |E_\lambda| ^{\frac1{q_1}}=\big( \lambda^{q_1} |E_\lambda|\big) ^{\frac1{q_1}}  \leq \big( \int_{E_\lambda} |u(x)|^{q_1}dx\big)^{\frac1{q_1}} ,
\end{align*}
 which implies that  $L^{q_1}(\Z^d)\subset   L^{q_1,\infty}(\Z^d)$. Moreover, letting $w_1(x)=(1+|x|)^{-\frac{d}{q_1}}$,
 then $w_1\in L^{q_1,\infty}(\Z^d)$, but it doesn't belong to $L^{q_1}(\Z^d)$. Therefore, $L^{q_1}(\Z^d)\subsetneqq   L^{q_1,\infty}(\Z^d)$. \hfill$\Box$\medskip

\subsection{Embedding  inequalities with Birman-Schwinger form   }

In this subsection, we assume  that 
\begin{itemize}
\item[  $(\bF_{0})$ ]
  for some $c>0$  
 \begin{equation}\label{embedding inequality-2}
 \|\Phi_{d,\beta} \ast f\|_{L^r(\Z^d)}\leq c\|  f\|_{L^q(\Z^d)}\ \,  \text{ for any $f\in L^{q}(\Z^d)$}
 \end{equation}
under the restriction that $\beta\in(0,\frac d{2})$ and %either
 $$\frac1r+\frac{2\beta}d\leq \frac1q \quad \text{with  }\ \,q, r\in(1,+\infty)  $$
 %or 
% $$\frac1r+\frac{2\beta}d<\frac1q \quad \text{with }\  \, q, r\in[1,+\infty].  $$

 \end{itemize}

\begin{lemma}\label{lm 2.2}
Assume that  $d\geq 1$, $\beta\in(0,\frac d2)$ and  $(\bF_0)$ holds. 
 
 Let 
 $v\in L^{p'}(\Z^d)$ and
 $Q\in L^{q_0,\infty}(\Z^d)$,
 where $p'=\frac{p}{p-1},\ p\in( 1,+\infty)$ and $q_0\in(1,+\infty)$ satisfy that %either
\begin{align}\label{eqq-01-1}
1<    \frac{q_0 p}{q_0p-q_0+1 } \leq \frac{2d}{d+2\beta }, 
 \end{align}
then   there exists $c>0$ independent of $v$ such that
 \begin{equation}\label{ei-1}
\Big|\int_{\Z^d} v\bK_{p,\beta}(v) dx\Big| \leq c\|v\|_{L^{p'}(\Z^d)}^2.
 \end{equation}
\end{lemma}
\noindent{\bf Proof. }  Let $p_1 \in \mathbb R$ be satisfying
\begin{align}\label{eqq-1}
1< p_1\leq \frac{2d}{d+2\beta},
\end{align}
which will be determinated below, then
$$%\begin{align}\label{eqq-1}
\frac1{p_1'}+\frac{2\beta}d\leq \frac1{p_1}.
$$%\end{align}
It follows by (\ref{embedding inequality-2}) and the weak H\"older inequality  that
\begin{align}
\Big| \int_{\Z^d} (Q^{\frac1p}v)\,  \Phi_{d,\beta}\ast (Q^{\frac1p}v ) dx\Big|  &\leq    \big\|Q^{\frac1p}v\big\|_{L^{p_1}(\Z^d)} \big\|\Phi_{d,\beta}\ast (Q^{\frac1p}v ) \big\|_{L^{p_1'}(\Z^d)}
\nonumber  \\[1mm] &\leq C \Big(\int_{\Z^d} Q^{\frac{p_1}p}|v|^{p_1} dx\Big)^{\frac{2}{p_1} }
\nonumber   \\[1mm]&\leq C \big\|Q^{\frac{p_1}p} \big\|_{L^{\theta',\infty}(\Z^d)}^{\frac{2}{p_1}}\,  \|v^{p_1}\|_{L^{\theta}(\Z^d)}^{\frac{2}{p_1}}
\nonumber  \\[1mm]&= C \big\|Q  \big\|_{L^{q_0,\infty}(\Z^d)}^{ 2p }  \,  \|v\|_{L^{p'}(\Z^d)}^{2 },   \label{kk-1}
\end{align}
where  $\theta\in(1,+\infty)$ has been chosen by
$$p_1\theta=p'\quad{\rm and}\quad  \frac{p_1}{p} \frac{\theta}{\theta-1}=q_0,   $$
when $q_0\in(1,+\infty)$, that is
$$ \theta=1+\frac{1 }{q_0(p-1)}  \quad{\rm and}\quad  p_1 =  \frac{q_0 p}{q_0(p-1)+1 }.  $$
%or by setting that $\theta=1$ when $q_0=+\infty$,  and $p_1=p'$ in this case.
  Taking (\ref{eqq-1}) into account,  we need   (\ref{eqq-01-1}). \hfill$\Box$\medskip

In the proof, the crucial point is (\ref{embedding inequality-2}) in the proof, so we have the following corollary. 
%\begin{corollary}\label{cr 2.2-}
% Assume  $d\geq 1$,     (\ref{embedding inequality-2}) holds by replacing $\Phi_{d,\beta}$ by $\Psi$. 
% Let  $v\in L^{p'}(\Z^d)$ and
 %$Q\in L^{q_0,\infty}(\Z^d)$,
 %where $p'\in[ 1,+\infty)$ and $q_0\in[1,+\infty]$ verify  either (\ref{eqq-01-1})
 % or (\ref{eqq-01-2}). 
 % Then   there exists $c>0$ such that
 %\begin{equation}\label{ei-1-psi}
%\Big|\int_{\Z^d}  (Q^{\frac1p}v)\, \Psi\ast (Q^{\frac1p}v ) dx\Big| \leq c\|v\|_{L^{p'}(\Z^d)}^2.
 %\end{equation}
%\end{corollary}

\begin{corollary}\label{lm 2.4}
 Let  $\beta\in(0,\frac d2)$, and  non-negative function $Q$ veriy that  for some  $\tilde \alpha\in[0,+\infty)$
\begin{align}\label{ekk-1}
\limsup_{|x|\to+\infty} Q(x)|x|^{\tilde\alpha}<+\infty
 \end{align}
 and
  $v\in L^{p'}(\Z^d)$ with
  that either
\begin{align}\label{eqq-01-1g}
1\leq    \frac{d p}{d p-d+\tilde\alpha } \leq \frac{2d}{d+2\beta}\quad {\rm for\ }\  \tilde\alpha \in(0,  d ]
 \end{align}
 or
 \begin{align}\label{eqq-01-2g}
  p' \geq 1\quad {\rm for\ }\   \tilde\alpha >d.
 \end{align}
   Then   there exists $c>0$ such that (\ref{ei-1}) holds true.
\end{corollary}
\noindent{\bf Proof. }    It follows by (\ref{ekk-1}) that
  $Q\in L^{\frac{d}{\tilde\alpha},\infty}(\Z^d)$  and then conditions
(\ref{eqq-01-1})  and (\ref{eqq-01-2}) are equivalent to (\ref{eqq-01-1g}) and (\ref{eqq-01-2g}) respectively. \hfill$\Box$\medskip

Here we show a sufficient condition for (\ref{embedding inequality-2}).

\begin{lemma}\label{lm 2.1}
Assume that $\beta\in(0,\frac d2)$, 
 $f\in L^{q}(\Z^d)$ with $1\leq q\leq \frac{d}{2\beta}$.  
 
If   there exists some $c>0$ such that 
 \begin{equation}\label{bound-fun-1}
 |\Phi_{d,\beta}(x,y)|\leq c (1+|x-y|)^{2\beta-d} \quad{\rm for}\  x,y\in\Z^d. 
   \end{equation}
 then  (\ref{embedding inequality-2}) holds  either for
 $$\frac1r+\frac{2\beta}d\leq  \frac1q\quad {\rm with}\ \,   q, r\in(1,+\infty) $$
 or for
   $$\frac1r+\frac{2\beta}d<  \frac1q\quad {\rm with}\ \,  q, r\in[1,+\infty].$$
\end{lemma}

\noindent{\bf Proof. }  We extend continuously the functions $\Phi_{d,\beta}$ and $f$  to $\R^n$, still denote by
    $\Phi_{d,\beta},f$ respectively,
    such that for $x\in\R^d$
    $$\min_{x'\in \Z^d, |x'-x|\leq \sqrt{d}}\Phi_{d,\beta} (x')\leq  \Phi_{d,\beta}(x)\leq \max_{x'\in \Z^d, |x'-x|\leq \sqrt{d}}\Phi_{d,\beta} (x')$$
    and
 $$\min_{x'\in \Z^d, |x'-x|\leq \sqrt{d}} f(x')   \leq f(x) \leq \max_{x'\in \Z^d, |x'-x|\leq \sqrt{d}} f(x').$$
 Indeed, we can systematically  do the extension first to the line segments connecting adjacent vertices that satisfy the aforementioned bounds, then to the squares bounded by these line segments, and finally to the 3-dimensional cubes formed by those squares—thereby obtaining the desired extension through a well-defined, hierarchical process to the $d$-dimensional cubes.

It follows by (\ref{fund-up1}) and Lemma \ref{lm 2.0} that  $f\in L^\infty(\R^d)\cap L^{q}(\R^d)$  and
\begin{align*}
\big|(\Phi_{d,\beta} \ast f)(x) \big|  \leq c_0 \int_{\Z^d}(1+ |x-y|)^{2\beta-d} |f(y)|dy \leq c_1\int_{\R^d}(1+ |x-y|)^{2\beta-d} |f(y)|dy .
\end{align*}

 Let $$\tilde\Phi_{d,\beta}(z)=(1+ |z|)^{2\beta-d},\quad\forall\, z\in\R^d, $$
     then
     $\tilde \Phi_{d,\beta}\in L^\infty(\R^d)\cap L^{s}(\R^d)$ for
     $s>\frac{d}{d-2\beta}\in(1,+\infty)$ for $\beta\in(0,\frac d2)$
    and 
      $\tilde \Phi_{d,\beta}\in L^\infty(\R^d)\cap L^{\frac{d}{d-2\beta},\infty}(\R^d)$ with 
     $\frac{d}{d-2\beta}\in(1,+\infty)$ for $\beta\in(0,\frac d2)$
   and 
   it follows by the Young's inequality for convolution that
\begin{align*}
\|\Phi_{d,\beta} \ast f\|_{L^r(\Z^d)} &\leq  c\|\tilde \Phi_{d,\beta} \ast f\|_{L^r(\R^d)}
\\[1mm]& \leq    c'\|f\|_{L^q(\R^d)} \| \tilde \Phi_{d,\beta}\|_{L^{s}(\R^d)}
 \leq c'' \|f\|_{L^q(\Z^d)},
\end{align*}
     where 
     $$\frac1q+  \frac1s=1+\frac1r\quad {\rm with}\quad s> \frac{d}{d-2\beta}\ \,  {\rm and}\ \, r, q\in [1,+\infty].   $$
   
By the weak type  Young's inequality for convolution that
\begin{align*}
\|\Phi_{d,\beta} \ast f\|_{L^r(\Z^d)} &\leq  c\|\tilde \Phi_{d,\beta} \ast f\|_{L^r(\R^d)}
\\[1mm]& \leq    c'\|f\|_{L^q(\R^d)} \| \tilde \Phi_{d,\beta}\|_{L^{s,\infty}(\R^d)}
 \leq c'' \|f\|_{L^q(\Z^d)},
\end{align*}
     where 
     $$\frac1q+  \frac1s=1+\frac1r\quad {\rm with}\quad s\geq \frac{d}{d-2\beta}\ \,  {\rm and}\ \, q,r\in(1,+\infty).   $$
   We complete the proof.  \hfill$\Box$ \medskip

Finally, we show positivity of with Birman-Schwinger form. 
\begin{lemma}\label{lm 2.3}
 Let $\Phi_{d,\beta}$ be the fundamental solution of $-\Delta \Phi_{d,\beta}(\cdot,y)=\delta_y$ in $\Omega$ for $y\in \Omega\subset \Z^d$,  subject to the zero Dirichlet condition $\Phi_{d,\beta}(\cdot,y)=0$ in $\Z^d\setminus \Omega$, if $\Omega\not=\Z^d$.  
 Assume that 
 $$\big|\int_{\Z^d}  v\bK_{p,\beta}(v)  dx\big|<+\infty$$
 Then  for any  $v\in L^{p'}(\Z^d)$, we have
  $$% \begin{equation}\label{ei-2}
  \int_{\Z^d}  v\bK_{p,\beta}(v)  dx\geq 0.
 $$%\end{equation}
If we assume more that
   \begin{align}\label{cond-p1}
   {\rm supp}\big(Q^{\frac1p} v\big)\cap \Omega \not=\emptyset,  
       \end{align}
 then
 $$% \begin{equation}\label{ei-2}
  \int_{\Z^d}  v\bK_{p,\beta}(v)  dx>0.
$$%\end{equation}
 \end{lemma}
 \noindent{\bf Proof. }
 Let $$u=\Phi_{d,\beta}\ast (Q^{\frac1p} v)\quad\text{ in $\Z^d$.}$$
Then we obtain that
  \begin{align*}
  \int_{\Z^d}  v\bK_{p,\beta}(v)\,  dx=\int_{\Omega} u (-\Delta) u\, dx
  =\int_{\Omega} |\nabla u|^2dx\geq 0.
    \end{align*}
  If $\Omega\not=\Z^d$,  by (\ref{cond-p1}), 
  ${\rm supp}\,u\subset \Omega$  and  $u$ is not a constant  in $\Z^d$, then
   \begin{align*}
  \int_{\Z^d}  v\bK_{p,\beta}(v)  dx =\int_{{\Omega}}|\nabla u|^2dx> 0.
    \end{align*}
    If $\Omega=\Z^d$ and $u$ is a constant in $\Z^d$, then $Q^{\frac1p}v=0$ in $\Z^d$, which
    contradicts  (\ref{cond-p1}), so if $u$ is not a constant, then \begin{align*}
  \int_{\Z^d}  v\bK_{p,\beta}(v)  dx =\int_{\Z^d}|\nabla u|^2dx> 0.
    \end{align*}
  We complete the proof.  \hfill$\Box$ 

    \setcounter{equation}{0}
\section{ Existence for  integral equations  }
In this subsection, we consider the existence of positive solution to the integral equations
   \begin{align}\label{eqq-3.2}
|v|^{p'-2}v=Q^{\frac1{p}} \Phi_{d,\beta}\ast (Q^{\frac1{p}}v) \qquad {\rm in}\ \, \Z^d,
\end{align}
where $p\geq 2$ and $p'=\frac{p}{p-1}$. In this section, we always assume that
\begin{itemize}
\item[  $(\bF_{1})$ ]
 $\Phi_{d,\beta}:\Z^d\times\Z^d\to[0,+\infty)$   satisfies     (\ref{embedding inequality-2}) with $ 0<\beta<\frac{d}2$ and $d\geq 1$ and 
$$\Phi_{d,\beta}(x,y)=\Phi_{d,\beta}(y,x)\quad{\rm for}\ \,  x,y\in\Z^d; $$

\item[  $(\bF_{2})$ ]
$$% \begin{equation}\label{ei-2}
  \int_{\Z^d}  v\Phi_{d,\beta}\ast v \,  dx\geq 0\quad \text{for all $v\in C_c(\Z^d)$}$$ and
$$ \int_{\Z^d}  v\Phi_{d,\beta}\ast v  dx>0\quad \text{for all $v\in C_c(\Z^d)$ and $\rm{supp}(v)\cap  \big\{z\in \Z^d: \exists y\in\Z^d:\,   \Phi_{d,\beta}(z,y)>0\big\}\not=\emptyset$}. 
 $$

 \end{itemize}

For $\beta\in(0,\frac d2)$, denote 
 $$2^*_{\beta,\alpha}:=\frac{2(d-\alpha)}{d-2\beta}.  $$
 Involving a nonnegative  $Q$, $(\bF_{2})$ implies that  
  \begin{equation}\label{ei-0-2}
  \int_{\Z^d}  v\bK_{p,\beta}(v)  dx\geq 0\quad \text{for all $v\in L^{p'}(\Z^d)$}  
  \end{equation}
   and
  \begin{equation}\label{ei-0-3}
   \int_{\Z^d}  v\bK_{p,\beta}(v)  dx>0\quad \text{for all $v\in L^{p'}(\Z^d)$ and $\rm{supp}(v)\cap \cM_0\not=\emptyset$}, 
  \end{equation}
 where
\begin{equation}\label{ei-0-4}
\cM_0:= \{x\in\Z^d:\,  \exists y\in\Z^d:   Q(x)Q(y)>0,\, \Phi_{d,\beta}(x,y)>0\}\not=\emptyset. 
  \end{equation}
Here we denote $m_0$  the number of the points in $\cM_0$, then we have that 
\begin{equation}\label{ei-0-5}
m_0=|\cM_0|\leq +\infty.
  \end{equation}

\subsection{Super-linear   case:  $p>2$ }

 In the super linear case,   we need to find out the sharp range of the exponent
  of the nonlinearity,  which depends on  the potential function.

\begin{itemize}
\item[ $(\bA_{\alpha,\beta,1})$]
Let
$$\alpha\in[0,+\infty), \quad  \beta\in(-\infty,\frac d2) \quad{\rm and}\quad p \in\big[2,+\infty\big)\cap  \big(2^*_{\beta,\alpha}, +\infty\big). $$

 If $\alpha=0$, we assume more that
  $$\limsup_{|x|\to+\infty} Q(x) <+\infty. $$
 \item[$(\bA_{\alpha,\beta,2})$] Let
  \begin{align}\label{ass-3.1}
\lim_{|x|\to+\infty} Q(x)|x|^{\alpha}=0
\end{align}
and
$$\alpha\in[0,+\infty), \quad  \beta\in(-\infty,\frac d2) \quad{\rm and}\quad p\in\big[2,+\infty\big)\cap  \big[2^*_{\beta,\alpha}, +\infty\big).    $$

 %\item[$(A1)$] Let    $Q\in L^{q_0, \infty}(\Z^d)$ with $q_0\in[1,+\infty]$ and
% $$p>\max\Big\{2,\,  \frac{2d}{d-2} \frac{q_0-1}{q_0}\Big\};   $$

\end{itemize}

\begin{theorem}\label{teo 3.1}
  Assume  that $d\geq 1$, $\beta\in(0,\frac d2)$, $\Phi_{d,\beta}$ satisfies $(\bF_1)$ and $(\bF_2)$ and 
  $\alpha\geq0$,  $p>2$ verify  either  $(\bA_{\alpha,\beta,1})$ or $(\bA_{\alpha,\beta,2})$.
  Then
  problem (\ref{eqq-3.2})
 has at least one nontrivial positive solution $v\in L^{p'}(\Z^d)$
 and
   $\displaystyle  \lim_{|x|\to+\infty}v(x) = 0.  $
    \end{theorem}

For the existence of solution of (\ref{eqq-3.2}), notice that (\ref{eqq-3.2}) has the variational structure in $L^{p'}(\Z^d)$ and the existence of solution will be derived by Mountain Pass Theorem (for instance, \cite[Theorem~6.1]{struwe}; see also \cite{ar, rabinowitz}).
Via the fundamental function, the above equations we consider  can be transformed
  into the integral models and  the associated energy functional could be established 
\begin{align}\label{eqq-3.3}
\cJ_0(v)=\frac1{p'}\int_{\Z^d} |v|^{p'} dx-\frac12 \int_{\Z^d}  v\bK_{p,\beta}(v)  dx\quad{\rm for}\ \,  v\in L^{p'}(\Z^d),
\end{align}
where $\bK_{p,\beta}$ is defined in (\ref{def-K}). Moreover,  we have that $\cJ_0\in C^1(L^{p'}(\Z^d),\R)$ and
\begin{align}\label{eqq-3.3-d}
\langle\cJ_0'(v),w\rangle= \int_{\Z^d}\big( |v|^{p'-2}v-\bK_{p,\beta}(v)\big) w dx \quad{\rm for}\ \,  v,w\in L^{p'}(\Z^d).
\end{align}

We need to prove the following

 \begin{proposition}\label{pr 3.1}
 Let the assumptions of Theorem \ref{teo 3.1} hold,  then 

 $(i)$ There exist $\delta>0$ and $\rho\in(0,1)$ such that
 $$\cJ_0(v)\geq \delta\quad \text{for all $v\in L^{p'}(\Z^d)$ with $\|v\|_{L^{p'}(\Z^d)}=\rho$}. $$

 $(ii)$ There is $v_0\in L^{p'}(\Z^d)$ such that $\|v_0\|_{L^{p'}(\Z^d)}>1$ and $\cJ_0(v_0)<0$.

 $(iii)$ Every Palais-Smale sequence $(v_n)_n$ of $\cJ_0$ verifying
 $$\cJ_0(v_n)\to c\not=0, $$
  up to translation,  has a converging subsequence in $L^{p'}(\Z^d)$.

 \end{proposition}

The proof of Proposition \ref{pr 3.1} is based on the following auxiliary  non-vanishing property, where the exact meaning can be stated as follows.
%We now introduce the  non-vanishing property.
%Recall that
%$$A^c=\Z^d\setminus A\quad {\rm for}\ \, A\subset \Z^d $$
%and
%$$\bQ_\ell(x_0)=\Big\{x=(x_1,\cdots, x_d)\in\Z^d:\,  \sum_{i=1}^d|x_i-(x_0)_i|\leq \ell \Big\} $$
%for $\ell\geq 1$ and $x_0=\big((x_0)_1,\cdots,(x_0)_d\big)\in\Z^d$.

\begin{lemma}\label{lm 3.1}
 Let the assumptions of Theorem \ref{teo 3.1} hold,  $(v_n)_n\subset L^{p'}(\Z^d)$ be a bounded sequence such that
$$\limsup_{n\to+\infty}  \int_{\Z^d} v_n\bK_{p,\beta}(v_n) dx   >0,  $$
then there are $R>0,\ n_0\geq1,\  \epsilon_0>0$ and $(x_n)_n\subset \Z^d$ such that, up to subsequence,
$$\int_{\bQ_R(x_n)} |v_n|^{p'} dx\geq \epsilon_0\quad {\rm for\ all}\ n\geq n_0. $$
\end{lemma}
\noindent{\bf Proof. }
We prove the following variant: {\it if for any $R>0$,
\begin{align}\label{eqq-3.4}
\lim_{n\to+\infty} \Big(\sup_{y\in \Z^d} \int_{\bQ_R(y)} |v_n|^{p'} dx\Big)=0,
\end{align}
 then
\begin{align}\label{eqq-3.5}
\lim_{n\to+\infty}  \int_{\Z^d}  v_n \bK_{p,\beta}(v_n) dx  =0.
\end{align}
 }

 {\it Part 1: }  Under the assumption $(\bA_{\alpha,\beta,1})$, $p \in\big[2,+\infty\big)\cap  \big(2^*_{\beta,\alpha}, +\infty\big)$,
% $Q_1$ is uniformly bounded in $\Z^d$, since
%$$p>\max\Big\{2,\,  \frac{2d}{d-2\beta} \frac{d-\alpha}{d }\Big\},  $$
 we can choose $\alpha_1<\alpha$   such that
$$p= \max\Big\{2,\,  2^*_{\beta,\alpha_1}\Big\}. $$
%Give   $\epsilon>0$ and let
Let $$Q_{\alpha_1} (x)=(1+|x|)^{- \alpha_1 },\quad Q_1(x)=\frac{Q(x)}{Q_{\alpha_1}(x)}\quad{\rm for}\ \, x\in\Z^d. $$
 Then $Q_1$ is uniformly bounded in $\Z^d$. Moreover there exists $C>0$ such that
 $$ Q_1(x)\leq C(1+ |x|)^{-\frac{\alpha-\alpha_1}{2}} \quad{\rm in} \;\, \Z^d,$$
and for any $R>1$,  it follows  by Lemma \ref{lm 2.1}  that, with $v=v_n$ for any $n$ and $p$ satisfies the assumption $(\bA_{\alpha,\beta,1})$,
 \begin{align*}
   \Big| \int_{\Z^d}\int_{\Z^d} & (Q^{\frac1{p}}v)(x) (Q^{\frac1{p}}v)(y)\Phi_{d,\beta}(x,y)1_{B_R(0)^c}(x-y) dxdy \Big|
 \\[1mm]&\leq
   \int_{\Z^d}\int_{\Z^d}|v(x)|  |v(y)|\big( Q(x)^{\frac1{p}}Q(y)^{\frac1{p}} \Phi_{d,\beta}(x,y)1_{B_R(0)^c}(x-y) \big) dxdy
   \\[1mm] &\leq C  \|Q_1\|_{L^\infty(\Z^d)}^2 (1+\frac R2)^{-\frac{\alpha-\alpha_1}{2}\frac1p}   \int_{\Z^d}\int_{\Z^d} Q_{\alpha_1}^{\frac1{p}}(x)Q_{\alpha_1}^{\frac1{p}}(y)|v(x)| \, |v(y)|  \Phi_{d,\beta}(x,y)   dxdy
    \\[1mm] &\leq C \|Q_1\|_{L^\infty(\Z^d)}^2 (1+\frac R2)^{-\frac{\alpha-\alpha_1}{2}\frac1p} \|v\|_{L^{p'}(\Z^d)} ^2.
 \end{align*}
  Then for any $\epsilon > 0$, there exist $R_\epsilon>1$ and $C>0$ such that  for $R\geq R_\epsilon$
  \begin{align} \label{part e-0}
 \Big|  \int_{\Z^d}\int_{\Z^d}(Q^{\frac1{p}}v)(x) (Q^{\frac1{p}}v)(y)\Phi_{d,\beta}(x,y)1_{\bQ_R(0)^c}(x-y) dxdy \Big| \leq C \epsilon.
 \end{align}

%\smallskip
Under the assumption (\ref{ass-3.1}) in $(\bA_{\alpha,\beta,2})$, we take $\alpha_1=\alpha$ and
 $$\lim_{|x|\to+\infty}Q_1(x)=0$$  and
for any $R>0$, it follows  by Lemma \ref{lm 2.1}  that
 \begin{align*}
   \Big| \int_{\Z^d}\int_{\Z^d} & (Q^{\frac1{p}}v)(x) (Q^{\frac1{p}}v)(y)\Phi_{d,\beta}(x,y)1_{B_R(0)^c}(x-y) dxdy \Big|
     \\[1mm] &\leq \Big(\sup_{z\in \bQ_{\frac R2}(0)^c} Q_1(z)\Big) \|Q_1\|_{L^\infty(\Z^d)} \int_{\Z^d}\int_{\Z^d} Q_{\alpha}^{\frac1p}(x)Q_{\alpha}^{\frac1p}(y)|v(x)| \, |v(y)|  \Phi_{d,\beta}(x,y)   dxdy
    \\[1mm] &\leq c\Big(\sup_{z\in \bQ_{\frac R2}(0)^c} Q_1(z)\Big) \|Q_1\|_{L^\infty(\Z^d)} \|v\|_{L^{p'}(\Z^d)} ^2.
 \end{align*}
There exists $R_\epsilon>1$ such that  for $R\geq R_\epsilon$
$$\sup_{z\in \bQ_{\frac R2}(0)^c} Q_1(z)\leq\epsilon, $$
which implies that for $R=R_\epsilon$ and $C>0$
  \begin{align} \label{part e-00}
 \Big|  \int_{\Z^d}\int_{\Z^d}(Q^{\frac1{p}}v)(x) (Q^{\frac1{p}}v)(y)\Phi_{d,\beta}(x,y)1_{\bQ_R(0)^c}(x-y) dxdy \Big| \leq C \epsilon.
 \end{align}

 {\it Part 2: }
 For $R=R_\epsilon$, there exists a sequence of  points $(z_\ell)_{\ell \geq 1} \subset \Z^d$ such that
 $$\bQ_R(z_\ell)\cap \bQ_R(z_{\ell'})=\emptyset\ \ {\rm if}\ \, \ell\not=\ell'\quad{\rm and}\quad   \Z^d=\bigcup_{\ell \geq 1} \bQ_R(z_\ell).$$

By (\ref{eqq-3.4}) with $R=R_\epsilon$,
 we obtain that
  \begin{align*}
   \Big| \int_{\Z^d}\int_{\Z^d} & (Q^{\frac1{p}}v)(x) (Q^{\frac1{p}}v)(y)\Phi_{d,\beta}(x,y)1_{\bQ_R(0)}(x-y) dxdy \Big|
    \\[1mm] &\leq  \sum^\infty_{\ell=1} \int_{\bQ_R(z_\ell)} \Big(\int_{\bQ_R(x)}  (Q^{\frac1{p}}|v|)(x) (Q^{\frac1{p}}|v|)(y)\Phi_{d,\beta}(x,y)1_{\bQ_R(0)}(x-y)dy\Big)  dx
    \\[1mm] &\leq 2\|Q\|_{L^\infty(\Z^d)}^2 \|\Phi_{d,\beta}\|_{L^\infty(\Z^d)}   \sum^\infty_{\ell=1} \int_{\bQ_{R}(z_\ell)} \Big(\int_{\bQ_{3R}(z_\ell)}  |v(x)||v(y)|dy  \Big) dx
             \\[1mm] &\leq C \|Q\|_{L^\infty(\Z^d)}^2 \|\Phi_{d,\beta}\|_{L^\infty(\Z^d)} R^{\frac{2d}{p}}  \sum^\infty_{\ell=1}  \Big(\int_{\bQ_{3R}(z_\ell)} |v(y)|^{p'}dy  \Big)^{\frac2{p'}}
       \\[1mm] &\leq C \|Q\|_{L^\infty(\Z^d)}^2 \|\Phi_{d,\beta}\|_{L^\infty(\Z^d)} R^{\frac{2d}{p}} \Big( \sup_{\ell\in\N} \int_{\bQ_{3R}(z_\ell)} |v(y)|^{p'}dy   \Big)^{\frac2{p'}-1}   \sum^\infty_{\ell=1}  \Big(\int_{\bQ_{3R}(z_\ell)} |v(y)|^{p'}dy  \Big)
         \\[1mm] &\leq C' \|Q\|_{L^\infty(\Z^d)}^2 \|\Phi_{d,\beta}\|_{L^\infty(\Z^d)} R^{\frac{2d}{p}} \Big( \sup_{z_\ell \in\Z^d} \int_{\bQ_{3R}(z_\ell)} |v(y)|^{p'}dy   \Big)^{\frac2{p'}-1}  \Big(  \int_{\Z^d} |v(y)|^{p'}dy  \Big),
 \end{align*}
 then by (\ref{eqq-3.4}),  there exists an integer $n_R>0$ such that for $n\geq n_R$
  \begin{align} \label{part e-1}
   \Big| \int_{\Z^d}\int_{\Z^d}(Q^{\frac1{p}}v)(x) (Q^{\frac1{p}}v)(y)\Phi_{d,\beta}(x,y)1_{\bQ_R(0)}(x-y) dxdy \Big| \leq \epsilon,
 \end{align}
 which,  together with (\ref{part e-0}),  implies that for any $\epsilon>0$,  there is $n_\epsilon>0$ such that
  \begin{align*}
   \Big| \int_{\Z^d}  v\bK_{p,\beta}(v) dx  \Big|=  \Big| \int_{\Z^d}\int_{\Z^d}(Q^{\frac1{p}}v)(x) (Q^{\frac1{p}}v)(y)\Phi_{d,\beta}(x,y) dxdy \Big| \leq \epsilon \quad {\rm for}\ \ n\geq n_\epsilon.
 \end{align*}
 Thus, we obtain (\ref{eqq-3.5}) as claimed. \hfill$\Box$\medskip

 \begin{lemma}\label{cr 3.1}
 Under the assumptions of Lemma \ref{lm 3.1},  suppose that $v_n\rightharpoonup v $ in $L^{p'}(\Z^d)$,
 then
 $$\int_{\Z^d} v_n\bK_{p,\beta}(v_n-v) dx\to 0\quad {\rm as}\quad n\to+\infty. $$

 \end{lemma}
{\bf Proof.}  For simplicity, we can assume that $v=0$.
Since $v_n\rightharpoonup 0 $ in $L^{p'}(\Z^d)$, then  $\|v_n\|_{L^{p'}(\Z^d)}$ is bounded,
 $v_n\to 0 $ in $L^{p'}_{\rm loc}(\Z^d)$,  that is, for any $R>1$ and any $y\in \Z^d$, we have that
 \begin{align}\label{eqs-3.1}
\lim_{n\to+\infty}   \int_{\bQ_R(y)} |v_n|^{p'} dx =0,
\end{align}

{\it Part I:}  Give   $\epsilon>0$ and recall
$$Q_{\alpha_1} (x)=(1+|x|)^{- \alpha_1 },\quad Q_1(x)=\frac{Q(x)}{Q_{\alpha_1}(x)}\quad{\rm for}\ \, x\in\Z^d,  $$
where $\alpha_1\leq \alpha$  such that
$p= \max\Big\{2,\, 2^*_{\beta, \alpha_1} \Big\}.$
%\frac{2d}{d-2\beta} \frac{d-\alpha_1}{d }\Big\} $$
Under the assumption $(\bA_{\alpha,\beta,1})$, $Q_1$ is uniformly bounded in $\Z^d$.

For any $R>1$, it follows  by Lemma \ref{lm 2.1}  that
 \begin{align*}
\Big|\int_{\Z^d} v_n 1_{B_R(0)^c} \bK_{p,\beta}(v_n) dx \Big| & =   \Big| \int_{\Z^d}\int_{\Z^d}(Q^{\frac1{p}}v_n)(x) (Q^{\frac1{p}}v_n)(y)\Phi_{d,\beta}(x,y)1_{B_R(0)^c}(x) dxdy \Big|
 \\[1mm]&\leq
  \int_{B_R(0)^c} \int_{\Z^d}|v_n(x)|  |v_n(y)|\big( Q(x)^{\frac1{p}}Q(y)^{\frac1{p}} \Phi_{d,\beta}(x,y) \big) dxdy
   \\[1mm] &\leq  c\|Q_1\|_{L^\infty(\Z^d)}^2  R^{-\frac{\alpha-\alpha_1}2\frac1p}  \int_{\Z^d}\int_{\Z^d} Q_{\alpha_1}(x)^{\frac1p}Q_{\alpha_1}(y)^{\frac1p}|v_n(x)|  \, |v_n(y)|  \Phi_{d,\beta}(x,y)   dxdy
    \\[1mm] &\leq c  \|Q_1\|_{L^\infty(\Z^d)}^2   \|v_n\|_{L^{p'}(\Z^d)}^2  (1+  R)^{-\frac{\alpha-\alpha_1}2\frac1p}  .
 \end{align*}
Then there exists $R_\epsilon>1$ such that  for $R\geq R_\epsilon$
  \begin{align} \label{part e-0+1}
 \Big|\int_{\Z^d} v_n 1_{B_R(0)^c} \bK_{p,\beta}(v_n) dx \Big|  \leq C \epsilon.
 \end{align}

Under the assumption $(\bA_{\alpha,\beta,2})$, we have that
 \begin{align*}
\Big|\int_{\Z^d} v_n 1_{B_R(0)^c} \bK_{p,\beta}(v_n) dx \Big| & =   \Big| \int_{\Z^d}\int_{\Z^d}(Q^{\frac1{p}}v_n)(x) (Q^{\frac1{p}}v_n)(y)\Phi_{d,\beta}(x,y)1_{B_R(0)^c}(x) dxdy \Big|
 \\[1mm]&\leq
  \int_{B_R(0)^c} \int_{\Z^d}|v_n(x)|  |v_n(y)|\big( Q(x)^{\frac1{p}}Q(y)^{\frac1{p}} \Phi_{d,\beta}(x,y) \big) dxdy
   \\[1mm] &\leq \big(\sup_{z\in \bQ_{R}(0)^c} Q_1(z)\big) \|Q_1\|_{L^\infty(\Z^d)} \int_{\Z^d}\int_{\Z^d}Q_{\alpha_1}^{\frac1p}(x)Q_{\alpha_1}^{\frac1p}(y) |v_n(x)|   |v_n(y)|  \Phi_{d,\beta}(x,y)   dxdy
    \\[1mm] &\leq c\big(\sup_{z\in \bQ_ R(0)^c} Q_1(z)\big) \|Q_1\|_{L^\infty(\Z^d)} \|v_n\|_{L^{p'}(\Z^d)}  \|v_n\|_{L^{p'}(\Z^d)}.
 \end{align*}
By (\ref{ass-3.1}),  for any $\epsilon>0$, there exists an integer, still denoted by $R_\epsilon>1$, such that  for $R\geq R_\epsilon$
$$\sup_{z\in \bQ_{\frac R2}(0)^c} Q_1(z)\leq\epsilon, $$
which implies that for $R=R_\epsilon$ and some $C>0$
  \begin{align} \label{part e-0+10}
 \Big|\int_{\Z^d} v_n 1_{B_R(0)^c} \bK_{p,\beta}(v_n) dx \Big|  \leq C \epsilon.
 \end{align}

{\it Part II:} \ For $R=R_\epsilon$,
 we obtain that
  \begin{align*}
 \Big|\int_{\Z^d} v_n 1_{B_R(0)} \bK_{p,\beta}(v_n) dx \Big| &= \Big| \int_{\Z^d}\int_{\Z^d}(Q^{\frac1{p}}v_n)(x) (Q^{\frac1{p}}v_n)(y)\Phi_{d,\beta}(x,y)1_{\bQ_R(0)}(x) dxdy \Big|
    \\[1mm] &\leq C\|Q\|_{L^\infty(\Z^d)}^{\frac{2}{p}} R^{\frac{2d}{p}} \Big(\int_{\Z^d} |v_n(x)|^{p'}1_{\bQ_R(0)}(x) dx\Big)^{\frac1{p'}} \Big(\int_{\Z^d} |v_n(y)|^{p'}  dy\Big)^{\frac1{p'}} ,
        \end{align*}
 then by (\ref{eqs-3.1}),  there exists an integer $n_R>0$ such that for $n\geq n_R$
  \begin{align} \label{part e-1+1}
 R^{\frac{2d}{p}} \Big(\int_{\Z^d} |v_n(x)|^{p'}1_{\bQ_R(0)}(x) dx\Big)^{\frac1{p'}}    \leq \epsilon,
 \end{align}
 which,  together with (\ref{part e-0}),  implies that for any $\epsilon>0$,  there is $n_\epsilon>0$ such that
  \begin{align*}
   \Big|   \int_{\Z^d} v_n\bK_{p,\beta}(v_n) dx   \Big| \leq C\epsilon \quad {\rm for}\ \ n\geq n_\epsilon.
 \end{align*}
 Thus, we obtain (\ref{eqs-3.1}). \hfill$\Box$\bigskip

\noindent {\bf Proof of Proposition \ref{pr 3.1}. } $(i)$ Since $p'\in(1,2)$ for $p>2$, it follows by Lemma \ref{lm 2.4} that for $\|v_n\|_{L^{p'}(\Z^d)} = \rho$,
 \begin{align*}
 \cJ_0(v)&=\frac1{p'} \rho^{p'}-\frac12\int_{\Z^d} v\bK_{p,\beta}(v) dx
 \\[1mm]&\geq \frac1{p'} \rho^{p'}-c \frac{\rho^2}2
  \\[1mm]&\geq \frac1{2p'} \rho^{p'}\quad {\rm for}\  \rho>0\ {\rm small\ enough}.
  \end{align*}

 $(ii)$ Take $v_t=t\delta_{x_0}$, where $x_0\in\Z^d$ such that  $\Phi_{d,\beta}(x_0,x_0)>0$, then
 \begin{align*}
 \cJ_0(v_t)&=\frac1{p'}t^{p'}  -\frac12   Q(0)^{\frac2p} \Phi_{d,\beta}(x_0,x_0)t^2<0 \quad
 {\rm if}\  t>1\ {\rm large\ enough}.
  \end{align*}

 $(iii)$ Let $(v_n)_n$ be a Palais-Smale sequence, i.e. there holds
 $\sup_n|\cJ_0(v_n)|<+\infty$ and
 $\cJ_0'(v_n)\to 0$ in $(L^{p'}(\Z^d))'=L^p(\Z^d)$ as $n\to+\infty$. We have
 \begin{align*}
+\infty> \sup_n|\cJ_0(v_n)|&\geq \cJ_0(v_n)
 \\[1mm]&=(\frac1{p'}-\frac12) \|v_n\|_{L^{p'}(\Z^d)}^{p'}+\frac12\cJ_0'(v_n)v_n
  \\[1mm]&\geq (\frac1{p'}-\frac12) \|v_n\|_{L^{p'}(\Z^d)}^{p'}- \frac12\|\cJ_0'(v_n)\|_{L^p(\Z^d)} \|v_n\|_{L^{p'}(\Z^d)}
  \\[1mm]&\geq (\frac1{p'}-\frac12-\epsilon) \|v_n\|_{L^{p'}(\Z^d)}^{p'}- \frac1{2\epsilon} \|\cJ_0'(v_n)\|_{L^p(\Z^d)}^p,
   \end{align*}
   where $\|\cJ_0'(v_n)\|_{L^p(\Z^d)}\to0$  as $n\to+\infty $.
 Then  $\|v_n\|_{L^{p'}(\Z^d)}$ is uniformly bounded.\smallskip

 Now we set the sequence $(v_n)_n$ in $L^{p'}(\Z^d)$ satisfying that
  $$\cJ_0(v_n)\to c\in\R\setminus\{0\}, \quad\cJ_0'(v_n)\to0\  {\rm in}\ L^{p}(\Z^d)\quad {\rm as}\ n\to+\infty, $$
 then
 \begin{align*}
(\frac1{p'}-\frac12)\int_{\Z^d}   v_n\bK_{p,\beta}( v_n) dx  = \cJ_0(v_n)-\frac1{p'}  \cJ_0'(v_n)v_n
\to c\quad {\rm as}\ n\to+\infty,
    \end{align*}
  and  there exists $n_0>1$ such that  for $n\geq n_0$
     \begin{align*}
 \int_{\Z^d} v_n\bK_{p,\beta}(v_n) dx  \not=0.
    \end{align*}

 Now we apply  Lemma \ref{lm 3.1} to obtain that, letting $\tilde v_n=v_n $ in $\Z^d$,  for some $R>1,\ \epsilon_0>0$
 \begin{align*}
 \int_{B_R(0)} |\tilde v_n|^{p'} dx\geq \epsilon_0\quad {\rm for\ all}\ n\geq n_0.
    \end{align*}
 Hence, up to a subsequence, we may assume $\tilde  v_n\rightharpoonup v\in L^{p'}(\Z^d)\setminus\{0\}$ as $n\to+\infty$. From the convexity of the function $t\mapsto |t|^{p'}$ and Lemma \ref{cr 3.1}, we obtain that
  \begin{align*}
 \frac1{p'} \|v\|_{L^{p'}(\Z^d)}^{p'} - \frac1{p'} \|\tilde v_n\|_{L^{p'}(\Z^d)}^{p'} &\geq \int_{\Z^d} |\tilde v_n|^{p'-2}\tilde v_n(v-\tilde v_n)
  \\[1mm]&=  \cJ_0'(\tilde v_n) (v-\tilde v_n) +\int_{\Z^d}  \tilde  v_n \bK_{p,\beta}  (v-\tilde v_n) dx
  \; \to0\quad {\rm as}\ \, n\to+\infty,
   \end{align*}
 then
  \begin{align*}
  \|v\|_{L^{p'}(\Z^d)}\geq \limsup_{n\to+\infty}   \|\tilde  v_n\|_{L^{p'}(\Z^d)}.
   \end{align*}
 Together with $\tilde  v_n\rightharpoonup v\in L^{p'}(\Z^d)\setminus\{0\}$, we derive that
 $$\tilde  v_n\to v\in L^{p'}(\Z^d) \quad {\rm as}\ n\to+\infty. $$
We complete the proof.   \hfill$\Box$\medskip

 \noindent{\bf  Proof of Theorem \ref{teo 3.1}. }
 %We employ the Mountain Pass Theorem to obtain the weak solution of (\ref{eq 1.1})
 %by considering  the associated  energy functional $\cJ_0\in C^1(L^{p'}(\Z^d),\R)$ defined by (\ref{eqq-3.3}).
%\begin{align*}
%\cJ_0(v)=\frac1{p'}\int_{\Z^d} |v|^{p'} dx-\frac12 \int_{\Z^d}  v\bK_{p,\beta}(v)  dx\quad{\rm for}\ \,  v\in L^{p'}(\Z^d).
%\end{align*}
We consider the critical level
$${\bf c}:=\inf_{\gamma\in\Gamma}\max_{t\in[0,1]} \cJ_0(\gamma(t)) ,$$
 where the functional $\cJ_0\in C^1(L^{p'}(\Z^d),\R)$ is defined by (\ref{eqq-3.3}) and
 $$ \Gamma=\{\gamma\in C([0,1],L^{p'}(\Z^d)):\, \gamma(0)=0,\ \cJ_0(\gamma(1))<0 \}. $$
From Proposition \ref{pr 3.1}, ${\bf c}>0$ and we may use Mountain Pass Theorem  to obtain that  there exists a point $v\in L^{p'}(\Z^d)$ achieving the   critical level ${\bf c}$
 and it verifies the equation
$$
|v|^{p'-2}v=Q^{\frac1{p}} \Phi_{d,\beta}\ast (Q^{\frac1{p}}v) \quad {\rm in}\ \, \Z^d.
$$
Since $Q,\Phi_{d,\beta}$ are nonnegative, then
$$ \int_{\Z^d}  |v|\bK_{p,\beta}(|v|)  dx\geq \int_{\Z^d}  v\bK_{p,\beta}(v)  dx $$
%{\color{blue}  
and  $\cJ_0(|v|)\leq \cJ_0(v)$ for $v\in L^{p'}(\Z^d)$.  Obviously, $\cJ_0(-v)=\cJ_0(v)$,  so if $v$ is critical point, then $|v|$ is also a critical point, so we can assume that $v$ doesn't change signs and  set $v\gneqq0$. \smallskip

  By Lemma \ref{lm 2.0} part $(i)$,  we have that $v(x)\to0$ as $|x|\to+\infty$ thanks to $v\in L^{p'}(\Z^d)$. \hfill$\Box$\medskip

\subsection{Linear case: $p=2$}
  For $p=2$, we have that $p'=2$ and  (\ref{eqq-3.2}) reduces to a linear model. To this end, we consider
  the solution $(\lambda,u)$ of a modified linear problem
    \begin{align}\label{eqq-3.2-l}
 v=\lambda \bK_{2,\beta}(v)\quad {\rm in}\ \, \Z^d.
\end{align}
Recall that 
$$%\begin{equation}\label{ei-0-4}
\cM_0:=\big\{x\in\Z^d:\, \exists y\in\Z^d,\,  Q(x)Q(y)>0,\, \Phi_{d,\beta}(x,y)>0\big\}, \qquad  m_0=|\cM_0|\in[1, +\infty].
$$%  \end{equation}

%In this subsection, we assume that  $\Phi_{d,\beta}$ is nonnegative.   

\begin{theorem}\label{teo 3.2}
  Assume  that $d\geq 1$, $\beta\in(0,\frac d2)$,  $\Phi_{d,\beta}$ satisfies $(\bF_1)$ and $(\bF_2)$.
  % and 
%$$ \Phi_{d,\beta}(y,x)=0 \quad{\rm if}\ \, \Phi_{d,\beta}(x,x)\Phi_{d,\beta}(y,y)=0.  $$
Let $\alpha > 0$ be such that
  $$2^*_{\beta,\alpha}< 2$$
  or
   $$2^*_{\beta,\alpha}= 2\quad {\rm and}   \quad \lim_{|x|\to+\infty} Q(x)|x|^{\alpha}=0. $$

  %If assume more that $m_0\in(0,+\infty]$,  
  Then   
   \begin{itemize}
\item[(i)] 
problem \eqref{eqq-3.2-l} admits   a  solution $(\lambda_1,\phi_1)\in (0,+\infty)\times L^2(\Z^d)$ with $\phi_1\geq0$, 
 where
 $$\lambda_1=\Big(\sup_{\|v\|_{L^2(\Z^d)}=1}\int_{\Z^d}v \bK_{2,\beta}(v) dx\Big)^{-1}>0. $$

Moreover,  problem \eqref{eqq-3.2-l} admits a sequence of the  eigenvalue $\lambda_{k}$ for any integer $k\in [2,m_0]\cap (0,+\infty)$ if $m_0\in[2,+\infty]$, 
 with corresponding eigenfunctions $\phi_{k}$
 which can be characterized as
\begin{equation}\label{Lambda-k-s}
\lambda_{k} =\Big( \sup_{v\in \cP_{k} }\int_{\Z^d}v \bK_{2,\beta}(v) dx\Big)^{-1},
\end{equation}
where 
\[
\cP_{k} := \Big\{ u\in  L^2(\Z^d): \int_{\Z^d}u\phi_{j}dx  =0 \ \text{ for \ } j=1,2,\cdots k-1 \ \text{ and } \ \|u\|_{L^2(\Z^d)}=1\Big\}.
\]

\item[(ii)] For $k\in[2,m_0]\cap (0,+\infty)$ if $m_0\geq 2$,    
\[
 \lambda_{k-1} \le \lambda_{k}  
\]
 and
$$
  \lim_{k\to \infty}\lambda_{k} = +\infty\quad  {\rm if}\ \, m_0=+\infty.
$$

Moreover, for any $k\in\Z$, 
$$\lim_{|x|\to+\infty}\phi_k(x) = 0.  $$

\item[(iii)] There hold %The sequence of eigenfunctions  $\{\phi_{k}\}_{k\in\N}$ corresponding to eigenvalues $\{\lambda_{k}(\Z^d)\}$ form a complete orthonormal basis of $L^2(\Z^d)$ and 
$$ \int_{\Z^d}\phi_k \phi_j dx=0\quad{\rm and}\quad   \int_{\Z^d}\phi_k \bK_{2,\beta}(\phi_j) dx=0\quad {\rm for}\ \, k\not=j.  $$

\end{itemize}

    \end{theorem}

\noindent  {\bf Proof. }  It is known that $L^2(\Z^d)$ is a Hilbert space with the inner product  $ \langle \cdot,\cdot\rangle$ given by
$$\langle u,v\rangle:=\int_{\Z^d} u(x)v(x)dx   $$
and denote 
$$\cI_\beta(v)=\int_{\Z^d}v \bK_{2,\beta}(v) dx,\quad v\in L^2(\Z^d). $$
The eigenvalue problem is equivalent to solve 
 \begin{align}\label{eqq-3.2-l-rep}
  \bK_{2,\beta}(v)=\mu v \quad {\rm in}\ \, \Z^d.
\end{align}

Note that
$$
  \bK_{2,\beta}(v) =   Q^{\frac1{2}} \Phi_{d,\beta}\ast (Q^{\frac1{2}}v)  \quad{\rm for} \ \, v\in L^2(\Z^d).
  $$
 We need to prove that $\bK_{2,\beta}:L^2(\Z^d) \to L^2(\Z^d)$ is a self-adjoint compact operator.

Under the assumptions of Theorem \ref{teo 3.2},  (\ref{ei-1}) with $p=2$ leads to
$$
|\cI_\beta(v)|= \Big|\int_{\Z^d} v\bK_{2,\beta}(v) dx\Big| \leq c\|v\|_{L^2(\Z^d)}^2\quad{\rm for} \ \, v\in L^2(\Z^d) .
$$
Obviously, we have that
 \begin{align*}
  \big\langle u, \bK_{2,\beta}(v) \big\rangle &=  \big\langle  \bK_{2,\beta}(u),v \big\rangle=\big\langle v, \bK_{2,\beta}(u) \big\rangle.
\end{align*}
 Now it follows  by  Lemma \ref{cr 3.1}   that
  $\bK_{2,\beta}\!: L^2(\Z^d)\to L^2(\Z^d)$ is compact. 
Then
  $$\mu_1:=\sup_{\|v\|_{L^2(\Z^d)}=1}\int_{\Z^d}v \bK_{2,\beta}(v) dx>0$$
  could be achieved by some $\phi_1\in L^2(\Z^d)$ with ${\rm supp}(\phi_1)\cap \cM_0\not=\emptyset$.
 Since $ \Phi_{d,\beta}>0$ and $Q\geq0$, we obtain that
 $$ \int_{\Z^d}|v| \bK_{2,\beta}(|v|) dx\geq \int_{\Z^d}v \bK_{2,\beta}(v) dx.$$
  So  we can assume $\phi_1\geq0$.
  
  If $m_0\in[2,+\infty]$,    inductively,  we  can  assume that $\phi_{1},\dots,\phi_{k} \in  L^2(\Z^d)$ and $\mu_{1} \le \dots \le \mu_{k} $ are already given for some integer $k \in[2,m_0]$  with the properties that for $j=1,\dots,k$, the function $\phi_{j}$ is a maximizer
of $\cI_\beta$ within the set
 \begin{align*}
\cP_{j}:= \{u\in  L^2(\Z^d) \::\:  \norm{u}_{ L^2(\Z^d)}=1,\: {\rm supp}(u)\subset \cM_0, \  \text{$  \int_{\Omega} u\phi_{n} \,dx  =0$ for $n=1,\dots j-1$} \},
 \end{align*}
$$ \mu_{j}  = \sup_{w\in \cP_{m,j}}  \cI_\beta(w) =\cI_\beta (\phi_{j})$$and
\begin{equation}
  \label{eq:inductive-eigenvalue}
  \cI_{\beta}(\phi_{j}, \varphi)=\mu_{j}   \int_{\Z^d} \phi_{j} \varphi  dx   \qquad \text{for all $\varphi\in L^2(\Z^d)$.}
\end{equation}
Since $\cP_{j}\subset \cP_{j-1}$, 
then 
$$\mu_{j}\geq \mu_{j+1}.$$
 
 By (\ref{ei-0-2}),  we have that 
 $$\mu_j=\int_{\Z^d} \phi_j \bK_{2,\beta}( \phi_j) dx\geq0.  $$ 
 If $\mu_j=0$, then 
 $$\int_{\Z^d} (\phi_jQ^{\frac12}) \Phi_{d,\beta}\ast (\phi_jQ^{\frac12})  dx=0$$
 and 
 $$\phi_jQ^{\frac12}=0\quad {\rm in\ } \{z\in\Z^d: \Phi_{d,\beta}(z,z)>0\}.  $$
 and $ \phi_j=0 $ in $\cM_0$, which is impossible. So $\mu_j>0$. 
 
Moreover,  when $m_0<+\infty$, then the eigenvalue problem is equivalent to that, 
letting $\bT=(t_1, \cdots, t_{m_0})$ and 
$$u(x)=\sum_{i=1}^{m_0}t_i \delta_{p_i},\quad p_i\in\cM_0, $$
 equation (\ref{eqq-3.2-l-rep}) is equivalent 
 to 
 $$\mu \bT= A\bT,\quad \text{for $i=1,\cdots, m_0$} $$ 
 where the Matrix 
 $A=(Q^{\frac12}(p_i)Q^{\frac12}(p_j) \Phi_{d,\beta}(p_i,p_j))_{m_0\times m_0}$ has $m_0$ many eigenvalues and 
 it follows by (\ref{ei-0-3}) that all eigenvalues are positive.

  Let 
  $$\cI_{\beta}(u, v  ) =\int_{\Z^d}u \bK_{2,\beta}(v) dx\quad u,v\in L^2(\Z^d).  $$
  Applying the compact embedding again, we obtain that  the value $\mu_{k+1}  $ is attained by a function {$\phi_{k+1} \in \cP_{m, k+1}$.} Thus, there exists a Lagrange multiplier $ \mu= \mu_{k+1} \in \R$ such that
\begin{equation}
  \label{eq:inductive-eigenvalue-k+1}
\cI_{\beta}(\phi_{k+1}, \varphi)=\mu  \int_{\Z^d} \phi_{k+1} \varphi  dx   \qquad \text{for all $\varphi \in \cP_{k+1} $}
\end{equation}
and by choosing $\varphi = \phi_{k+1}$.  Moreover, for $j=1,\dots,k$, we have, by (\ref{eq:inductive-eigenvalue}) and the definition of $\cP_{k+1}$, that
\begin{align*}
\cI_{\beta}(\phi_{k+1}, \phi_{j})&=\cI_{\beta}(\phi_{j}, \phi_{k+1}) 
\\&= \mu_{j}   \int_{\Z^d} \phi_{j} \phi_{ k+1} \,dx  
 \\& = 0
\\&= \mu_{k+1}  \int_{\Z^d} \phi_{j} \phi_{ k+1} \,dx.
\end{align*}
Hence (\ref{eq:inductive-eigenvalue-k+1}) holds with $\mu=  \mu_{k+1} $ for all $\varphi\in L^2(\Z^d) $.
Inductively, we have now constructed a normalized sequence $(\phi_{k})_{k\in\N}$ in $L^2(\Z^d)$ and a non-increasing sequence $\{\mu_{k}  \}_{k\in\N}$ in $\R$ such that $\phi_k$ is an eigenfunction of (\ref{eq 1.1}) corresponding to $\mu =  \mu_{k}  $ for every $k \in \N$. %Moreover, by construction, the sequence $\{\phi_k\}_k$ forms an orthonormal system of $L^2(\Z^d)$.  

{\it We now show  $\lim \limits_{k\to+\infty} \mu_{k}  =0$ if $m_0=+\infty$. }  Assume by contradiction that $c:= \lim \limits_{k \to \infty}\mu_{k} >0$, then we deduce that
$\cI_{\beta} (\phi_{k},\phi_{k})_m\ge c$ for every $k \in \N$. Hence the sequence $(\phi_{k})_{k\in \N}$ is bounded in $L^2(\Z^d) $, and therefore by Rellich compactness theorem, $(\phi_{k})_k$
contains a convergent subsequence $(\phi_{k_j})_j$ in $L^2(\Z^d)$. However, this   is impossible since the functions $\{\phi_{k_j}\}_{j \in \N}$ are $L^2(\Z^d)$-orthogonal.  
  
  Finally, all properties of $\lambda_k$ can be derived by equality that  $\lambda_k=\mu_k^{-1}$. 
   \hfill$\Box$\medskip

\noindent{\bf Proof of Theorem \ref{teo 1-eig}. }
$(i)$ For given $\Omega\subset \Z^d$, if $\Phi_{d,\Omega}$ is the fundamental solution of $-\Delta_{\Z^d}$ under the zero Dirichlet boundary condition, i.e.
$$ \left\{%\arraycolsep=1pt
\begin{array}{lll}
 -\Delta_{\Z^d}  \Phi_{d,\Omega}(\cdot,y)=\delta_{y}  \qquad
 {\rm in}\ \  \Omega , \\[2mm]
 \phantom{ \quad\quad \   }
 \displaystyle \Phi_{d,\beta}(\cdot,y)=0  \ \,   \qquad {\rm on}\ \,   \Z^d\setminus \Omega\ \  \text{ if  $\Omega\not=\Z^d$},
 %\\[2mm]
% \phantom{   }
%\displaystyle \lim_{x\in\Omega,|x|\to+\infty}u(x)=0 \ \, \text{ if  $\Omega$ is not finite}, 
 \end{array}
 \right.
$$
where $y\in \Omega$,  then  eigenvalue problem
 \begin{equation}\label{eq 1.1-eign}
 \left\{%\arraycolsep=1pt
\begin{array}{lll}
 -\Delta_{\Z^d}  u=\lambda Q  u  \quad
 &{\rm in}\ \  \Omega , \\[2mm]
 \phantom{ \quad\quad \   }
 \displaystyle u=0    &{\rm on}\ \,  \Z^d\setminus \Omega\ \ \text{ if  $\Omega\not=\Z^d$},
 \end{array}
 \right.
\end{equation}
 turns to the  problem (\ref{eqq-3.2-l}) with $\beta=1$ and it follows by Theorem \ref{teo 3.2}. \hfill$\Box$

   \subsection{Sub-linear case: $p\in(1,2)$}

  For $p\in(1,2)$,   we consider
  the positive solution $u$ of the  sublinear problem
    \begin{align}\label{eqq-3.2-subl}
 u=\Phi_{d,\beta}\ast (Q |u|^{p-2}u)\quad {\rm in}\ \, \Z^d.
\end{align}

  \begin{theorem}\label{teo 3.3}
  Assume  that   $d\geq 1$, $\beta\in(0,\frac d2)$,  $\Phi_{d,\beta}$ satisfies $(\bF_1), (\bF_2)$ and 
   $\alpha\in\R$ and $p\in(1,2)$.
  If there exists $\bar u\gneqq0$ in $\Z^d$ such that
  $$\bar u\geq  \Phi_{d,\beta}\ast   (Q \bar u^{p-1})  \quad {\rm in}\ \, \Z^d. $$
   Then    problem (\ref{eqq-3.2-subl})
 has a positive solution  $u$ with $\displaystyle \lim_{|x|\to+\infty}u(x) = 0.  $
    \end{theorem}
\noindent{\bf Proof. }  
%{\it Existence: }
Let $x_0\in\Z^d$ satisfy
 $$\bar u(x_0)>0, \ \ \Phi_{d,\beta}(x_0, x_0)>0\quad\text{ and}\quad Q(x_0)>0,  $$
then
$$\bar u(x)\geq   \big(Q(x_0) \bar u(x_0)^{p-1}\big) \Phi_{d,\beta}(x,x_0)   \qquad {\rm for}\ \, x\in\Z^d. $$

We construct a subsolution $\underbar u\leq \bar u$ in $\Z^d$.
Let
$$w_t(x)=t\Phi_{d,\beta}(x, x_0)\quad {\rm in}\ \Z^d.  $$
Clearly for $t\in(0,t_1]$ with $t_1:= Q(x_0) \bar u(x_0)^{p-1}$,
 $$w_t(x)\leq \bar u(x) \qquad {\rm for}\ \, x\in\Z^d. $$
Note that
$$ \Phi_{d,\beta}\ast   (Q w_t^{p-1})  \leq \Phi_{d,\beta}\ast   (Q \bar u^{p-1}) \leq \bar u\quad\ {\rm in}\ \, \Z^d.$$
  and
   \begin{align*}
 w_t  -\Phi_{d,\beta}\ast   (Q w_t^{p-1})  &\leq   t\Phi_{d,\beta}(\cdot, x_0) -t^{p-1} \Phi_{d,\beta}(x_0, x_0)\Phi_{d,\beta}(\cdot, x_0)
 \leq 0
\end{align*}
  if $t>0$ small enough since $p<2$.
  Thus there is $t_2 \in (0, t_1]$ such that
  $$w_{t_2}\leq \Phi_{d,\beta}\ast   (Q w_{t_2}^{p-1})\quad {\rm in}\ \Z^d. $$
  That is $\underbar u := w_{t_2}$ is a subsolution.
  Now we set that $u_0=\underbar u$ and
  $$u_n=\Phi_{d,\beta}\ast   (Q u_{n-1}^{p-1})\quad {\rm in}\ \Z^d, $$
  then the mapping $n\to u_n$ is nondecreasing and bounded by $\bar u$.
  Therefore, there exists $u\in C(\Z^d)$ such that
  $$\underbar u \leq  u\leq \bar u\qquad {\rm in}\ \Z^d,$$
  $$\lim_{n\to+\infty}u_n(x)=u(x) \qquad {\rm for}\ \, x\in\Z^d$$
 and
  $$\lim_{n\to+\infty} \Phi_{d,\beta}\ast   (Q u_{n-1}^{p-1})=\Phi_{d,\beta}\ast   (Q u^{p-1})\quad\ \  {\rm in}\ \Z^d.$$
  So $u$ is a solution of (\ref{eqq-3.2-subl}) and this ends the proof. \smallskip

 \medskip

      \setcounter{equation}{0}
   \section{ In whole space $\Z^d$ }

\subsection{Existence  }

To show the existence in sub-linear case, we need the following lemmas.

  \begin{lemma}\label{pr 5.1-w}
Let   $g_\tau\in C (\Z^d)$ with $\tau\in(2-d,0)$    be a nonnegative function such that
\begin{align}\label{g bound1w}
\frac1{c_0}(1+ |x|)^{\tau-2} \leq  g_\tau(x)\leq c_0 (1+ |x|)^{\tau-2},\quad \ \forall \, x\in \Z^d \setminus B_{n_0}
\end{align}
    for   some $n_0>0$   and  $c_0\geq 1$.
Then the Poisson problem
 \begin{equation}\label{eq 3.1-poissonw}
\left\{%\arraycolsep=1pt
\begin{array}{lll}
\ -\Delta u     =  g_\tau   \quad\
    {\rm in}\ \,   \Z^d , \\[2mm]
 \phantom{   }
\displaystyle \lim_{|x|\to+\infty}u(x)= 0
 \end{array}
 \right.
\end{equation}
has a unique positive solution $v_\tau$ such that   for some $c\geq 1$,
 \begin{equation}\label{eq 2.1-homxxx}
 \frac1c (1+ |x|)^{\tau} \leq v_\tau (x)\leq c (1+  |x|)^{\tau}\quad \ \forall \, x\in \Z^d .
 \end{equation}
\end{lemma}

\noindent{\bf Proof. } For simplicity, we write $g=g_\tau$.    Let
$$ v_g(x)=(\Phi_d\ast g)(x)\quad {\rm for}\ \, x\in\Z^d,  $$
which is well-defined by (\ref{g bound1w}),  and is a solution of  (\ref{eq 3.1-poissonw}).
Obviously, $v_g$ is positive. We can   define
$$v_n=\Phi_d\ast g_n\quad {\rm in}\  \Z^d,  $$
where $g_n=g\chi_{B_n}$. Direct computation shows that
$$v_n\to v_g\quad {\rm locally}\; {\rm in} \ \Z^d\ \  {\rm as} \ \, n\to+\infty.   $$

Recall that for $\tau<0$ and
 $$  \bar v_{\tau}(x):= (1+ |x|)^{\tau}\quad{\rm for}\  x\in \Z^d\setminus B_2(0),$$
with $|x|$ large enough,
\begin{align}\label{coll-1}
  \Delta_x   \bar v_\tau(x)& =  \tau(d-2+\tau)  |x|^{\tau-2}+    O( |x|^{\tau-3}).
 \end{align}
Thus there is $n_0\geq 1$ such that
 $$ \frac1c  |x|^{\tau-2} \leq -\Delta_x   \bar v_\tau(x) \leq c  |x|^{\tau-2}\quad {\rm for}\ x\in\Z^d\setminus B_{n_0}. $$
 Observe that  for some $t_0>1$
$$\frac1{t_0}n_0^{2-d} \leq  v_g(x)\leq t_0n_0^{2-d} \quad {\rm for}\ \, x\in \Z^d,\  n_0-1\leq |x|\leq n_0+1.  $$
It follows by the comparison principle that
$$   v_n(x)\leq t_0  \bar v_{\tau}(x)\quad \  \forall \, x\in \Z^d\setminus B_{n_0}. $$
So is $v_g$. Again by applying the comparison principle, we can get that for some suitable $t_0>1$
$$\frac1{t_0} \bar v_{\tau}(x)\leq  v_g(x) \quad \  \forall \, x\in \Z^d\setminus B_{n_0}. $$
We complete the proof.   \hfill$\Box$\medskip

 \begin{lemma}\label{lm 4.2-mm}

Let $g_\sigma\in C(\Z^d)$ with $\sigma>0$ satisfy
$$\frac1c (1+|x|)^{ -d}\big(\ln (e+|x|^2)\big)^{ \sigma-1} \leq g_\sigma(x)\leq c(1+|x|)^{ -d}\big(\ln (e+|x|^2)\big)^{ \sigma-1}\quad {\rm for\ }\,  \Z^d\setminus B_{n_0}$$
for some $c>1$ and $n_0>0$.
Then the Poisson problem
 \begin{equation}\label{eq 3.1-poissonw1}
\left\{%\arraycolsep=1pt
\begin{array}{lll}
-\Delta u     =  g_\sigma   \quad\
    {\rm in}\ \,   \Z^d , \\[2mm]
 \phantom{   }
\displaystyle \lim_{|x|\to+\infty}u(x)= 0
 \end{array}
 \right.
\end{equation}
has a unique positive solution $v_{\sigma}$ such that   for some $c\geq 1$
 \begin{equation}\label{eq 2.1-homxxx1}
\frac1c (e+|x|)^{ 2-d} \big(\ln (e+|x|)\big)^{  \sigma}  \leq v_{\sigma}(x)\leq c (e+|x|)^{ 2-d} \big(\ln (e+|x|)\big)^{  \sigma} \quad{\rm for}\ \forall \, x\in \Z^d .
 \end{equation}

\end{lemma}
{\bf Proof. }  The existence and uniqueness  are standard. We only need to show (\ref{eq 2.1-homxxx1}).

For $\sigma>0$, let $\varphi_{0, \sigma}\in C^2(\R_+)$ be
$$\varphi_{0, \sigma} (t):=(e+t)^{\frac12 (2-d)} \big(\ln (e+t)\big)^{  \sigma}\quad{\rm for}\ \ \forall\, t \in \R_+,$$
where $\R_+=[0,+\infty)$.
 Let also
 $$\psi_{0, \sigma}(x):=\varphi_{0, \sigma}(|x|^2), $$
then  the bound (\ref{eq 2.1-homxxx1})  is equivalent to that for    $  \sigma  >0$,  $r_0>1$ and $c>1$,
 \begin{equation}\label{eq 2.1-0-homxxx1}
\frac1c  |x| ^{ -d}\big(\ln (e+|x|^2)\big)^{ \sigma-1}\leq   -\Delta     \psi_{0, \sigma}(x)  \leq c  |x| ^{ -d}\big(\ln (e+|x|^2)\big)^{ \sigma-1}\quad \text{for $|x|>r_0$}.
 \end{equation}
Direct computation shows that
$$
\varphi_{0, \sigma} '(t)=\frac12(2-d) (e+t)^{ -\frac{d}2} \big(\ln (e+t)\big)^{ \sigma}+ \sigma (e+t)^{ -\frac{d}2}  \big(\ln (e+t)\big)^{ \sigma-1},$$
\begin{align*}
\varphi_{0, \sigma}''(t)&= (e+t)^{-\frac12d-1} \big(\ln (e+t)\big)^{\sigma}  \Big[ \frac14(2-d)(-d)+ \frac12\sigma(2-d)  \big(\ln (e+t)\big)^{-1}
 + \sigma(\sigma-1)  \big(\ln (e+t)\big)^{-2}\Big].
\end{align*}
Then for $x\in\Z^d$ with $|x|$ large, we have  that
\begin{align*}
 \Delta \psi_{0, \sigma}(x)   &= \sum_{y\sim x}\big(\psi_{0, \sigma}(y)-\psi_{0, \sigma}(x)\big)
 \\[1mm]&= \sum_{y\sim x} \bigg\{\Big[\frac{1}2(2-d) + \sigma    \big(\ln (e+|x|^2)\big)^{-1} \Big](e+|x|^2)^{-\frac d2} \big(\ln (e+|x|^2)\big)^{ \sigma } (|y|^2-|x|^2)
 \\&\qquad\quad\ +\frac12
 \Big( \frac14(2-d)(-d)+\frac12\sigma(2-d)  \big(\ln (e+|x|^2)\big)^{-1}
 + \sigma(\sigma-1)  \big(\ln (e+|x|^2)\big)^{-2}\Big)
   \\&\qquad\qquad\quad\cdot(e+|x|^2)^{-\frac{d+2}2} \big(\ln (e+|x|^2)\big)^{ \sigma}  (|y|^2-|x|^2)^2 \bigg\}\big(1+o(1)\big) \allowdisplaybreaks
  \\[1mm]&=  \bigg\{    \Big[ (2-d)d   +2 \sigma   \big(\ln (e+|x|^2)\big)^{-1} \Big](e+|x|^2)^{-\frac{d}2} \big(\ln (e+|x|^2)\big)^{ \sigma  }
 \\&\qquad +
   \Big( (2-d)(-d)+2 \sigma(1-d)  \big(\ln (e+|x|^2)\big)^{-1}
 + 4\sigma(\sigma-1)  \big(\ln (e+|x|^2)\big)^{-2}\Big)
   \\&\qquad\qquad\quad\cdot(e+|x|^2)^{-\frac{d}2} \big(\ln (e+|x|^2)\big)^{\sigma}      \bigg\}\big(1+o(1)\big).
\end{align*}
Thus for $|x|$ large enough, we have that
\begin{align*}
 -\Delta \psi_{0, \sigma}(x)  &=|x| ^{ -d}\big(\ln (e+|x|^2)\big)^{\sigma-1}  \Big( \beta_{1}(\sigma)
 +\beta_{2}(\sigma)  \big(\ln (e+|x|^2)\big)^{ -1}\Big) \big(1+o(1)\big),
\end{align*}
where
\begin{align}\label{ex-con-1}
\beta_{1}(\sigma)=  2\sigma  (d-2 )\quad{\rm and}\quad  \beta_{2}(\sigma)=-4\sigma(\sigma-1) .
\end{align}

For  $ \sigma  >0$ , then $\beta_{1}(\sigma)>0$  and
there exists $r_0>1$ such that for $|x|>r_0$,
$$\frac12\beta_{1}(\sigma)  |x| ^{ -d}\big(\ln (e+|x|^2)\big)^{ \sigma-1}\leq  -\Delta    \psi_{0, \sigma}(x)  \leq 2\beta_{1}(\sigma)  |x| ^{ -d}\big(\ln (e+|x|^2)\big)^{ \sigma-1}.  $$
The proof ends. \hfill$\Box$\medskip

  \noindent{\bf Proof of Theorem \ref{teo 1}. }  {\it Part $(i)$: }
    It is known that the fundamental solution  $\Phi_d(\cdot-y)$  of $-\Delta$ satisfies
 \begin{equation}\label{eq 2.1-F0}
\left\{%\arraycolsep=1pt
\begin{array}{lll}
  -\Delta u   = \delta_y \quad\ \,
 {\rm in}\ \  \Z^d , \\[2mm]
 \phantom{    }
 \displaystyle \lim_{|x|_{_Q}\to+\infty} u(x)=0,
 \end{array}
 \right.
\end{equation}
where   $y\in \Z^d$ and $\delta_y$ is the Dirac mass at $y$. When $d\geq 3$,  from \cite{LL,EG},    the fundamental solution  $\Phi_d$ has   the following asymptotic behaviors:
  \begin{equation}\label{fund-00}
\lim_{|x|_{_Q}\to+\infty}   \Phi_d (x,y)|x-y|^{d-2}=\varpi_d >0
  \end{equation}
  and
   \begin{equation}\label{fund-0}
   0<   \Phi_d (x,y) \leq  c_1(1+|x-y|)^{2-d}\quad\text{ in  \ $\Z^d $}.
     \end{equation}
%where    $ \varpi_d >0$.
Thus the original equation  (\ref{eq 1.1}) turns to the following integral equation
  \begin{align} \label{int eq-1-0-0}
  u=\Phi_{d}\ast (Q|u|^{p-2}u)\quad {\rm in}\ \, \Z^d.
  \end{align}
This is equivalent to considering the equation
 \begin{align} \label{int eq-1}
|v|^{p'-2}v=Q^{\frac1{p}} \Phi_{d}\ast (Q^{\frac1{p}}v) \quad {\rm in}\ \, \Z^d,
\end{align}
%where \begin{align*}
%v:=Q^{\frac1{p'}}|u|^{p-2}u\quad {\rm in}\ \, \Z^d,
%\end{align*}
%for $u\in L^{p}(\Z^d,Qdx)$.
%then $v$ satisfies
We employ Theorem \ref{teo 3.1} with  $\beta=1$,  $d\geq 3$
and
$\Phi_{d,1}=\Phi_{d}\quad {\rm in}\ \, \Z^d\times \Z^d $
to obtain that Eq.(\ref{int eq-1}) has a nonnegative nontrivial solution $v\in L^{p'}(\Z^d)$. In our setting,  we mention that
assumptions   $(\bA_{\alpha,\beta,1})$ and $(\bA_{\alpha,\beta,2})$ reduce to
$(A1)$ or $(A2)$ respectively.
Therefore
  $$u:= \Phi_d \ast (Q^{\frac{1}{p}}v)  \quad {\rm in}\ \Z^d $$
  is a solution of (\ref{int eq-1-0-0})
 %then
 %replacing $v=Q^{\frac1{p'}}|u|^{p-2}u$,
 %we obtain that
 %$$u=\Phi_{d}\ast (Q|u|^{p-2}u)\quad {\rm in}\ \, \Z^d   $$
 and
  $$ \int_{\Z^d}Q|u|^p dx=\int_{\Z^d} \big(Q^{\frac1{p'}}|u|^{p-1}\big)^{p'} dx=\int_{\Z^d} |v|^{p'} dx<+\infty,   $$
 which implies that
  $u\in L^{p}(\Z^d,Qdx)$ is a solution of (\ref{eq 1.1}), so that $u(x)\to 0$ as $|x|\to+\infty$ by  Lemma \ref{lm 2.0} part $(i)$.
 Moreover it follows by the strong maximum principle that $u>0$ in $\Z^d$.
\smallskip

%By Lemma \ref{lm 2.0} part $(i)$, we have that $u(x)\to0$ as $|x|\to+\infty$ thanks to $u\in L^{p}(\Z^d,Qdx)$.

\vskip2mm
{\it Part $(ii)$: } We first consider the case: $p-1\in  (0,1)$ with $\alpha>2$.  Let
$$\bar u_t=t \Phi_d\ast (1+|\cdot|)^{\tau_p-2} \quad{\rm in}\ \,  \Z^d, $$
where
$$\tau_p=\max\Big\{-\frac{\alpha-2}{2-p},\frac{2-d}{2}\Big\}\quad \text{ for $\alpha>2$}. $$
Then $\tau_p\in(2-d,0)$  and  $(p-1)\tau_p-\alpha \leq \tau_p-2$.
From Lemma \ref{pr 5.1-w} we have that
$$\frac1c t (1+|x|)^{\tau_p}\leq \bar u_t(x)\leq c t (1+|x|)^{\tau_p}\quad{\rm for}\ \, x\in \Z^d. $$
Note that for  $x\in\Z^d$,
\begin{align*}
Q(x) \bar u_t(x)^{p-1} &\leq C t^{p-1} (1+|x|)^{ (p-1)\tau_p-\alpha}
\\[1mm]&\leq C t^{p-1} (1+|x|)^{ \tau_p-2}
\\[1mm]&\leq C t^{p-2}( -\Delta \bar u_t),
\end{align*}
for all $t\geq t_1$, where $t_1\ge 1$ is taken such that $$C t_1^{p-2}\leq 1. $$
It follows by Theorem \ref{teo 3.3} that  problem (\ref{eq 1.1})
has  a unique positive solution $u$ such that for some $c>0$
$$0<u(x)\leq c t_1 (1+|x|)^{\tau_p}\quad {\rm for}\ x\in\Z^d. $$
We complete the proof.    \hfill$\Box$\medskip

\subsection{Nonexistence}
   This subsection  is devoted to the nonexistence of solution to (\ref{eq 1.1}).

\begin{proposition} \label{teo 2-wh}
 Assume that $d\geq 3$ and
 $$Q(x)\geq c(1+|x|)^{-\alpha},\quad \forall\, x\in\Z^d\setminus B_{n_0}$$
for some $c>0$,  $n_0>1$ and  $\alpha\in(-\infty,2)$.  Then  for $p-1\in(0, \frac{d-\alpha}{d-2}]$,  problem (\ref{eq 1.1}) has no positive solutions.

%There hold then \smallskip

%$(i)$ If  $p-1\in(0, \frac{d-\alpha}{d-2})$,   then
 % problem (\ref{eq 1.1}) has no positive solutions.

%$(ii)$ If $p-1=\frac{d-\alpha}{d-2}>1$, then  problem (\ref{eq 1.1}) has no positive solutions.
 \end{proposition}

 To show the nonexistence results, we need the following auxiliary lemmas.

  \begin{lemma}\label{pr 3.2w}
  Let $d\ge 3$ and nonnegative function $f\in C(\Z^d_+)$  verify that
  %for some $m\geq 1$
\begin{equation}\label{con 2.1-1w}
\lim_{n\to+\infty}\int_{   B_n(0)   }f(x) (1+ |x|)^{2-d} dx=+\infty.
\end{equation}
Then the homogeneous problem
\begin{equation}\label{eq 1.1 EH1}
\left\{%\arraycolsep=1pt
\begin{array}{lll}
-\Delta u    \geq  f   \quad
   &{\rm in}\ \  \Z^d  , \\[2mm]
 \phantom{ --  }
u\geq 0 \quad &{\rm   in}\ \    \Z^d
 \end{array}
 \right.
\end{equation}
 has no solutions.
 \end{lemma}
\noindent{\bf Proof.} We assume by contradiction that there exists a nonnegative solution $u_0$  of (\ref{eq 1.1 EH1}).
Then the strong maximum principle implies that $u_0>0$ in $\Z^d$.

 Let  $v_{n,f}$ be the minimal positive solution of
 \begin{equation}\label{eq 3.1-fn}
\left\{%\arraycolsep=1pt
\begin{array}{lll}
-\Delta u   = f_n   \quad
   {\rm in}\ \ \Z^d   , \\[2mm]
 \phantom{   }
\displaystyle \lim_{|x|\to+\infty}u(x)= 0 ,
 \end{array}
 \right.
\end{equation}
where  $f_n=f\chi_{ B_n(0)  }$. Here $\chi_{ B_n(0)  }$ is the indicator function of $ B_n(0)$.
%$\chi_{A}=1$ in $A$, $\chi_{A}=0$ otherwise.

By the comparison principle, we have that
$$0\leq v_{n,f}\leq u_0\quad{\rm in}\ \, \Z^d$$
and
$$v_{n,f}(x)=\sum_{z\in \Z^d}\Phi_{d}(x,z)f_n(z),\quad\forall\,  x\in\Z^d.  $$
There is some $c>1$ such that
$$\frac1c |x|^{2-d}\leq v_{n,f}(x)\leq c |x|^{2-d} \quad {\rm for}\ \, x\in\Z^d$$
and it follows by (\ref{con 2.1-1w}) and the comparison principle that  there exists $c>0$ such that for $n>4$
\begin{align*}
u_0(0)  \geq v_{n,f}(0)&=\int_{\Z^d }\Phi_{d}(0,z)f_n(z) dz
 \\& \geq c\int_{ B_{n}\setminus B_{4} }   |z|^{2-d}f(z) dz
\to+\infty\quad{\rm as}\ \ n\to+\infty,
\end{align*}
which is impossible. The nonexistence conclusion follows.  \hfill$\Box$\medskip

 \begin{lemma}\label{lm 2.1-ew}
 Let $d\geq 3$ and $\alpha<2$, $q\in(0,\frac{d-\alpha}{d-2})$
  and $\{\tau_j\}_j$ be a sequence defined by
  $$\tau_0=2-d<0,\qquad  \tau_{j+1}=  \tau_jq +2-\alpha,\quad j\in\N_+,$$
  where  $\N_+$ be the set of nonnegative integers.

 Then the map $j\in\N\to \tau_j$ is  strictly increasing and  for any $\bar \tau>\tau_0$ if $q\geq1$ or
 for any $\bar \tau\in(\tau_0, \frac{2-\alpha}{1-q})$ if $q\in(0,1)$,   there exists $j_0\in\N $  such that
 \begin{equation}\label{2.3-0}
 \tau_{j_0}\ge \bar \tau\quad {\rm and}\quad \tau_{j_0-1}<\bar \tau.
 \end{equation}
 \end{lemma}

 \noindent{\bf Proof. }
 First we have
 $$\tau_1-\tau_0= 2-\alpha+\tau_0(q-1)>0$$
 since $q\in(0,\frac{d-\alpha}{d-2})$, and by definition,
 Under the assumptions, we have that
 $$%\begin{equation}\label{2.1-increasing}
 2 -\alpha+\tau_0(q-1)>0,
 $$%\end{equation}
 then

%and
 \begin{equation}\label{2.2}
 \tau_j-\tau_{j-1} = q(\tau_{j-1}-\tau_{j-2})=q^{j-1} (\tau_1-\tau_0)>0.
 \end{equation}
 Then the sequence $\{\tau_j\}_j$ is strictly increasing. Moreover, if $q\ge1 $, the conclusion (\ref{2.3-0}) is straightforward.
 If $q\in(0,1)$,
   in the case when $\tau_1\ge0$, we are done, and in the case when $\tau_1<0$,
 it follows from (\ref{2.2}) that
 \begin{eqnarray*}
 \tau_j  &=&  \frac{1-q^j}{1-q}(\tau_1-\tau_0)+\tau_0\\
 &\to&\frac{1}{1-q}(\tau_1-\tau_0)+\tau_0=\frac{2 -\alpha}{1-q}\quad\  {\rm as}\ \, j\to+\infty,
 \end{eqnarray*}
 then there exists $j_0>0$ such that (\ref{2.3-0}) holds.
   \hfill$\Box$\medskip

\noindent{\bf Proof of Proposition \ref{teo 2-wh}. }    By contradiction, let $u_0\in C(\Z^d)$  be a nonnegative nonzero solution of  (\ref{eq 1.1}).   By the maximum principle, we obtain that
$$u_0>0\quad {\rm in}\ \, \Z^d. $$
% {\it  We claim that there exists $c_0>0$ and $n_0\geq 1$ such that
%$$  u_0(x)\geq c_0  (1+ |x|)^{2-d}\quad {\rm for\ all}\ \,  x\in \Z^d  .$$}
Moreover, from  the comparison principle, there exists $d_0>0$  and $n_0\geq 1$ such that
%\begin{equation}\label{sec 3-1.0-0}
$$
u_0(x)\geq  \frac{u_0(0)}{\Phi_{d}(0, 0)}\, \Phi_{d}(x, 0)\geq d_0 (1+ |x|)^{2-d}\quad {\rm for}\ \, x\in  \Z^d.
$$%\end{equation}
%Let $\tau_0=2-d<0$ satisfy that
Therefore
\begin{equation}\label{sec 3-1.0}
-\Delta u_0(x) =Q(x)u_0^{p-1}
%\geq  d_0^{p-1}   |x|^{-\alpha-q(d-2)}
\geq d_0^{p-1}    |x|^{\tau_1-2},
\quad\forall\, x\in \Z^d\setminus B_{n_0},
\end{equation}
where  $\tau_1$ is given by the previous lemma. Let $q:=p-1$, we recall that for $q \in (0, \frac{d-\alpha}{d-2})$, it holds that
$\tau_1 - \tau_0 = \tau_0 (q-1) +2-\alpha > 0$.
%$q=p-1$ and we have used the assumption on $Q(x)$.
% \quad  \tau_1=-q(d-2)-\alpha+2.$$
 Recall that $\alpha <2$, so we distinguish three cases.
 %for $q\in(0,  \frac{d-\alpha}{d-2})$, it holds that
%$$\tau_1-\tau_0=-q(d-2)-\alpha+2 +(d-2)>0.$$

\smallskip

{\bf Case 1}:    $q\in(0, \frac{2-\alpha}{d-2}]$.
%with $\alpha\in(-\infty,2)$.
Note that $q(2-d)-\alpha\geq  -2$, then
%\begin{equation}\label{sec 3-1.0-0}
%Q(x)  u_0(x)^{q}\geq  d_0^q (1+  |x|)^{(2-d)q-\alpha},\quad\forall x\in  \Z^d \setminus B_{n_0}
%\end{equation}
a contradiction follows by Lemma  \ref{pr 3.2w}  with $f(x)=d_0^q (1+|x|)^{q(2-d)-\alpha }$.   \smallskip

{\bf Case 2}:  $q\in\big(\frac{2-\alpha}{d-2},   \frac{d-\alpha}{d-2}\big)\subset  (0,+\infty)$.
%with $\alpha\in(-\infty,d)$.
 By Lemma \ref{pr 5.1-w}, there exists $d_1>0$   such that
$$u_0(x)\geq d_1(1+|x|)^{\tau_1},\quad\forall\, x\in   \Z^d, $$
with $\tau_1 \in(2-d, 0)$.  If $q\tau_1-\alpha \geq   -2$, we are done by Lemma  \ref{pr 3.2w}. Otherwise, we claim that the
iteration must be stopped after a finite number of times.  In fact, if $q \in [1, +\infty) \cap \big(\frac{2-\alpha}{d-2},   \frac{d-\alpha}{d-2}\big)$, since
$\tau_j\to+\infty$, then there exists $j_0\in\N$ such that
$q\tau_{j_0} -\alpha \geq  -2$, then a contradiction could be derived as in Case 1.  For $q\in(0,1)\cap \big(\frac{2-\alpha}{d-2},   \frac{d-\alpha}{d-2}\big)$, one has that  $\tau_j\to \tilde \tau_q:=\frac{2-\alpha}{1-q}>0$ as $j\to+\infty$,
then
$$-\alpha +q \tilde \tau_q>-2 $$
 if $q<\frac{d-\alpha}{d-2}$.
As a consequence,  there exists $j_0\in\N$ such that
 $q\tau_{j_0} -\alpha \leq  -2$ and  $q\tau_{j_0+1}-\alpha  \geq  -2$. This means we again get a contradiction and we are done. \smallskip

%Recall that
%$$\tau_{j+1}:=q\tau_{j} -\alpha+2,\quad\forall\,  j\in\N_+,$$
 %which is an increasing sequence.

 %If $\tau_{j+1}=\tau_jq +\alpha-2\in  (0,d-2)$,
%it follows by Proposition \ref{pr 5.1-w} that there exist integer  $d_{j}>0$  such that
%$$u_0(x)\geq d_{j}(1+ |x|)^{\tau_{j+1}}\quad {\rm in}\ \,  \Z^d. $$
 %  If $q\tau_{j+1}-\alpha \geq   -2$, we are done by Lemma  \ref{pr 3.2w}.\smallskip

 %Now we claim that the iteration must stop after a finite number of times. It infers by   Lemma \ref{lm 2.1-ew} that
% $j\mapsto \tau_j$ is strictly increasing thanks to $ 0<q<\frac{d-\alpha}{d-2}. $

%Note that for $q\in[1,+\infty)\cap \big(\frac{2-\alpha}{d-2},   \frac{d-\alpha}{d-2}\big)$, $\tau_j\to+\infty$, then there exists $j_0\in\N$ such that
%$q\tau_{j_0+1} \geq  -2$ and a contradiction could be derived for
%$1\leq q<\frac{d}{d-1}$.

% For $q\in(0,1)\cap \big(\frac{2-\alpha}{d-2},   \frac{d-\alpha}{d-2}\big)$, one has  $\tau_j\to \tilde \tau_q:=\frac{2-\alpha}{1-q}>0$ as $j\to+\infty$, since $\alpha<2$.   Then

%Thus,
%  there exists $j_0\in\N$ such that
 % $q\tau_{j_0} -\alpha \leq  -2$ and  $q\tau_{j_0+1}-\alpha  \geq  -2$. This means we can get a contradiction and we are done. \smallskip

{\bf Case 3}:   $q= \frac{d-\alpha}{d-2}>1$.
%with $\alpha\in(-\infty,2)$. }
%From (\ref{sec 3-1.0-0}),
In this case, we have in fact that
\begin{equation}\label{sec 3-2.0-0}
Q(x)  u_0(x)^{q}\geq  d_0^q (1+  |x|)^{-d},\quad\forall x\in  \Z^d \setminus B_{n_0}.
\end{equation}
From Lemma \ref{lm 4.2-mm} with $\sigma=1$, we have that
 \begin{equation}\label{sec 3-2.0-0}
  v_{\sigma}(x)\geq c (e+|x|)^{ 2-d}  \ln (e+|x|) \quad{\rm for}\ \forall \, x\in \Z^d,
  \end{equation}
  where $v_\sigma$ is the solution of (\ref{eq 3.1-poissonw1}).
  Now (\ref{sec 3-2.0-0}) implies that
 $$H_0(x):=Q(x)  u_0(x)^{q-1}\geq  d_0^q (1+  |x|)^{-2} \big(\ln (e+|x|)\big)^{q-1},\quad\forall x\in  \Z^d \setminus B_{n_0}.$$
Then there exists $n_0\geq 1$ such that 
$$H_0(x)\geq  \frac12 \bar H_0 u_0\quad {\rm in}\ \Z^d\setminus B_{n_0} $$
 Then  we can write
 $$-\Delta u_0 \geq \frac12 \bar H_0 u_0\quad {\rm in}\ \Z^d\setminus B_{n_0}.  $$
 On the other hand, from  (\ref{coll-1}), 
 \begin{align*} 
  -\Delta    \bar v_\tau(x)-\frac12 \bar H_0 \bar v_\tau(x)& = -\big( \tau(d-2+\tau)+c_0\big)  |x|^{\tau-2}+    O( |x|^{\tau-3})
 \end{align*}
where   $c_0=\frac12\inf_{x\in \Z^d}\frac{\bar H_0(x)}{(1+|x|)^{-2}}$ and $\bar v_\tau(x)=(1+|x|)^{\tau}$ for $x\in\Z^d$.  
 Then  there exists $\tau\in (2-d,0)$ and $n_1\geq n_0$ such that  
  \begin{align*} 
  -\Delta    \bar v_\tau(x)\leq \frac12 \bar H_0 \bar v_\tau(x) \quad \quad {\rm in}\ \Z^d\setminus B_{n_1}
 \end{align*}
If there exists $\sigma_1>0$ such that 
$$\sigma_1v_\tau\leq u_0\quad {\rm in}\ \, \bar B_{n_1}$$
then it follows by Lemma \ref{com lm-h1} with $\mu=-\frac12$
 $$u_0(x)\geq \sigma_1 v_\tau(x) {\rm in}\,\ \Z^d. $$
 Next, we can continue the iteration in {\bf case 2} to obtain a contradiction.    \hfill$\Box$ \medskip

\noindent {\bf Proof of Theorem \ref{teo 1}      Part $(iii)$. } It follows by Proposition \ref{teo 2-wh} directly.

    \setcounter{equation}{0}
   \section{ Lane-Emden equations in unbounded domains   }

\subsection{In half Space $\Z^d_+$}
  We consider the fundamental solution  of $-\Delta$  in the half space under the zero Dirichlet boundary condition, i.e.
 \begin{equation}\label{eq 2.1-F}
\left\{%\arraycolsep=1pt
\begin{array}{lll}
\quad   -\Delta u   = \delta_y \quad\ \,
  {\rm in}\ \  \Z^d_+ ,
 \\[2mm]
 \phantom{  --}
\quad u=0\quad\ \  {\rm on} \ \  \partial \Z^d_+,
 \\[2mm]
 \phantom{    }
 \displaystyle \lim_{x\in\Z^d_+,|x| \to+\infty} u(x)=0,
 \end{array}
 \right.
\end{equation}
where   $y\in \Z^d_+$ and $\delta_y$ is the Dirac mass at $y$. Then the existence and its asymptotic behaviors at infinite can be stated as follows.

\begin{proposition}\label{lm fund-1}
Let $d\geq 2$, 
then  (\ref{eq 2.1-F})  has a unique solution $\Phi_{d,+}$.

Furthermore,   $(i)$
$$  \Phi_{d,+} (x,y) =  \Phi_{d,+} (y,x) \quad\text{for  \ $(x,y)\in \Z^d_+\times  \Z^d_+ $}   $$
and 
$$  \Phi_{d,+} (x,y) >0 \quad\text{for  \ $x,y\in \Z^d_+  $}.    $$

 $(ii)$  When $d\geq3$,   there exists $c_1>1$ such that  for $x,y\in \Z^d_+$ 
  \begin{align}\label{fund-0-hf0}
 \frac1{ c_1}  (1+|x-y|)^{2-d} \min\{1,\frac{x_1y_1}{1+|x-y|^2}\}  \leq \Phi_{d,+} (x,y)\leq c_1  (1+|x-y|)^{2-d}\min\{1,\frac{x_1y_1}{1+|x-y|^2}\}. 
 \end{align}

      \end{proposition} 
\noindent{\bf Proof of Proposition \ref{lm fund-1}. }
Recall that  when $d\geq 3$,  from \cite{LL,EG},    the fundamental solution  $\Phi_{d}$ of $-\Delta$ in $\Z^d$ has   the following asymptotic behaviors:
  \begin{equation}\label{fund-00}
\lim_{|x|_{_Q}\to+\infty}   \Phi_{d} (x,y)|x-y|^{d-2}=\varpi_d
  \end{equation}
  and
   \begin{equation}\label{fund-0}
   0<   \Phi_{d} (x,y) \leq  c_1(1+|x-y|)^{2-d}\quad\text{ in  \ $\Z^d $},
     \end{equation}
where    $ \varpi_d >0$.   Moreover, by \cite[Theorem 2]{U} (also see \cite[Theorem 1]{KS}),
%the fundamental solution in $\Z^d$
$\Phi_{d}$ has the following asymptotic behavior at infinity:
% We conclude that  for $d\geq 3$,
\begin{align}\label{asad-1}
\Phi_{d}(x)=\varpi_d |x|^{2-d} +O(|x|^{-d}) \quad {\rm as}\ |x| \to+\infty.
  \end{align}
While for $d=2$,   the fundamental solutions in the whole space $\Z^2$ are different.   It was proven by \cite[Theorem 7.3]{Ke} that
\begin{equation}\label{eq 2.1}
\left\{%\arraycolsep=1pt
\begin{array}{lll}
-\Delta u = \delta_0 \quad
   &{\rm in}\ \  \Z^2, \\[2mm]
 \phantom{ \   }
  u(0)=0
 \end{array}
 \right.
\end{equation}
  has a unique nonpositive solution  $\Phi_2$   satisfying
 \begin{equation}\label{fund-1-d=2}
  \Phi_2(x)= -\frac1{2\pi}\ln |x|-\frac{\gamma_0}{2}+O(|x|^{-1})\quad {\rm as}\ \ |x|\to\infty,
  \end{equation}
 where $\gamma_0=\frac{1}{\pi}(\gamma_E+\frac12\ln 2)$ with the Euler constant $\gamma_E.$
  \smallskip

\noindent {\it Uniqueness. }
%Let $ w_{1}, w_{2}$ be two solutions of  (\ref{eq 2.1-F}),
%then letting $w=w_1-w_2$, it is a solution of
%$$
%\left\{%\arraycolsep=1pt
%\begin{array}{lll}
%  -\Delta w   =0  \quad\ \,
% {\rm in}\ \  \Z^d_+ ,
% \\[2mm]
% \phantom{  -- }
%w=0\quad\ \,  {\rm on} \ \  \partial \Z^d_+,
% \\[2mm]
% \phantom{    }
% \displaystyle \lim_{\Z^d_+, \ |x|\to+\infty} w(x)=0.
% \end{array}
% \right.
%$$
%By the maximum principle, we get $w=0$ in $\Z^d_+$.
The uniqueness follows by the the maximum principle. \smallskip

 \noindent {\it Existence and properties.  }    For $d\geq 2$,  let
 %it is known that
\begin{equation}\label{eq 2.1-hf}
\Phi_{d,+}(x,y):=\left\{%\arraycolsep=1pt
\begin{array}{lll}
\Phi_{d}(x-y) -\Phi_{d}(x-y^*)\qquad {\rm for}\ \, x, y\in\Z^d_+, \\[3mm]
0\qquad\qquad\qquad\  {\rm for}\ \, x\in\partial \Z^d_+\ \ {\rm or}\ \, y\in\partial \Z^d_+,
\end{array}
 \right.
\end{equation}
where  $y^*=(-y_1,y')$ with $y'=(y_2,\cdots,y_d)$ and
$\Phi_{d}$ is the fundamental solution of  $-\Delta$ in $\Z^d$. Of course, we can get $\Phi_{d,+}(\cdot,y)=0$ on $\partial \Z^d_+$.
Note that
\begin{align*}
\Phi_{d,+}(x,y) =\Phi_{d}(x-y) -\Phi_{d}(x-y^*) &=\Phi_{d}(x-y) -\Phi_{d}\big((x-y^*)^*\big)
\\[1mm]&=\Phi_{d}(y-x) -\Phi_{d}(y-x^*)
\\[1mm]&=\Phi_{d,+}(y,x),
   \end{align*}
since $\Phi_{d}(x)=\Phi_{d}(z)$  for $|z|=|x|$.

Finally, (\ref{fund-0-hf0}) follows by \cite[Theorem 8.3]{GHS} with $k=1$.  \hfill$\Box$\medskip

  \begin{lemma}\label{lm 5.2-mm}

  Let $d\geq 3$,  
$$\cA_0=\Big\{(x_1,x')\in\R^d:\, x_1\geq |x'| \Big\}  $$
and
   $g_\tau  \in C(\Z^d)$
 be  a nonnegative function with a parameter $\tau <-1$.  

$(i)$  If there is $\tau\in(-d,-1)$ such that
$$ g_{\tau} (x)\geq  |x|^{\tau} \quad {\rm for\ }\,  (\Z^d\cap \cA_0)\setminus B_{n_0}$$
for some   $n_0>1$,
  then there exists $c>0$ such that
   \begin{equation}\label{eq 5.1-homx1x}
\Phi_{d,+}\ast g_\tau (x)\geq c x_1(1+|x|)^{1+\tau}    \quad {\rm for\ }\,  \Z^d_+ .
 \end{equation}

 $(ii)$
 If  there is $\tau\leq -d$ such that
$$ g_{\tau} (x)\geq   |x|^{\tau}   \quad {\rm for\ }\,  (\Z^d\cap \cA_0)\setminus B_{n_0}, $$
  then there exists $c>0$ such that
   \begin{equation}\label{eq 5.2-homx1x}
\Phi_{d,+}\ast g_\tau(x)\geq c x_1(1+|x|)^{1-d}   \quad {\rm for\ }\,  \Z^d_+ .
 \end{equation}

   $(iii)$
 If  there is $\tau<-1$ such that
$$ g_{\tau} (x)\leq  (1+ |x|)^{\tau}   \quad {\rm for\ }\,  \Z^d_+,  $$
  then there exists $c>0$ such that for $x\in \Z^d_+$
   \begin{equation}\label{eq 5.2-homx1xupp}
 \Phi_{d,+}\ast g_\tau (x)\leq  \left\{%\arraycolsep=1pt
\begin{array}{lll}
  c x_1  (1+|x|)^{1+\tau}\ln (e+|x|)      \quad 
 &{\rm if}\ \   \tau\in(-1-d,-1),   \\[2mm]
 \phantom{     }
 \displaystyle    c x_1 |x|^{-d} &{\rm if}\ \   \tau<-1-d.  
 \end{array}
 \right.   
 \end{equation}

\end{lemma}
\noindent{\bf Proof. } For given $x\in \Z^d_+$ and any $y\in (\Z^d\cap \cA_0)$ and $|y|>4|x|$,
 we see that 
 $$\sqrt{2}y_1\geq |y|\geq y_1,\quad |x|\geq x_1 $$
 and
$$1\leq x_1\leq |x|\leq \frac14 |y|\leq \frac{\sqrt{2}}4y_1, $$
then
\begin{align*}
\frac{x_1y_1}{1+|x-y|^2} \leq  \frac{\frac{\sqrt{2}}4y_1^2 }{1+\frac{9}{16}|y|^2} \leq  \frac{\frac{\sqrt{2}}4y_1^2 }{1+\frac{9}{16} y_1^2} <1 
\end{align*}
and by (\ref{fund-0-hf0})
there exists $c>0$ such that 
$$\Phi_{d,+}(x,y)\geq c x_1 y_1 |x-y|^{-d}\geq c' x_1   |y|^{1-d}.   $$

 $(i)$ For $x\in \Z^d_+$,  we have that
 \begin{align*}
\Phi_{d,+}\ast g_\tau (x) =\int_{\Z^d_+}\Phi_{d,+}(x,y) g_\tau(y)dy
&\geq c   x_1 \int_{\cA_0\setminus B_{4|x|}}y_1  |x-y|^{-d} |y|^{\tau}dy
\\[1mm]&\geq c   x_1 |x|^{1-\tau} \int_{\cA_0 \setminus B_{4}} z_1 |e_1-z|^{-d} |z|^{\tau}dz,
  \end{align*}
  which implies  (\ref{eq 5.1-homx1x}), where we recall $e_1=(1,0,\cdots,0)\in\R^d$.  \smallskip

  $(ii)$    For $x\in \Z^d_+$,  we have that
 \begin{align*}
\Phi_{d,+}\ast g_\tau (x) 
&\geq c   x_1 \int_{\cA_0\cap \big(B_{8|x|} \setminus B_{4|x|}\big)}  |x-y|^{1-d} |y|^{\tau}  dy
\\[1mm]&\geq c'   x_1 |x|^{1-d}  \int_{\cA_0 (B_8\setminus B_4) }  |z|^{\tau}  dz 
\\[1mm]&\geq c''  x_1 |x|^{1-d}. 
  \end{align*}
  Thus, together with $\bar  v_{\tau}>0$ in $\Z^d_+$,   we obtain  (\ref{eq 5.2-homx1x}).    \smallskip

 $(iii)$ For $x\in \Z^d_+$,  it follows by the second inequality of (\ref{fund-0-hf0}) that 
 \begin{align*}
 \Phi_{d,+}\ast g_\tau (x) & 
 \leq c  x_1 \int_{\R^d}  |y_1| (1+|x-y|) ^{-d} (1+|y|)^{\tau} dy.
  \end{align*}
  
  For $|x|\leq 8$, we can get
  $\bar  v_{\tau}(x)$ is bounded.

  Now we set $|x|>8$.  Let 
  $$\cK_\tau(x,y)= |y_1| (1+|x-y|) ^{-d} (1+|y|)^{\tau},  $$
 then  by  direct computations, for $\tau<-1$
  \begin{align*}
  \int_{\R^d \setminus B_{2|x|}}  \cK_\tau(x,y) dy &\leq c    \int_{\R^d \setminus B_{2|x|}}  |x-y|^{-d} |y|^{1+\tau}dy
\\[1mm]&= c     |x|^{1+\tau}  \int_{\R^d\setminus B_{2}}  |e_x-z|^{-d} |z|^{1+\tau}  dz
 \leq  c'  |x|^{1+\tau}  ,
 \end{align*}
  \begin{align*}
  \int_{  B_{\frac12|x|}}    \cK_\tau(x,y)  dy
 &\leq c    |x|^{-d}   \int_{  B_{\frac12|x|}} |y_1|  (1+|y|)^{\tau}  dy
 \\&\leq c'   |x|^{-d}    \int_0^{\frac12|x|}  (1+r)^{\tau+d}   dr
 \\&\leq \left\{%\arraycolsep=1pt
\begin{array}{lll}
  c''  |x|^{1+\tau}   \quad\ \,
 &{\rm if}\ \   \tau\in(-1-d,-1), \\[2mm]
 \phantom{     }
 \displaystyle    c''  |x|^{-d} &{\rm if}\ \   \tau\leq -1-d
 \end{array}
 \right.  
 \end{align*}
  and
   \begin{align*}
  \int_{  B_{2|x|}\setminus B_{\frac12|x|}}  \cK_\tau(x,y) dy &\leq c    |x|^{1+\tau}     \int_ {B_{2|x|}\setminus B_{\frac12|x|}}   (1+|x-y|)^{-d}    dy
\\[1mm]&\leq c    |x|^{1+\tau}       \int_ {B_{4|x|}(x)} (1+|x-y|)^{-d}    dy
\\[1mm]&\leq c     |x|^{1+\tau}\ln |x|.      
 \end{align*}
As a consequence, we derive   (\ref{eq 5.2-homx1xupp}).  \hfill$\Box$\medskip

 To show the non-existence, we need the following auxiliary lemmas.

  \begin{lemma}\label{pr 5.2w}
  Let $d\ge 2$ and nonnegative function $f\in C(\Z^d_+)$  verify that either 
\begin{equation}\label{con 5.1-1w}
\lim_{n\to+\infty}\int_{ \cA_0 \setminus   B_n(0)   }f(x) (1+ |x|)^{1-d} dx=+\infty
\end{equation}
or for some $\tau\geq -1$, $c>0$ and $n>1$
\begin{equation}\label{con 5.1-1w1}
 f(x) \geq c  |x|^{\tau}\quad {\rm for}\ x\in(\cA_0\cap \Z^d)\setminus  B_n(0). 
\end{equation}

Then the homogeneous problem
\begin{equation}\label{eq 1.1 EH1-5}
\left\{%\arraycolsep=1pt
\begin{array}{lll}
-\Delta u    \geq  f   \quad
   &{\rm in}\ \  \Z^d_+  , \\[2mm]
 \phantom{ --  }
u\geq 0 \quad &{\rm   in}\ \    \Z^d_+
 \end{array}
 \right.
\end{equation}
 has no solutions.
 \end{lemma}
\noindent{\bf Proof. } By contradiction, we assume that  $u_0$ is a nonnegative solution of  (\ref{eq 1.1 EH1-5}), then
strong maximum principle implies that $u_0>0$ in $\Z^d_+$.

 Let  $v_{n,f}$ be the unqiue positive solution of
 \begin{equation}\label{eq 5.1-fn}
\left\{%\arraycolsep=1pt
\begin{array}{lll}
-\Delta u   = f_n   \quad
   {\rm in}\ \ \Z^d_+   , \\[2mm]
 \phantom{  --}
 u=0\quad
   {\rm on}\ \ \partial \Z^d_+,
 \\[2mm]
 \phantom{   }
\displaystyle \lim_{x\in\Z^d_+, |x|\to+\infty}u(x)= 0 ,
 \end{array}
 \right.
\end{equation}
where  $f_n=f\chi_{ B_n(0)  }$.
By the comparison principle, we have that
$$0\leq v_{n,f}\leq u_0\quad{\rm in}\ \, \Z^d_+$$
and
$$v_{n,f}(x)=\int_{ \Z^d_+}\Phi_{d,+}(x,z)f_n(z)dz,\quad\forall\,  x\in\Z^d_+.  $$
It follows by (\ref{con 5.1-1w}) or (\ref{con 5.1-1w1})  that  there exists $c>0$ such that for $n>4$
\begin{align*}
u_0(e_1)  \geq v_{n,f}(e_1)&=\int_{\Z^d }\Phi_{d,+}(e_1,z)f_n(z) dz
 \\& \geq c\int_{\cA_0\cap( B_{n}\setminus B_{4}) }   |z|^{1-d}f(z) dz
\to+\infty\quad{\rm as}\ \ n\to+\infty,
\end{align*}
which is impossible. The non-existence part follows.  \hfill$\Box$\medskip

\begin{lemma}\label{lm 5.1-ew}
Let $d\geq 2$ and $\alpha<d$, $q\in(0,\frac{d+1-\alpha}{d-1})$
 and $\{\tau_j\}_j$ be a sequence defined by
 $$\tau_0=1-d<0,\qquad  \tau_{j+1}=  \tau_jq -\alpha+2,\quad j\in\N_+,$$
 where  $\N_+$ be the set of positive integers.

Then $j\in\N\to \tau_j$ is  strictly increasing and   for any $\bar \tau>\tau_0$ if $q\geq1$ or
for any $\bar \tau\in(\tau_0, \frac{2-\alpha}{1-q})$ if $q
\in(0,1)$,   there exists $j_0\in\N $  such that
\begin{equation}\label{2.3}
 \tau_{j_0}\ge \bar \tau\quad {\rm and}\quad \tau_{j_0-1}<\bar \tau.
\end{equation}
\end{lemma}
The proof is similarly to Lemma \ref{lm 2.1-ew} and we omit it. \medskip

 \noindent{\bf Proof of Theorem \ref{teo 1-hf}. }
 {\it Part $(i)$:  Existence in the Sobolev super critical case: }
  We do the zero extension of   $\Phi_{d,+}$ in $(\Z^d\times\Z^d)\setminus (\Z^d_+\times\Z^d_+)$ and we still denote it by $\Phi_{d,+}$, even extension for $Q$ as following
 $$Q(x_1,x' )= 0\quad {\rm for}\ x_1\leq 0,  \,  x'\in\Z^{d-1}.  $$
Then the original equation  (\ref{eq 1.1-half}) turns to the following integral equation
  $$u=\Phi_{d,+}\ast (Q|u|^{p-2}u)\quad {\rm in}\ \, \Z^d. $$
  Now let
\begin{align*}
v=Q^{\frac1{p'}}|u|^{p-2}u\quad {\rm in}\ \, \Z^d,
\end{align*}
then
 \begin{align} \label{int eq-1-hf}
|v|^{p'-2}v=Q^{\frac1{p}} \Phi_{d,+}\ast (Q^{\frac1{p}}v) \quad {\rm in}\ \, \Z^d.
\end{align}

We employ Theorem \ref{teo 3.1} with  $\beta=1$,  $d\geq 3$
and replace $\Phi_{d,1}$ by $\Phi_{d,+}$
  to obtain that (\ref{int eq-1-hf}) has a nonnegative nontrivial solution $v$.
  Here $(\bA_{\alpha,\beta,1})$ and $(\bA_{\alpha,\beta,2})$ become $(\bA_1)$, $(\bA_2)$ respectively.

 Now we let  $$u= \Phi_{d,+}\ast (Q^{\frac{1}{p}}v)\quad {\rm in}\ \Z^d, $$
 then
 $$u=\Phi_{d,+}\ast (Q|u|^{p-2}u)\quad {\rm in}\ \, \Z^d $$
 and
 $$\int_{\Z^d_+} u^pQ dx=\int_{\Z^d_+} \big(Q^{\frac1{p'}}u^{p-1}\big)^{p'} dx=\int_{\Z^d} |v|^{p'}dx<+\infty.  $$
So
  \begin{equation}\label{eq 1.1-half-ex}
 \left\{%\arraycolsep=1pt
\begin{array}{lll}
 -\Delta  u=Q |u|^{p-2}u  \quad
 &{\rm in}\ \  \Z^d_+ , \\[2mm]
 \phantom{ \ \   }
 \displaystyle u=0    &{\rm on}\ \,   \Z^d\setminus \Z^d_+
 \end{array}
 \right.
\end{equation}
 and  $u$ is a nonnegative nontrivial solution of (\ref{eq 1.1-quad}).
 It follows by the strong maximum principle  that $u>0$ in $\Z^d_+$. \smallskip

 \smallskip

{\it Part $(ii)$: } We first consider the case: $p-1\in (0,1)$ with $\alpha>2$.   Let
$$\bar v_p(x)=   (1+|x|)^{\tau_p}\quad{\rm for}\ \, x\in \Z^d_+ $$
where we set 
$$\tau_p \in\big(-1-d, \, -1  \big) $$
and 
$$ \tau_p>\frac{2(p-1)-\alpha}{2-p}$$
So we need 
$$\frac{2(p-1)-\alpha}{2-p} <-1$$
which holds for $\alpha>2$ and $p\in(1,2)$. 
Note that  $\frac{2(p-1)-\alpha}{2-p} >-1-d$ requires 
$  p-1>\frac{d+1-\alpha}{d-1}.$ Here $\frac{d+1-\alpha}{d-1}>0$ for $\alpha\in(2,d+1)$.
So $\frac{2(p-1)-\alpha}{2-p} \leq -1-d$, we choose $\tau_p=-d$.  As a consequence, our requirement becomes $p\in(1,2)$.

Note that  $\Phi_{d,+}\ast  \bar v_{p}$ in $\Z^d_+$ is well-defined 
  and 
  $$\tau_p >  (p-1)(\tau_p+2) -\alpha. $$
For $t>0$, denote
$$\bar u_t=t \Phi_{d,+}\ast  \bar v_{p} \quad{\rm in}\ \,  \Z^d_+. $$
From Lemma \ref{lm 5.2-mm} part $(iii)$,  we have that
$$  \bar u_t(x)\leq c  tx_1 (1+|x|)^{1+\tau_p} \ln (e+|x|) \quad{\rm for}\ \, x\in \Z^d_+. $$
  Note that for  $x\in\Z^d_+$,
\begin{align*}
Q(x) \bar u_t(x)^{p-1} &\leq   (ct)^{p-1} (1+|x|)^{ (p-1)(\tau_p+2)-\alpha} \big(\ln (e+|x|)\big)^{p-1}
\\[1mm]&\leq C t^{p-1}  (1+|x|)^{ \tau_p }
\\[1mm]& \leq C t^{p-2}   (t\bar v_p)
  = C t^{p-2}( -\Delta \bar u_t),
\end{align*}
for $t\geq t_1$, where
$t_1\geq 1$ is taken such that $C t_1^{p-2}\leq 1. $
It follows by Theorem \ref{teo 3.3} that  problem (\ref{eq 1.1-half})
has  a unique positive solution $u$ such that for some $c>0$
$$0<u(x)\leq c t_1x_1 (1+|x|)^{1+\tau_p} \ln (e+|x|) \quad {\rm for}\ x\in\Z^d_+. $$

{\it Part $(iii)$. }    By contradiction, suppose that there is a $u_0\in C(\Z^d_+)$, a nonnegative nonzero function verifying    (\ref{eq 1.1-half}).   By the maximum principle, we obtain that $u_0>0\quad {\rm in}\ \, \Z^d_+. $

%The part of  proof is similar to the one of Proposition \ref{teo 2-wh}. For convenience of readers, we show the proof.

From  the comparison principle, there exists $d_0>0$  and $n_0\geq 1$ such that
%\begin{equation}\label{sec 3-1.0-0}
$$
u_0(x)\geq  \frac{u_0(e_1)}{\Phi_{d,+}(e_1, e_1)}\, \Phi_{d,+}(x, e_1)\geq d_0 x_1(1+ |x|)^{-d}\geq d_0  (1+ |x|)^{1-d}\quad {\rm for}\ \, x\in  \cA_0\cap \Z^d.
$$%\end{equation}
Let $\tau_0=1-d<-1$ satisfy that
\begin{equation}\label{sec 3-5.1.0}
-\Delta u_0(x)   \geq  d_0^q   x_1^q |x|^{-\alpha-qd} \geq d_0^q    |x|^{\tau_1-2},\quad\forall\, x\in ( \cA_0\cap \Z^d)\setminus B_{n_0},
\end{equation}
where   $$q=p-1,   \quad  \tau_1=-q(d-1)-\alpha+2.$$
 Thus, for $q\in(0,  \frac{d+1-\alpha}{d-1})$, it holds that
$$\tau_1-\tau_0=-q(d-1)-\alpha+2 +(d-1)>0.$$

\smallskip

{\bf Case 1}:  $q\in(0, \frac{1-\alpha}{d-1}]$.
%with $\alpha\in(-\infty,1)$. }
Note that $q(1-d)-\alpha\geq  -1$, then
\begin{equation}\label{sec 3-5.1.0-0}
Q(x)  u_0(x)^{q}\geq  d_0^q (1+  |x|)^{(1-d)q-\alpha},\quad\forall x\in  (\cA_0\cap \Z^d) \setminus B_{n_0}
\end{equation}
and a contradiction follows by Lemma  \ref{pr 5.2w}  with $f(x)=d_0^q (1+|x|)^{q(1-d) -\alpha}$ for $x\in  \cA_0\cap \Z^d$.   \smallskip

{\bf Case 2}:     $q\in\big(\frac{1-\alpha}{d-1},   \frac{d+1-\alpha}{d-1}\big)\subset  (0,+\infty)$.
%with $\alpha\in(-\infty,d)$. }
 By Lemma \ref{lm 5.2-mm} part $(i)$, there exists $d_1>0$   such that
$$u_0(x)\geq d_1(1+|x|)^{\tau_1},\quad\forall\, x\in    \cA_0\cap \Z^d, $$
where
$\tau_1:=-q(d-1)-\alpha+2\in(1-d,-1)$.

Recall that
$$\tau_{j+1}:=q\tau_{j} -\alpha+2,\quad\forall\,  j\in\N_+,$$
 which is an increasing sequence.
 
 Assume that 
 $$u_0(x)\geq d_{j}(1+ |x|)^{\tau_{j}}\quad {\rm for}\ \,  x\in \cA_0\cap \Z^d, $$
then 
$$Q(x)  u_0(x)^{q}\geq  d_0^q (1+  |x|)^{\tau_j q-\alpha},\quad\forall x\in  (\cA_0\cap \Z^d) \setminus B_{n_0}. $$

 If $\tau_{j+1}=\tau_jq -\alpha+2\in  (-1,1-d)$,
it follows by Lemma \ref{lm 5.2-mm} part $(i)$ that there exist integer  $d_{j}>0$  such that
$$u_0(x)\geq d_{j+1}(1+ |x|)^{\tau_{j+1}}\quad {\rm in}\ \,  \cA_0\cap \Z^d . $$
   If $q\tau_{j+1}-\alpha \geq   -1$, we are done by Lemma  \ref{pr 5.2w}. If not, we continue the above process. \smallskip

% Now we claim that the iteration must be stopped after a finite number of times. It infers by
According to Lemma \ref{lm 5.1-ew},
% that $j\mapsto \tau_j$ is strictly increasing thanks to $ 0<q<\frac{d-\alpha}{d-1}. $
we know that for $q\in[1,+\infty)\cap \big(\frac{1-\alpha}{d-1},   \frac{d+1-\alpha}{d-1}\big)$, $\tau_j\to+\infty$, then there exists $j_0\in\N$ such that
$q\tau_{j_0+1} \geq  -1$ and a contradiction could be derived for
$1\leq q<\frac{d}{d-1}$.
While for $q\in(0,1)\cap \big(\frac{1-\alpha}{d-1},   \frac{d+1-\alpha}{d-1}\big)$,  $\tau_j\to \tilde \tau_q:=\frac{2-\alpha}{1-q}>0$ as $j\to+\infty$,  then there exists $j_0\in\N$ such that
 $q\tau_{j_0} -\alpha \leq  -1$ and  $q\tau_{j_0+1}-\alpha  \geq  -1$. Thus we get a contradiction and we are done. 
 \hfill$\Box$\medskip

 \subsection{ In quadrant  $ \Z^d_{*}$ }

To prove Theorem \ref{teo 1-qd},   using  the zero extension technique and  Theorem \ref{teo 3.1}, we only need to show the existence of the fundamental solution $\Phi_{d,*}$  of $-\Delta$  in the quadrant space under the zero Dirichlet boundary condition, as well as the estimates of $\Phi_{d,*}$ at infinity,   i.e.
 \begin{equation}\label{eq 2.1-E}
\left\{%\arraycolsep=1pt
\begin{array}{lll}
  -\Delta u   = \delta_y \quad\ \,
  {\rm in}\ \  \Z^d_{*} ,
 \\[2mm]
 \phantom{  -- }
u=0\quad\ \  {\rm on} \ \  \partial \Z^d_{*},
 \\[2mm]
 \phantom{    }
 \displaystyle \lim_{z\in \Z^d_{*}, |x| \to+\infty} u(x)=0,
 \end{array}
 \right.
\end{equation}
where   $y\in \Z^d_{*}$ and $\delta_y$ is the Dirac mass at $y$.
\begin{proposition}\label{lm fund-2}
Let $d\geq 3$ and $y\in \Z^d_{*}$,  then  (\ref{eq 2.1-E})  has a solution $\Phi_{d,*}$. Moreover,  for  \ $(x,y)\in \Z^d_{*}\times \Z^d_{*} $ 
$$ \Phi_{d,*} (x,y)= \Phi_{d,*} (y,x)>0 $$
    and for some $c>1$ 
   \begin{align}
&\qquad\frac1c (1+|x-y|)^{2-d}\min\big\{1,\frac{x_1y_1}{1+|x-y|^2}\big\} \min\big\{1,\frac{x_2y_2}{1+|x-y|^2}\big\}   
 \nonumber
\\[1mm]& \leq   \Phi_{d,*} (x,y) 
 \leq c (1+|x-y|)^{2-d}\min\big\{1,\frac{x_1y_1}{1+|x-y|^2}\big\} \min\big\{1,\frac{x_2y_2}{1+|x-y|^2}\big\}. \label{fund-1-qd}
   \end{align}

\end{proposition}

 \noindent{\bf Proof. }  {\it Uniqueness. } The uniqueness comes from the maximum principle directly. \smallskip
 %Let $ w_{1}, w_{2}$ be two solutions of  (\ref{eq 2.1-F}),
%then letting $w=w_1-w_2$,
%it is a solution of
%$$
%\left\{%\arraycolsep=1pt
%\begin{array}{lll}
 % -\Delta w   =0  \quad\ \,
% {\rm in}\ \  \Z^d_{*} ,
 %\\[2mm]
 %\phantom{  -- }
%w=0\quad\ \,  {\rm on} \ \  \partial \Z^d_{*},
 %\\[2mm]
 %\phantom{    }
 %\displaystyle \lim_{x\in \Z^d_{*}, |x|\to+\infty} w(x)=0
 %\end{array}
 %\right.
%$$
%and the maximum principle implies the uniqueness directly. \smallskip

  {\it Existence and properties.  }
Denote
\begin{equation}\label{eq 2.1-qd}
\Phi_{d,*}(x,y)=
\left\{%\arraycolsep=1pt
\begin{array}{lll}
 \Phi_{d}(x-y) -\Phi_{d}(x-y^*)  -\Phi_{d}(x-y^\#)+\Phi_{d}(x-(y^*)^\#)  \quad {\rm for}\ \, x\in \Z^d_{*}, \\[3.5mm]
  0\quad \qquad {\rm for}\ \, x\in\partial \Z^d_{*}\ \ {\rm or}\ \, y\in\partial \Z^d_{*},
   \end{array}
 \right.
\end{equation}
where $y^\#=(y_1,-y_2)$ if $d=2$,  $y^\#=(y_1,-y_2,y'')$ with  $y''=(y_3,\cdots,y_d)$ if $d\geq 3$.  Of course, we can get $\Phi_{d,*}(\cdot,y)=0$ on $\partial \Z^d_*$
and $-\Delta_x\Phi_{d,*}=\delta_y$ for $y\in \Z^d_{*}$.
Direct computation shows that $\Phi_{d,*}$ is the fundamental solution of $-\Delta$ in $ \Z^d_{*}$, i.e.  it verifies (\ref{eq 2.1-E}).

Note that for $x,y\in \Z^d_{*}$,
\begin{align*}
\Phi_{d,*}(x,y)&= \Phi_{d}(x-y) -\Phi_{d}(x-y^*)  -\Phi_{d}(x-y^\#)+\Phi_{d}(x-(y^*)^\#) 
\\[1mm]&= \Phi_{d}(y-x) -\Phi_{d}(x^*-y)  -\Phi_{d}(x^\#-y)+\Phi_{d}((x^*)^\#-y) 
\\[1mm]&=\Phi_{d,*}(y,x)
   \end{align*}
by the fact that $\Phi_{d}(x)=\Phi_{d}(z)$  for $|z|=|x|$.

 When $d\geq 2$,  since $\Phi_{d,+}$ decays at infinity, we have that
 $$|\Phi_{d,*}(x,y)|\leq \Phi_{d,+}(x,y),\quad {\rm for}\ (x,y)\in \Z^d_{*}\times \Z^d_{*},$$
then by comparison principle, we have that
$\Phi_{d,*}$ is positive in $ \Z^d_{*}$.

 \smallskip

 {\it Now we show (\ref{fund-1-qd}).}  It follows by \cite[Theorem 8.3]{GHS} with $k=2$. 
 \hfill$\Box$\medskip

 To show the existence, we need the following auxiliary lemmas.

  \begin{lemma}\label{lm 6.2-mm}

  Let $d\geq 3$, 
  $$\cA_1=\Big\{(x_1,x_2,x'')\in\R^d:\,     x_1x_2\geq \frac14(x_1^2+x_2^2) > \frac14 \sum_{j=3}^{d}x_j^2 \Big\} $$
  and
 $g_\tau \in C(\Z^d_*)$   be a nonnegative function.

$(i)$  If there is $\tau\in(-d,-2)$ such that
$$ g_{\tau} (x)\geq  (1+|x|)^{\tau} \quad {\rm for\ }\, x \in  (\Z^d\cap \cA_1)\setminus B_{n_0}$$
for some   $n_0>1$,
  then there exists $c>0$ such that
   \begin{equation}\label{eq 6.1-homx1x}
 \Phi_{d,*}\ast g_{\tau}(x)\geq c x_1x_2(1+|x|)^{\tau}    \quad {\rm for\ }\, x \in  \Z^d_* .
 \end{equation}

 $(ii)$
 If  $\tau\leq -d$
$$ g_{\tau} (x)\geq   |x|^{\tau}    \quad {\rm for\ }\, x \in  (\Z^d\cap \cA_1)\setminus B_{n_0}$$
for some   $n_0>e$,
  then there exists $c>0$ such that
   \begin{equation}\label{eq 6.2-homx1x}
 \Phi_{d,*}\ast g_{\tau} (x)\geq c x_1x_2 (1+|x|)^{-d}     \quad {\rm for\ }\,  x \in  \Z^d_*.
 \end{equation}

   $(iii)$
   If  there is $\tau<-2$ such that
$$ g_{\tau} (x)\leq  (1+ |x|)^{\tau}   \quad {\rm for\ }\,  \Z^d_+,  $$
  then there exists $c>0$ such that for $x\in \Z^d_*$
   \begin{equation}\label{eq 6.2-homx1xupp}
 \Phi_{d,*}\ast g_\tau (x)\leq  \left\{%\arraycolsep=1pt
\begin{array}{lll}
  c x_1x_2   |x|^{ \tau}      \quad 
 &{\rm if}\ \   \tau\in(-2-d,-2),   \\[2mm]
 c x_1x_2   |x|^{ \tau}\ln (e+|x|)      \quad 
 &{\rm if}\ \   \tau=-2-d,   \\[2mm]
 \phantom{     }
 \displaystyle    c x_1x_2 |x|^{-d-2} &{\rm if}\ \   \tau<-2-d.  
 \end{array}
 \right.   
 \end{equation}

 \end{lemma}
\noindent{\bf Proof. }    For given $x\in \Z^d_+$ and any $y\in (\Z^d\cap \cA_1)$ and $|y|>4|x|$,
 we see that 
 $$ 8y_1y_2\geq  2(y_1^2+y_2^2)\geq |y|^2\geq y_1^2+y_2^2\geq 2y_1y_2,\quad |x|\geq x_1x_2 $$
 and
$$1\leq x_1x_2 \leq \frac12  |x|^2 \leq \frac1{32} |y|^2 \leq \frac14 y_1y_2, $$
then
\begin{align*}
\frac{x_1x_2y_1y_2}{1+|x-y|^4} \leq   \frac{\frac14(y_1y_2)^2 }{1+(\frac{3}{4})^4 |y|^4} \leq    \frac{ (y_1y_2)^2 }{1+ 3^4 (y_1y_2)^2} <1 
\end{align*}
and by (\ref{fund-1-qd})
there exists $c>0$ such that 
$$\Phi_{d,*}(x,y)\geq c x_1x_2 y_1y_2 |x-y|^{-2-d}\geq c' x_1x_2   |y|^{-d}.   $$

 $(i)$ For $x\in \Z^d_*$,  we have that
 \begin{align*}
\Phi_{d,*}\ast g_\tau (x) =\int_{\Z^d_*}\Phi_{d,*}(x,y) g_\tau(y)dy
&\geq c   x_1x_2 \int_{\cA_1\setminus B_{4|x|}}  |y|^{-d} |y|^{\tau}dy
 \geq c   x_1x_2 |x|^{\tau}, 
  \end{align*}
  which implies  (\ref{eq 6.1-homx1x}).  \smallskip

  $(ii)$    For $x\in \Z^d_*$,  we have that
 \begin{align*}
\Phi_{d,*}\ast g_\tau (x) 
&\geq c   x_1x_2 \int_{\cA_1\cap \big(B_{8|x|} \setminus B_{4|x|}\big)}  |y|^{-d} |y|^{\tau}  dy
\\[1mm]&\geq c'   x_1x_2  |x|^{-d}  \int_{ B_8\setminus B_4  }  |z|^{\tau}  dz 
\\[1mm]&\geq c''  x_1x_2 |x|^{-d}. 
  \end{align*}
  Thus, together with $\bar  v_{\tau}>0$ in $\Z^d_*$,   we obtain  (\ref{eq 5.2-homx1x}).    \smallskip

 $(iii)$ For $x\in \Z^d_*$,  it follows by the second inequality of (\ref{fund-1-qd}) that 
 \begin{align*}
 \Phi_{d,*}\ast g_\tau (x) & 
 \leq c  x_1x_2  \int_{\R^d}  |y_1||y_2| (1+|x-y|) ^{-d-2} (1+|y|)^{\tau} dy.
  \end{align*}
  
  For $|x|\leq 9$, we can get
  $\bar  v_{\tau}(x)$ is bounded.

  Now we set $|x|>9$.  Reset
  $$\cK_\tau(x,y)= |y_1| |y_2|(1+|x-y|) ^{-d-2} (1+|y|)^{\tau},  $$
 then, taking 
 $$t=|x|-|x|^{\frac{\tau}{\tau+2}} \in(\frac12 |x|,|x|)$$ 
 for $\tau<-2$
  \begin{align*}
  \int_{\R^d \setminus B_{|x|+t }}  \cK_\tau(x,y) dy \leq      \int_{\R^d \setminus B_{|x|+t}}  |x-y|^{-d-2} |y|^{2+\tau}dy
 &\leq      (|x|+t)^{2+\tau}   \int_{\R^d \setminus B_{t}}  |z|^{-d-2}   dz
\\[1mm]&\leq   c  |x|^{2+\tau}t^{-2}
\\[1mm]&\leq   c  |x|^{\tau} (\ln |x|)^{2},  
 \end{align*}
  \begin{align*}
  \int_{  B_{|x|-t}}    \cK_\tau(x,y)  dy
 &\leq      t^{-d-2}   \int_{  B_{|x|-t} } |y_1| |y_2|  (1+|y|)^{\tau}  dy
 \\[1mm]&\leq c    |x|^{-d-2}    \int_0^{{|x|-t}}  (1+r)^{\tau+d+1}   dr
 \\[1mm] &\leq \left\{%\arraycolsep=1pt
\begin{array}{lll}
  c'    |x| ^{\tau-\frac{2d}{\tau+2}-2}    \quad\ \,
 &{\rm if}\ \   \tau\in(-2-d,-2), \\[2mm]
 \phantom{     }
 \displaystyle    c''  |x|^{-d-2}  \ln(|x|)  &{\rm if}\ \   \tau= -2-d, \\[2mm]
 \phantom{     }
 \displaystyle    c''  |x| ^{-d-2}   &{\rm if}\ \   \tau< -2-d, 
 \end{array}
 \right.  
 \end{align*}
 with $-\frac{2d}{\tau+2}-2<0$ for $\tau\in(-2-d,-2)$, 
  \begin{align*}
  \int_{  B_{t}(x)}  \cK_\tau(x,y) dy &\leq c (|x|- t)^{2+\tau}      \int_{  B_{t}(x)}    (1+|x-y|)^{-2-d}    dy
 \leq c  |x|^{\tau}     \int_{ \R^d}    (1+|y|)^{-2-d}    dy
 \leq c'  |x|^{\tau}        
 \end{align*}
 and
 \begin{align*}
  \int_{  (B_{|x|+t}\setminus B_{|x|-t})\setminus B_{t}(x)}  \cK_\tau(x,y) dy &\leq c   ( |x|-t)^{2+\tau}  t^ {-2-d}   | (B_{|x|+t}\setminus B_{|x|-t})\setminus B_{t}(x) |
 \leq c'    |x|^{\tau-2}.
 \end{align*}
 As a consequence, we derive   (\ref{eq 6.2-homx1xupp}).  \hfill$\Box$\medskip

  \begin{lemma}\label{pr 6.2w}
  Let $d\ge 2$ and $f\in C(\Z^d_*)$ be a nonnegative function verifying that  
\begin{equation}\label{con 6.1-1w}
\lim_{n\to+\infty}\int_{ \cA_1 \setminus   B_n(0)   }f(x) (1+ |x|)^{-d} dx=+\infty,
\end{equation}
where $\cA_1$ is given in Lemma \ref{lm 6.2-mm}. Then the homogeneous problem
\begin{equation}\label{eq 1.1 EH1-6}
\left\{%\arraycolsep=1pt
\begin{array}{lll}
-\Delta u    \geq  f   \quad
   &{\rm in}\ \  \Z^d_* , \\[2mm]
 \phantom{ --  }
u\geq 0 \quad &{\rm   in}\ \    \Z^d_*
 \end{array}
 \right.
\end{equation}
 has no solutions.
 \end{lemma}
\noindent{\bf Proof. } The proof is by contradiction and very similar to that of Lemma \ref{pr 5.2w}, where we replace $e_1$ by $e_{11} = (1,1, \cdots, 0) \in
\Z^d$, $\cA_0$ by $\cA_1$ and $|z|^{1-d}$ is replaced by $|z|^{-2-d}$.  \hfill$\Box$\medskip

%Let  $u_0$ be a nonnegative solution of  (\ref{eq 1.1 EH1}), then we can assume that $u_0>0$ in $\Z^d_*$ by the strong maximum principle.

 %Let  $v_{n,f}$ be the minimal positive solution of
 %\begin{equation}\label{eq 5.1-fn}
%\left\{%\arraycolsep=1pt
%\begin{array}{lll}
%-\Delta u   = f_n   \quad
%   {\rm in}\ \ \Z^d_*   , \\[2mm]
% \phantom{  --}
% u=0\quad
%   {\rm on}\ \ \partial \Z^d_*,
% \\[2mm]
% \phantom{   }
%\displaystyle \lim_{x\in\Z^d_*, |x|\to+\infty}u(x)= 0 ,
% \end{array}
% \right.
%\end{equation}
%where  $f_n=f\chi_{ B_n(0)  }$.
%By comparison principle, we have that
%$$0\leq v_{n,f}\leq u_0\quad{\rm in}\ \, \Z^d_*$$
%and
%$$v_{n,f}(x)=\int_{ \Z^d_+}\Phi_{d}(x,z)f_n(z),\quad\forall\,  x\in\Z^d_*.  $$
%There is $c_1>1$ such that
%$$\frac{1}{c_1} x_1|x|^{-d}\leq v_{n,f}(x)\leq c_1 x_1|x|^{-d} \quad {\rm for}\ \, x\in\Z^d_*$$
%and it follows by (\ref{con 2.1-1w}) and the comparison principle that  there exists $C >0$ such that for $n>4$
%\begin{align*}
%u_0(e_{11})  \geq v_{n,f}(e_{11})&=\int_{\Z^d }\Phi_{d,*}(e_{11},z)f_n(z) dz
% \\& \geq c\int_{\cA_1\cap( B_{n}\setminus B_{4}) }   |z|^{-d}f_n(z) dz
%\to+\infty\quad{\rm as}\ \ n\to+\infty,
%\end{align*}
%which is impossible, where $e_{11}=(1,1,0,\cdots,0)\in\Z^d$.  The non-existence follows.  \hfill$\Box$\medskip

\begin{lemma}\label{lm 6.1-ew}
Let $d\geq 2$ and $\alpha<d$, $q\in(0,\frac{d-\alpha}{d})$
 and $\{\tau_j\}_j$ be a sequence defined by
 $$\tau_0=-d<0,\qquad  \tau_{j+1}=  \tau_jq +\alpha,\quad j\in\N_+,$$
 where  $\N_+$ be the set of positive integers.

Then $j\in\N\to \tau_j$ is  strictly increasing and
for any $\bar \tau\in(\tau_0, \frac{-\alpha}{1-q})$,   there exists $j_0\in\N $  such that
\begin{equation}\label{2.3-6}
 \tau_{j_0}\ge \bar \tau\quad {\rm and}\quad \tau_{j_0-1}<\bar \tau.
\end{equation}
\end{lemma}
The proof is similar to Lemma \ref{lm 2.1-ew} and we omit it. \medskip

\noindent{\bf  Proof of Theorem \ref{teo 1-qd}. }  
{\it Part $(i)$:  Existence in the Sobolev super critical case: }
  We do the zero extension of   $\Phi_{d,*}$ in $(\Z^d\times\Z^d)\setminus (\Z^d_*\times\Z^d_*)$ and we still denote it by $\Phi_{d,*}$, even extension for $Q$ as following
 $$Q(x)= 0\quad {\rm for}\ x\in\Z^d\setminus \Z^d_*.  $$
Then the original equation  (\ref{eq 1.1-quad}) turns to the following integral equation
  $$u=\Phi_{d,*}\ast (Q|u|^{p-2}u)\quad {\rm in}\ \, \Z^d. $$
  Now let
\begin{align*}
v=Q^{\frac1{p'}}|u|^{p-2}u\quad {\rm in}\ \, \Z^d,
\end{align*}
then
 \begin{align} \label{int eq-1-hf}
|v|^{p'-2}v=Q^{\frac1{p}} \Phi_{d,*}\ast (Q^{\frac1{p}}v) \quad {\rm in}\ \, \Z^d.
\end{align}

We employ Theorem \ref{teo 3.1} with  $\beta=1$,  $d\geq 3$
and replace $\Phi_{d,1}$ by $\Phi_{d,*}$
  to obtain that (\ref{int eq-1-hf}) has a nonnegative nontrivial solution $v$.
  Here $(\bA_{\alpha,\beta,1})$ and $(\bA_{\alpha,\beta,2})$ become $(\bA_1)$, $(\bA_2)$ respectively.

 Now we let  $$u= \Phi_{d,*}\ast (Q^{\frac{1}{p}}v)\quad {\rm in}\ \Z^d, $$
 then
 $$u=\Phi_{d,*}\ast (Q|u|^{p-2}u)\quad {\rm in}\ \, \Z^d $$
 and
 $$\int_{\Z^d_*} u^pQ dx=\int_{\Z^d_*} \big(Q^{\frac1{p'}}u^{p-1}\big)^{p'} dx=\int_{\Z^d} |v|^{p'}dx<+\infty.  $$
So
  \begin{equation}\label{eq 1.1-half-ex}
 \left\{%\arraycolsep=1pt
\begin{array}{lll}
 -\Delta  u=Q |u|^{p-2}u  \quad
 &{\rm in}\ \  \Z^d_* , \\[2mm]
 \phantom{ \ \   }
 \displaystyle u=0    &{\rm on}\ \,   \Z^d\setminus \Z^d_*
 \end{array}
 \right.
\end{equation}
 and  $u$ is a solution of (\ref{eq 1.1-quad}).
 It follows by the strong maximum principle  that $u>0$ in $\Z^d_*$. \smallskip

{\it Part $(ii)$: Existence for sublinear case. }  We first consider the case: $p-1\in (0,1)$ with $\alpha>2$.   Let
$$\bar v_p(x)=   (1+|x|)^{\tau_p}\quad{\rm for}\ \, x\in \Z^d_* $$
where we set 
$$\tau_p \in\Big(-2-d, \,  -2  \Big)\quad {\rm and}\quad \tau_p > \frac{2(p-1)-\alpha}{2-p} $$

 Note that for $\alpha>2$ and $p\in(1,2)$, there holds $\frac{2(p-1)-\alpha}{2-p}  <-2$,  
 and if  $p-1<\frac{d+2-\alpha}{d}$, then  $ \frac{2(p-1)-\alpha}{2-p} >-2-d$
 Here $\frac{d+2-\alpha}{d}>0$ for $\alpha\in(2,2+d)$. So when $p-1\geq \frac{d+2-\alpha}{d}$, we choose $\tau_p=-1-d$.  
 So our requirements become  that $p\in(1,2)$ and $\alpha>2$.

For $t>0$, denote
$$\bar u_t=t \Phi_{d,*}\ast  \bar v_{p} \quad{\rm in}\ \,  \Z^d_*. $$
From Lemma \ref{lm 5.2-mm} part $(iii)$,  we have that
$$  \bar u_t(x)\leq c  tx_1x_2 (1+|x|)^{\tau_p}  \quad{\rm for}\ \, x\in \Z^d_*. $$
  Note that for  $x\in\Z^d_*$,
\begin{align*}
Q(x) \bar u_t(x)^{p-1} &\leq   (ct)^{p-1} (1+|x|)^{ (p-1)(\tau_p+2)-\alpha}  
\\[1mm]&\leq C t^{p-1}  (1+|x|)^{ \tau_p }
 \leq C t^{p-2}   (t\bar v_p)
  = C t^{p-2}( -\Delta \bar u_t),
\end{align*}
for $t\geq t_1$, where 
$$\tau_p >  (p-1)(\tau_p+2) -\alpha $$ and 
$t_1\geq 1$ is taken such that $C t_1^{p-2}\leq 1. $
It follows by Theorem \ref{teo 3.3} that  problem (\ref{eq 1.1-quad})
has  a unique positive solution $u$ such that for some $c>0$
$$0<u(x)\leq c t_1x_1x_2 (1+|x|)^{1+\tau_p}  \quad {\rm for}\ x\in\Z^d_*. $$

{\it Part $(iii)$. }  The proof is similar to the one of Proposition \ref{teo 2-wh}. Suppose that, by contradiction, there is a positive solution $u_0$  of  (\ref{eq 1.1-quad}).

 From  the comparison principle, there exists $d_0>0$  and $n_0\geq 1$ such that
%\begin{equation}\label{sec 3-1.0-0}
$$
u_0(x)\geq  \frac{u_0(e_{11})}{\Phi_{d,+}(e_{11}, e_{11})}\, \Phi_{d,*}(x, e_{11})\geq d_0 (1+ |x|)^{-d}\quad {\rm for}\ \, x\in  \cA_1\cap \Z^d.
$$%\end{equation}

Let $\tau_0=-d<0$ satisfy that
\begin{equation}\label{sec 3-5.1.0}
-\Delta u_0(x)   \geq  d_0^q     |x|^{\tau_1-2 },\quad\forall\, x\in (\cA_1\cap \Z^d)\setminus B_{n_0},
\end{equation}
where   $$\tau_1=2-\alpha -qd\quad {\rm and}\quad q=p-1.$$
 Thus, for $q\in(0,  \frac{d+2-\alpha}{d})$, it holds that
$$\tau_1-\tau_0=2-qd-\alpha  +d>0.$$

\smallskip

{\bf Case 1}:  $q\in(0, \frac{-\alpha}{d})$.
% with $\alpha\in(-\infty,0)$. }
Note that $-qd-\alpha >  0$, then
\begin{equation}\label{sec 3-5.1.0-0}
Q(x)  u_0(x)^{q}\geq  d_0^q (1+  |x|)^{-dq-\alpha},\quad\forall x\in (\cA_1\cap \Z^d) \setminus B_{n_0}
\end{equation}
and a contradiction follows by Lemma  \ref{pr 6.2w}  with $f(x)=d_0^q (1+|x|)^{-dq-\alpha}$.   \smallskip

{\bf Case 2}:     $q\in\big(\frac{-\alpha}{d},   \frac{d+2-\alpha}{d}\big)\subset  (0,+\infty)$.
% with $\alpha\in(-\infty,d)$. }
 By Proposition \ref{lm fund-2}, there exists $d_1>0$   such that
$$u_0(x)\geq d_1(1+|x|)^{\tau_1},\quad\forall\, x\in   \Z^d\cap \cA_1, $$
where
$\tau_1:=2-qd-\alpha \in(-2-d,-2)$.

%Recall that
%$$\tau_{j+1}:=q\tau_{j} -\alpha,\quad\forall\,  j\in\N_+,$$
% which is an increasing sequence.

 If for some $j\in\N_+,$ there holds 
 $$\tau_{j+1}=2-\alpha+ \tau_jq  \in  (-2-d,-2),$$
  then
it follows by Proposition \ref{lm fund-2} that there exist integer  $d_{j}>0$  such that
$$u_0(x)\geq d_{j}(1+ |x|)^{\tau_{j+1}}\quad {\rm in}\ \,  \cA_1\cap \Z^d. $$
   If $q\tau_{j+1}-\alpha \geq0$, we are done by Lemma  \ref{pr 6.2w}.\smallskip

% Now we claim that the iteration must stop after a finite number of times. It infers by   Lemma \ref{lm 6.1-ew} that
% $j\mapsto \tau_j$ is strictly increasing thanks to $ 0<q<\frac{d-\alpha}{d}. $

According to Lemma \ref{lm 6.1-ew} and by similar argument, we can get that for $q\in[1,+\infty)\cap \big(\frac{-\alpha}{d},   \frac{d+2-\alpha}{d}\big)$, $\tau_j\to+\infty$, then there exists $j_0\in\N$ such that
$q\tau_{j_0+1} -\alpha \geq  0$ and a contradiction could be derived for
$1\leq q<\frac{d+2-\alpha}{d}$.
For $q\in(0,1)\cap \big(\frac{-\alpha}{d},   \frac{d+2-\alpha}{d}\big)$, then  $\tau_j\to \tilde \tau_q:=\frac{2-\alpha}{1-q}>0$ as $j\to+\infty$,  thus there exists $j_0\in\N$ such that
 $q\tau_{j_0} -\alpha<  0$ and  $q\tau_{j_0+1}-\alpha  \geq  0$. Again we get a contradiction and we are done.  
   \hfill$\Box$

  \bigskip
 
  \bigskip

 {\small
 
 \noindent {\bf  Conflicts of interest:} The authors declare that there is no conflict of interest.\medskip

 \noindent {\bf  Data availability:}  This paper has no associated data.\medskip

 \noindent {\bf Acknowledgements:} We are grateful to Yuhua Sun and Yanyu Qin for their valuable comments and suggestions on this paper. \smallskip

  H. Chen is supported by  NSFC (No.   12361043).
 
 B. Hua is supported by NSFC (no.12371056), and by Shanghai Science and Technology Program (No. 22JC1400100).

 F. Zhou is supported by Science and Technology Commission of Shanghai Municipality (No. 22DZ2229014)
 and also NSFC (No. 12071189).

}

  % \end{document}
  \bigskip\medskip

  \small
  
\noindent   Huyuan Chen:    Center for Mathematics and Interdisciplinary Sciences,
   
	\qquad \qquad \quad   	Fudan University, Shanghai 200433, PR China \smallskip
	
\qquad \qquad \quad  Shanghai Institute for Mathematics and Interdisciplinary Sciences,  

	\qquad \qquad  \quad  	 Shanghai 200433, PR China\smallskip
	
\noindent    Email address: chenhuyuan@simis.cn, chenhuyuan@yeah.net
		 
		 \bigskip

\noindent   Bobo Hua:  School of Mathematical Sciences, LMNS, Fudan University,

\qquad \quad \ \   Shanghai, 200433, P.R. China\smallskip

\noindent    Email address:  bobohua@fudan.edu.cn\bigskip

\noindent Feng Zhou:  CPDE, School of Mathematical Sciences, East China Normal University, 

 \qquad \qquad\!\!  Shanghai 200241, PR China \smallskip
 
  \qquad \qquad\!\!  NYU-ECNU Institute of Mathematical Sciences at NYU-Shanghai, 
  
\qquad \qquad\!\!  Shanghai 200120, PR China\smallskip
  
\noindent      Email address:  fzhou@math.ecnu.edu.cn
  
    \end{document}